\documentclass[10pt]{article}
\usepackage{amsmath,amscd,amsfonts,amssymb,amsthm,enumerate,mathtools,textcomp,xr,
makeidx,stmaryrd,microtype,indentfirst,graphicx}
\usepackage[mathscr]{eucal}

\usepackage[english]{babel}
\usepackage[all,2cell]{xy}
\usepackage[colorlinks,plainpages,urlcolor=blue,linkcolor=reference,citecolor=citation]{hyperref}
\usepackage[abbrev,lite,alphabetic]{amsrefs}
\usepackage{cleveref}
\usepackage{chngcntr}

\usepackage{color}
\definecolor{reference}{rgb}{0.20,0.36,0.74}
\definecolor{citation}{rgb}{0,.40,.80}

\usepackage[left=2cm, top=2cm, right=2cm, bottom=2cm]{geometry}
\usepackage{tocloft}
\UseTwocells
\xyoption{2cell}{\UseTwocells}
\hypersetup{linktocpage}

\makeatletter
\DeclareRobustCommand{\em}{%
  \@nomath\em \if b\expandafter\@car\f@series\@nil
  \normalfont \else \bfseries \fi}
\makeatother

\renewcommand{\phi}{\varphi}

\theoremstyle{definition}
\newtheorem{Theor}{Theorem}[subsection]
\newtheorem*{thm*}{Theorem}
\newtheorem{Lemma}[Theor]{Lemma}
\newtheorem{Prop}[Theor]{Proposition}
\newtheorem*{prop*}{Proposition}
\newtheorem{Rem}[Theor]{Remark}

\newtheorem{Def}[Theor]{Definition}
\newtheorem*{Conv}{Conventions}
\newtheorem*{Conven}{Convention}
\newtheorem{warning}[Theor]{Warning}
\newtheorem{construction}[Theor]{Construction}
\newtheorem{Cor}[Theor]{Corollary}

\newtheorem{Ex}[Theor]{Example}
\newtheorem{Not}[Theor]{Notation}
\newtheorem{variant}[Theor]{Variant}
\newtheorem*{Ackn}{Acknowledgments}
\newtheorem{assumption}[Theor]{Assumption}
\newtheorem{notations}[Theor]{Notations}

\newcommand{\red}{\operatorname{\sf red}}
\DeclareMathOperator{\Pic}{Pic}
\DeclareMathOperator{\Corr}{Corr}
\newcommand{\all}{\operatorname{\sf all}}
\newcommand{\proper}{\operatorname{\sf proper}}
\newcommand{\la}{\operatorname{\langle}}
\newcommand{\ra}{\operatorname{\rangle}}

\newcommand{\Stab}{\operatorname{\sf Stab}}
\newcommand{\coStab}{\operatorname{\sf coStab}}
\DeclareMathOperator{\coCAlg}{coCAlg}

\newcommand{\Tr}{\operatorname{\sf Tr}}

\DeclareMathOperator{\Vect}{Vect}
\DeclareMathOperator{\ch}{ch}
\newcommand{\aft}{\operatorname{\sf aft}}
\newcommand{\laft}{\operatorname{\sf laft}}

\DeclareMathOperator{\Rep}{Rep}
\DeclareMathOperator{\Hom}{Hom}
\newcommand{\HHom}{ {\mathcal Hom} }
\newcommand{\Funct}{\operatorname{\sf Funct}}
\DeclareMathOperator{\Map}{Map}

\newcommand{\Mmod}{\operatorname{\sf Mod}}
\DeclareMathOperator{\Cat}{Cat}
\DeclareMathOperator{\QCoh}{QCoh}
\DeclareMathOperator{\Coh}{Coh}
\DeclareMathOperator{\ICoh}{ICoh}
\DeclareMathOperator{\Perf}{Perf}
\DeclareMathOperator{\perf}{perf}

\newcommand{\st}{\operatorname{\sf st}}
\newcommand{\pr}{\operatorname{\sf pr}}

\DeclareMathOperator{\Sch}{Sch}
\DeclareMathOperator{\Aff}{Aff}
\DeclareMathOperator{\Spec}{Spec}

\newcommand{\ii}{\operatorname{(\infty,1)}}

\newcommand{\mscr}{\mathscr}
\newcommand{\bb}{\mathbb}
\newcommand{\mcal}{\mathcal}

\DeclareMathOperator{\Id}{Id}
\DeclareMathOperator{\Grp}{Grp}

\DeclareMathOperator{\Ind}{Ind}
\DeclareMathOperator{\Psh}{Psh}

\DeclareMathOperator{\PreStack}{PreStack}
\DeclareMathOperator{\CAlg}{CAlg}
\DeclareMathOperator{\Alg}{Alg}
\newcommand{\op}{\operatorname{\sf op}}

\newcommand{\Ll}{\operatorname{\sf L}}

\DeclareMathOperator{\End}{End}
\DeclareMathOperator{\Aut}{Aut}
\DeclareMathOperator{\Ad}{Ad}

\DeclareMathOperator{\GL}{GL}
\newcommand{\gl}{\mathfrak{gl}}

\DeclareMathOperator{\Lie}{Lie}
\DeclareMathOperator{\LAlg}{LAlg}
\DeclareMathOperator{\Mod}{Mod}
\DeclareMathOperator{\At}{At}
\DeclareMathOperator{\Sym}{Sym}

\definecolor{note_color}{rgb}{0.0,0.9,0.0}

\newcommand{\EEnd}{\mathcal End}
\DeclareMathOperator{\FormModuli}{\widehat{Moduli}}
\DeclareMathOperator{\td}{td}
\DeclareMathOperator{\ad}{ad}
\DeclareMathOperator{\coMod}{coMod}

\DeclareMathOperator{\cofib}{cofib}
\DeclareMathOperator{\Prim}{Prim}
\DeclareMathOperator{\HopfAlg}{HopfAlg}

\DeclareMathOperator{\Pro}{Pro}
\DeclareMathOperator{\Free}{Free}

\DeclareMathOperator{\colim}{colim}
\newcommand{\ssl}{\scriptscriptstyle}
\newcommand{\sotimes}{\overset{!}{\otimes}}
\newcommand{\sboxtimes}{\overset{!}{\boxtimes}}
\renewcommand{\epsilon}{\varepsilon}
\DeclareMathOperator{\inv}{inv}
\newcommand{\prolim}[1][]{\lim\limits_{\xleftarrow[#1]{}}}

\DeclareMathOperator{\orient}{u}

\begin{document}
\title{\textbf{Equivariant Grothendieck-Riemann-Roch theorem via formal deformation theory}}
\date{}
\author{Grigory~Kondyrev, Artem~Prikhodko}
\maketitle

\begin{abstract}
We use the formalism of traces in higher categories to prove a common generalization of the holomorphic Atiyah-Bott fixed point formula and the Grothendieck-Riemann-Roch theorem. The proof is quite different from the original one proposed by Grothendieck et al.: it relies on the interplay between self dualities of quasi- and ind- coherent sheaves on $X$ and formal deformation theory of Gaitsgory-Rozenblyum. In particular, we give a description of the Todd class in terms of the difference of two formal group structures on the derived loop scheme $\mathcal LX$. The equivariant case is reduced to the non-equivariant one by a variant of the Atiyah-Bott localization theorem.
\end{abstract}

\tableofcontents

\section{Introduction}
\begin{Conven} For the rest of the document we will assume that we work over some base field $k$ of characteristic zero.
\end{Conven}
In \cite{We} the formalism of traces in symmetric monoidal $(\infty, 2)$-categories was used to prove the following classical result
\begin{thm*}[Holomorphic Atiyah-Bott fixed point formula]
Let $X$ be a smooth proper $k$-scheme with an endomorphism $\xymatrix{X \ar[r]^-{g} & X}$ such that the graph of $g$ intersects the diagonal in $X\times X$ transversely. Then for a lax $g$-equivariant perfect sheaf $(E,t)$ (i.e. a sheaf $E\in \QCoh(X)^{\perf}$ equipped with a map $\xymatrix{E \ar[r]^-{t} & g_*E}$) there is an equality
$$
\Tr_{\Vect_k}\left(\xymatrix{\Gamma(X,E) \ar[r]^{\Gamma(X,t)} & \Gamma(X,E)}\right) = \sum_{x=g(x)} \frac{\Tr_{\Vect_k}\left(\xymatrix{E_x \simeq E_{g(x)} \ar[r]^-{t_x} & E_x}\right)}{\det(1-d_x g)}
$$
where $\xymatrix{\mathbb T_{X,x} \ar[r]^-{d_x g} & \mathbb T_{X,g(x)} \simeq \mathbb T_{X,x}}$ is the induced map on tangent spaces.
\end{thm*}

In turns out that the transversality assumption in the theorem above can be considerably weakened. For example, the following extreme opposite case (when the equivariant structure is trivial, i.e. $g= \Id_X$, $t= \Id_E$) was known a decade before the Atiyah-Bott formula
\begin{thm*}[Hirzebruch-Riemann-Roch]
Let $X$ be a smooth proper scheme over $k$. Then for every perfect sheaf $E$ on $X$ there is an equality
$$
\Tr\left(\xymatrix{\Gamma(X, E) \ar[rr]^-{\Id_{\Gamma(X,E)}} && \Gamma(X, E)} \right) = \int_X \ch(E) \td_X
$$
where $\ch(-)$ and $\td_X$ are Chern character and Todd class respectively.
\end{thm*}

The goal of this work is to provide a common generalization of the two theorems above as well as their relative versions naturally suggested by the formalism of traces. In order to state it we first need to introduce a bit of notations:
\begin{notations}
Let $X$ be a smooth scheme equipped with an endomorphism $\xymatrix{X \ar[r]^-{g} & X}$ such that the reduced classical scheme $\overline{X^g} := \mathcal H^0(X^g)^{\red}$ is smooth (but not necessary connected). We will denote by $\xymatrix{j\colon \overline{X^g} \ar[r] & X}$ the canonical embedding and by $\mathcal N_g^\vee$ its conormal bundle. Note that the action of $g$ on $\Omega_X^1$ in particular restricts to an endomorphism $\xymatrix{g^*_{|\mathcal N_g}\colon \mathcal N^\vee_g \ar[r] & \mathcal N^\vee_g}$.
\end{notations}
We then have
\begin{thm*}[Equivariant Grothendieck-Riemann-Roch, {\Cref{equiv_GRR}}]\label{intro:main_application}
Let
$$\xymatrix{X\ar@(ld,lu)^{g_X} \ar[rr]^-f && Y\ar@(rd,ru)_{g_Y}}$$
be an equivariant morphism between smooth proper schemes such that
\begin{itemize}
\item Reduced fixed loci $\overline{X^{g_X}}$ and $\overline{Y^{g_Y}}$ are smooth.

\item The induced morphisms on conormal bundles $1 - (g_X^*)_{|\mathcal N_{g_X}^\vee}$ and $1 - (g_Y^*)_{|\mathcal N_{g_Y}^\vee}$ are invertible.
\end{itemize}
Then for a lax $g_X$-equivariant perfect sheaf $E$ on $X$ there is an equality
$$
(\overline{f^g})_* \left(\ch(E,t)\frac{\td_{\overline{X^{g_X}}}}{e_{g_X}}\right) = \ch(f_*(E,t))\frac{\td_{\overline{Y^{g_Y}}}}{e_{g_Y}} \quad \in\quad \bigoplus_p H^{p,p}(\overline{Y^{g_Y}}),
$$
where $\xymatrix{\overline{X^{g_X}} \ar[r]^-{\overline{f^g}} & \overline{Y^{g_Y}}}$ is the induced map on reduced fixed loci, $\ch(-,-)$ is an equivariant Chern character (see Construction \ref{chern} and Proposition \ref{prop:semi_eq_Chern}), $\td_-$ are usual Todd classes and $e_{g_-}$ are equivariant Euler classes (see \Cref{def:euler_class} and \Cref{euler_as_sum}).
\end{thm*}

\medskip

If equivariant structures on $X$ and $Y$ are trivial (so $e_{g_X} = 1, e_{g_Y} = 1$), the theorem above reduces to the usual Grothendieck-Riemann-Roch theorem. On the other hand, if $Y=*$ and $X^{g_X}$ is discrete, then $\td_{X^{g_X}} = 1$ and we recover the holomorphic Atiyah Bott-formula (see corollaries \ref{cor:GRR} and \ref{cor:atiayh_bott} for more details).
\begin{Rem}
Even in the case of trivial equivariant structures our proof of Grothendieck-Riemann-Roch is quite different from the original approach due to Grothendieck et al.: it is valid for arbitrary smooth proper $k$-schemes $X$ and in particular does not rely on the trick of factoring a projective morphism into a composition of a closed embedding and a projection. On the other hand, due to heavy usage of deformation theory our proof works only in characteristic zero and we consider the Chern character and the Todd class as elements of Hodge cohomology $H^{*,*}(X)$, not of the Chow ring.
\end{Rem}
 
We now explain the key steps in the proof of the Equivariant Grothendieck-Riemann-Roch theorem above. 

\paragraph{Non-equivariant part.}
Let $\xymatrix{X \ar[r]^-{f} & Y}$ be a morphism of smooth proper $k$-schemes and let $E$ be a perfect sheaf on $X$. Recall that the Grothendieck-Riemann-Roch theorem asserts an equality 
$$f_*(\ch(E)\td_X) = \ch(f_*(E))\td_Y \quad \in \quad \bigoplus_p H^p(Y, \Omega_Y^p).$$
 Our first goal is to give a proof of the Grothendieck-Riemann-Roch theorem using the formalism of traces. Note that the formula above is equivalent to the commutativity of the diagram
\begin{align}\label{comutativity_to_prove}
\xymatrix{
k \ar[rr]_-{\ch(E)} \ar@/^2pc/[rrrr]^{\ch(f_*(E))} && \displaystyle \bigoplus_p H^p(X, \Omega_X^p) \ar[d]^-{\cdot \td_X}_-{\sim}  && \displaystyle \bigoplus_p H^p(Y, \Omega_Y^p) \ar[d]^-{\cdot \td_Y}_-{\sim}
\\
 && \displaystyle \bigoplus_p H^p(X, \Omega_X^p) \ar[rr]_-{f_*} && \displaystyle \bigoplus_p H^p(Y, \Omega_Y^p)}
\end{align}
of $k$-vector spaces. We now recall
\begin{Prop}[{\cite[Definition 1.2.6]{We}}] \label{2trace}
Let $\mathscr C, \mathscr D \in \Cat_k$ be a pair of dualizable $k$-linear presentable categories. Suppose we are given a (not necessarily commutative) diagram
$$\xymatrix{
\mathscr C \ar@/_/[dd]_-{\varphi} \ar[rr]^-{F_{\mathscr C}} && \mathscr C \ar@/_/[dd]_-{\varphi} \ar@2[ddll]_-{T}
\\
\\
\mathscr D \ar@/_/[uu]_-{\psi} \ar[rr]_-{F_{\mathscr D}} && \mathscr D \ar@/_/[uu]_-{\psi},
}$$
where $\varphi$ is left adjoint to $\psi$ and $\xymatrix{ \varphi \circ F_{\mathscr C} \ar[r]^-{T} & F_{\mathscr D} \circ \varphi}$ is a (not necessary invertible) natural transformation. Then there exists a natural morphism
$$\xymatrix{
\Tr_{2\Cat_k}(F_{\mathscr C}) \ar[rrr]^-{\Tr_{2\Cat_k}(\varphi,T)} &&& \Tr_{\Cat_k}(F_{\mathscr D})
}$$
in the $\ii$-category $\Hom_{2\Cat_k}(\Vect_k ,\Vect_k) \simeq \Vect_k$ called the \emph{morphism of traces induced by $T$} (in the case when $T=\Id_{\varphi \circ F_{\mathscr C}}$ we will further frequently use the notation $\Tr_{2\Cat_k}(\varphi):=\Tr_{2\Cat_k}(\varphi,\Id_{\varphi \circ F_{\mathscr C}})$). Moreover the morphism of traces is functorial in an appropriate sense, see \cite[Proposition 1.2.11.]{We}.
\end{Prop}

We refer the reader to \cite[Section 1.2]{We} for a detailed discussion of the formalism of traces (see also \cite{HSS_ntraces} for a more general treatment in the context of $(\infty, n)$-categories and \cite{BN} for traces in the context of derived algebraic geometry). The main idea of our proof of the Grothendieck-Riemann-Roch theorem is that one can obtain commutativity of the diagram (\ref{comutativity_to_prove}) above as a corollary of functoriality of the construction of the morphism of traces. Namely, let $2\Cat_k$ be the $(\infty,2)$-category of $k$-linear stable presentable categories and continuous functors. Recall that to each derived scheme $Z$ we can associate the $\ii$-category $\QCoh(Z) \in 2\Cat_k$ of quasi-coherent sheaves on $Z$ (see \cite[Chapter 3]{GaitsRozI}) and a closely related category $\ICoh(Z)$ of ind-coherent sheaves on $Z$ (see \cite[Chapter 4]{GaitsRozI}). Then by applying functoriality of traces to the diagram
$$\xymatrix{
\Vect_k \ar[rr]^-{E\otimes -} && \QCoh(X) \ar[d]^{\otimes \mathcal O_X}_\sim \ar[rr]^-{f_*} && \QCoh(Y)\ar[d]^-{\otimes \mathcal O_Y}_\sim 
\\
 && \ICoh(X) \ar[rr]_-{f_*} && \ICoh(Y)
}$$
in $2\Cat_k$ we obtain a commutative diagram of traces
$$\xymatrix{
k \ar[rr]_-{\Tr_{2\Cat_k}(E\otimes-)} \ar@/^2pc/[rrrr]^{\Tr_{2\Cat_k}(f_*(E)\otimes-)} && \Tr_{2\Cat_k}(\Id_{\QCoh(X)}) \ar[d]^-{\Tr_{2\Cat_k}(-\otimes \mathcal O_X)} \ar[rr] && \Tr_{2\Cat_k}(\Id_{\QCoh(Y)}) \ar[d]^-{\Tr_{2\Cat_k}(-\otimes \mathcal O_Y)}
\\
&& \Tr_{2\Cat_k}(\Id_{\ICoh(X)}) \ar[rr]_-{\Tr_{2\Cat_k}(f_*)} && \Tr_{2\Cat_k}(\Id_{\ICoh(X)})
}$$
in $\Vect_k$. The bulk of the paper will be devoted to identifying the morphism of traces in the diagram above with their classical counterparts. Namely
\begin{itemize}
\item Section \ref{sect:Chern}: we prove
 $$\pi_0 \Tr_{2\Cat_k}(\Id_{\QCoh(X)}) \simeq \bigoplus_p H^p(X, \Omega_X^p) \qquad \pi_0 \Tr_{2\Cat_k}(\Id_{\QCoh(Y)}) \simeq \bigoplus_p H^p(Y, \Omega_Y^p).$$
Moreover, under the isomorphisms above $\Tr(E\otimes-)$ and $\Tr(f_*(E)\otimes-)$ will coincide with the usual Chern characters of $E$ and $f_*(E)$ respectively. In fact, our description of $\Tr(E\otimes-)$ will be in terms of the Atiyah class of $E$ and is closely related to the description given in \cite{Markarian}. 

\item Section \ref{sect:icoh}: we have 
$$\pi_0 \Tr_{2\Cat_k}(\Id_{\ICoh(X)}) \simeq \bigoplus_p H^p(X, \Omega_X^p)^\vee \qquad \pi_0 \Tr_{2\Cat_k}(\Id_{\ICoh(Y)}) \simeq \bigoplus_p H^p(Y, \Omega_Y^p)^\vee.$$
Moreover, under the two isomorphisms above the morphism of traces induced by the pushforward functor $\xymatrix{\ICoh(X) \ar[r]^-{f_*} & \ICoh(Y)}$ coincides with the usual pushforward in homology (defined as the Poincar\'e dual of the pullback).

\item Sections \ref{sect:CY}, \ref{sect:todd}: under the Poincar\'e self-duality
$$\bigoplus_p H^p(X, \Omega_X^p) \simeq \bigoplus_p H^p(X, \Omega_X^p)^\vee$$
the morphism $\Tr_{2\Cat_k}(-\otimes \mathcal O_X)$ is given by the multiplication with the Todd class $\td_X$ and analogously for $Y$.
\end{itemize}
Using these identifications and the commutative diagram of traces above, one immediately concludes the Grothendieck-Riemann-Roch theorem.

\paragraph{Equivariant part.}
We now discuss how one can get an equivariant version of the GRR theorem. Let $X$ be a smooth proper scheme with an endomorphism $g$ and $(E,t)$ be a lax $g$-equivariant perfect sheaf on $X$. First, it turns out that if $g = \Id_X$ there is simple description of the equivariant Chern character $\ch(E,t)$ (see Construction \ref{chern}) in the spirit of the Chern-Weil theory (see \Cref{ChernWeild_and_Atiyah}):

\begin{Ex}
Let $X$ be a smooth proper scheme with a trivial equivariant structure. Then under the identification
$$\Tr_{2\Cat_k}(\Id_{\QCoh(X)}) \simeq \bigoplus_p H^p(X, \Omega_X^p)$$
the equivariant Chern character is equal to $\Tr(e^{\At(E)} \circ t)$, where $\At(E)$ is the Atiyah class of $E$ (see Corollary \ref{prop:semi_eq_Chern}).
\end{Ex}

Now if the equivariant structure on $X$ is non-trivial, we in general do not have a convenient description of $\Tr_{2\Cat_k}(g_*)$ and therefore of $\ch(E,t)$. To circumvent this, we shall reduce the situation to the non-equivariant case by restricting along the (classical, reduced) fixed point locus $\xymatrix{j\colon \overline{X^g} \ar[r] & X}$ of the endomorphism $g$, like in the Atiyah-Bott localization theorem. Note that this indeed gives a good description of the Chern character: since the restriction of the lax $g$-equivariant sheaf $(E,t)$ to $\overline{X^g}$ is naturally equivariant with respect to the trivial equivariant structure on $\overline{X^g}$, by the example above we have a grasp on 
$$
\ch(j^*(E,t)) \in \pi_0\Tr_{2\Cat_k}(\Id_{\QCoh(\overline{X^g})}) \simeq \bigoplus_p H^p(\overline{X^g}, \Omega_{\overline{X^g}}^p).
$$
Moreover, we will show that under reasonable assumptions the morphism $j^*$ in fact does not loose any information:
\begin{thm*}[Localization theorem, {\ref{localization_thm}}]
Assume that the reduced classical fixed locus $\overline{X^g}$ is smooth and let $\mathcal N_g^\vee$ be the conormal bundle of $j$. Then the induced map
$$\xymatrix{\pi_0\Tr_{2\Cat_k}(g_*) \ar[rr]^-{\Tr_{2\Cat_k}(j^*)} && \pi_0\Tr(\Id_{\QCoh(\overline{X^g})})} \simeq \bigoplus_{p} H^p(\overline{X^g}, \Omega_{\overline{X^g}}^p)$$
is an equivalence if and only if the determinant $\det(1 - g^*_{|\mathcal N_g^\vee})$ is invertible.
\end{thm*}

Let now
$$\xymatrix{X\ar@(ld,lu)^{g_X} \ar[rr]^-f && Y\ar@(rd,ru)_{g_Y}}$$
be an equivariant morphism from $(X,g_X)$ to $(Y,g_Y)$ satisfying assumptions of the Equivariant Grothendieck-Riemann-Roch theorem. Then by applying $\Tr_{2\Cat_k}$ to the commutative diagram
$$\xymatrix{
\Vect_k \ar[rr]^-E && \QCoh(X)\ar@(ul,ur)^{g_{X*}} \ar[d]^{\otimes \mathcal O_X}_\sim \ar[rr]^-{f_*} && \QCoh(Y)\ar@(ul,ur)^{g_{Y*}} \ar[d]^-{\otimes \mathcal O_Y}_\sim 
\\
 && \ICoh(X)\ar@(dl,dr)_{g_{X*}} \ar[rr]_-{f_*} && \ICoh(Y)\ar@(dl,dr)_{g_{Y*}}
}$$
and using the identification from the localization theorem $\pi_0\Tr_{2\Cat_k}((g_X)_*) \simeq \bigoplus_{p} H^{p,p}(\overline{X^{g_X}})$ and analogously for $Y$, we obtain a commutative diagram
$$\xymatrix{
k \ar[rr]_-{\ch(E,t)} \ar@/^2pc/[rrrr]^{\ch(f_*(E,t))} && \displaystyle\bigoplus_p H^{p,p}(\overline{X^{g_X}}) \ar[d]^-{\cdot \td_{g_X}}_-{\sim} \ar[rr] && \displaystyle\bigoplus_p H^{p,p}(\overline{Y^{g_Y}}) \ar[d]^-{\cdot \td_{g_Y}}_-{\sim}
\\
 && \displaystyle\bigoplus_p H^{p,p}(\overline{X^{g_X}}) \ar[rr]_-{(\overline{f^g})_*} && \displaystyle\bigoplus_p H^{p,p}(\overline{Y^{g_Y}})
}$$
in $\Vect_k$ for some equivalences $\td_{g_X}$ and $\td_{g_Y}$ (see \Cref{equiv_todd}). Finally, by applying the commutative diagram above to the special case when the morphism is $\overline{X^{g_X}} \xymatrix{\ar@{^(->}[r] &} X$ and using the non-equivariant part it is straightforward to identify $\td_{g_X}$ with $\frac{\td_{\overline{X^{g_X}}}}{e_{g_X}}$ (\Cref{eq_Todd_in_classical_terms}), where $e_{g_X}$ is the Euler class (see \Cref{def:euler_class} and \Cref{euler_as_sum}) and similarly for $Y$. 

\smallskip 

\begin{Rem}
Even if one is only interested in the equivariant part, in order to reduce it to the non-equivariant version in our approach it is still necessary to identify the morphism of traces
$$\xymatrix{
\Tr_{2\Cat_k} (\Id_{\QCoh(\overline{X^g})}) \ar[rrr]^{\Tr_{2\Cat_k}(-\otimes \mathcal O_{\overline{X^g}})} &&& \Tr_{2\Cat_k}(\Id_{\ICoh(\overline{X^g})})
}$$
in classical terms. In other words, we really need to reprove (non-equivariant) Grothendieck-Riemann-Roch theorem in the language of traces.
\end{Rem}

\subsection*{Relation to previous work}

Our identification of the morphism of traces in the non-equivariant case is closely related to the work of Markarian \cite{Markarian}, where the role of the canonical action of the Lie algebra $\mathbb T_X[-1]$ on any object of $\QCoh(X)$ (via Atiyah class) was emphasized. Moreover, he provides an interpretation of the Todd class $\td_X$ as an invariant volume form with respect to the Hopf algebra structure on Hochschild homology. The key difference of our approach is the systematic use of the formalism of traces, derived algebraic geometry, ind-coherent sheaves and deformation theory of Gaitsgory-Rozenblyum, allowing us to give a precise geometric interpretation of Markarian's ideas in terms of formal groups over $X$ and to generalize them in some directions. The idea of such an interpretation of $\td_X$ was explained to us by Dennis Gaitsgory who in turn learned it from Maxim Kontsevich, according to whom this idea ultimately goes back to Boris Feigin.

Ideas similar to our trace formalism were used in a series of papers \cite{BFN}, \cite{BN}, \cite{BN_sec} to derive many interesting trace formulas in various contexts (such as D-modules on schemes or quasi-coherent sheaves on stacks). However, as explained in \cite[Remark 1.6]{BN}, some additional work is needed to fully recover the usual GRR-theorem from their results. Among other things, in this paper we provide a necessary identifications of morphisms of traces in classical context by studying the formal geometry of the derived loop scheme.

Localization theory relating equivariant cohomology and ordinary cohomology of the fixed locus was studied for a long time already, see e.g. \cite{AtiyahBott_Original} for the de Rham variant of the theory. For various reasons the methods of \cite{AtiyahBott_Original} don't apply directly in our context (e.g. since equivariant cohomology of a point is too small in the case of a plain endomorphism to invert something there), but the idea that equivariant cohomology and ordinary cohomology of the fixed locus agree up to localization of the Euler class is ultimately motivated by their work.

Formality of derived fibered products was thoroughly studied in \cite{ACH_formality}. As an application they prove (\cite[Corollary 1.12]{ACH_formality}) that for a finite order automorphism the derived fixed locus is always formal. Hence our localization theorem \ref{localization_thm} may be considered as a generalization of this result, providing formality criterion for an arbitrary endomorphism.

A similar result to that of ours was obtained in \cite{Donovan_equivariant_finite_auto} over an arbitrary algebraically closed field but with an additional assumption that both $X,Y$ are projective and that both $g_X, g_Y$ are automorphisms of finite order coprime to the characteristic of the base field. Donovan's proof is quite different from ours and is much closer to the original approach due to Grothendieck: it relies on the fact that an equivariant projective morphism $\xymatrix{f: (X, g_X) \ar[r] & (Y, g_Y)}$ (with $g_X, g_Y$ being automorphisms of finite orders) can be factored into an equivariant closed embedding into relative projective space followed by projection. We do not see how to generalize this approach to the case of a general endomorphism. On the other hand, Donovan's formula holds in the Chow ring, while our proof gives equality only in Hodge cohomology.

One can also extend the classical Grothendieck-Riemann-Roch theorem in a non-commutative direction by studying traces of maps on Hochschild homology induced by functors between nice enough (e.g. smooth and proper) categories (see \cite{Shklyarov}, \cite{Lunts}, \cite{Polishchuk}, \cite{DG}). But while the Chern character makes perfect sense in the non-commutative context, it was already pointed out by Shklyarov that the Todd class seems to be of commutative nature and is missing in general non-commutative GRR-like theorems. Hence the \Cref{equiv_GRR} may be considered as a refinement of the more general non-commutative versions under additional geometricity assumptions on categories and endofunctors.

Finally, there is a categorified version of the GRR-theorem conjectured in \cite{TV_CGRR_in_progress}, where (among many other things) the role of the derived loops construction was empasized. See \cite{TV_CGRR} and {\cite{HSS_RR} for proofs and interesting applications.

\begin{Ackn}
The idea to use the categorical trace to prove the Atiyah-Bott formula and Grothendieck-Riemann-Roch was suggested to us by Dennis Gaitsgory. We are grateful to him for guidance throughout this project. We also would like to thank Chris Brav for carefully reading drafts of this text and providing many useful suggestions. The second author wants to thank Pavel Popov for helpful discussions about localization theorems.

The first author was supported by the HSE University Basic Research Program, Russian Academic Excellence Project '5-100'. The second author is partially supported by Laboratory of Mirror Symmetry NRU HSE, RF Government grant, ag. № 14.641.31.0001 and by the Simons Foundation.
\end{Ackn}

\medskip
\begin{Conv}\
\\
1) All the categories we work with are assumed to be $\ii$-categories. For an $\ii$-category $\mscr{C}$ we will denote by $(\mscr{C})^{\simeq}$ the underlying $\infty$-groupoid of $\mscr{C}$ obtained by discarding all the non-invertible morphisms from $\mscr{C}$.
\\
\\
2) We will denote by $\mcal{S}$ the symmetric monoidal $\ii$-category of spaces. For a field $k$ we will denote by ${\Vect_k}$ the stable symmetric monoidal $\ii$-category of unbounded cochain complexes over $k$ up to quasi-isomorphism with the canonical $(\infty,1)$-enhancement. We will also denote by $\Vect_k^{\heartsuit}$ the ordinary category of $k$-vector spaces considered as an $\ii$-category.
\\
\\
3) We will denote by $\Pr^{\Ll}_{\infty}$ the $\ii$-category of presentable $\ii$-categories and continuous functors with the symmetric monoidal structure from \cite[Proposition 4.8.1.15.]{HA}. Similarly, we will denote by $\Pr^{\Ll,\st}_{\infty}$ the $\ii$-category of stable presentable $\ii$-categories and continuous functors considered as a symmetric monoidal $\ii$-category with the monoidal structure inherited from $\Pr^{\Ll}_{\infty}$.
\\
\\
4) Notice that $\Vect_k$ is a commutative algebra object in $\Pr^{\Ll,\st}_{\infty}$. By \cite[Theorem 4.5.2.1.]{HA} it follows that the presentable stable $\ii$-category of $k$-linear presentable $\ii$-categories and $k$-linear functors $\Cat_k:=\Mod_{\Vect_k}(\Pr^{\Ll,\st}_{\infty})$ admits natural symmetric monoidal structure. We will also denote by $2\Cat_k$, the symmetric monoidal $(\infty,2)$-category of $k$-linear presentable $\ii$-categories and continuous $k$-linear functors so that the underlying $\ii$-category of $2\Cat_k$ is precisely $\Cat_k$.
\\
\\
5) We will denote by $\PreStack$ the $\ii$-category of functors $\Funct(\CAlg^{\leq 0}_k, \mcal{S})$, where $\CAlg_k^{\leq 0}:=\CAlg(\Vect_k)^{\leq 0}$ is the $\ii$-category of connective commutative algebras in $\Vect_k$. For a prestack $X \in \PreStack_k$ we will denote the $k$-linear symmetric monoidal $\ii$-category of unbounded complexes of quasi-coherent sheaves on $X$ by $\QCoh(X) \in \CAlg(\Cat_k)$. We refer the reader to \cite[Part I]{GaitsRozI} for an introduction to the basic concepts of derived algebraic geometry.
\\
\\
6) By 'scheme' we will always mean a derived schemes (in the sense of \cite{GaitsRozI}) if not explicitly stated otherwise. For a smooth classical scheme $X$ we will sometimes denote its Hodge cohomology $H^q(X, \Omega_X^p)$ by $H^{p,q}(X)$.
\end{Conv}

\section{Categorical Chern character}\label{sect:Chern}
Let $X$ be an almost finite type scheme (see \cite[Section 3.5.1]{GaitsRozI}) with an endomorphism $\xymatrix{X \ar[r]^-{g} & X}$ and let $(E,t) \in \QCoh(X)$ be a lax $g$-equivariant compact quasi-coherent sheaf on $X$ (i.e. a sheaf $E \in \QCoh(X)$ together with a morphism $\xymatrix{E \ar[r]^-{t} & g_* E}$). Our goal in this section is to describe the categorical Chern character $\ch(E, t) \in \Tr_{2\Cat_k}(g_*)$ obtained by applying the formalism of traces to the commutative diagram
$$\xymatrix{
\Vect_k \ar@/_/[dd]_-{- \otimes E} \ar[rr]^-{\Id_{\Vect_k}} && \Vect_k \ar@/_/[dd]_-{- \otimes E}  \ar@2[ddll]^-{T}
\\
\\
\QCoh(X) \ar@/_/[uu]_-{\Hom(E, -)} \ar[rr]_-{g_*} && \QCoh(X) \ar@/_/[uu]_-{\Hom(E, -)}
}$$
in more concrete terms. This will be done in several steps:
\begin{itemize}
\item First, using that the assignment $\xymatrix{X \ar@{|->}[r] & \QCoh(X)}$ lifts to a functor from an appropriate category of correspondences we will identify
$$\Tr_{2\Cat_k}(g_*) \simeq \Gamma(X^g, \mathcal O_{X^g})$$
where $X^g$ is the derived fixed locus of $g$.

\item Specializing to the case when $g=\Id_X$ so that $X^g \simeq \mathcal LX$ is the so-called \emph{inertia group of $X$} we will show that $E$ admits a canonical action of $\mathcal LX$ and compute $\ch(E,t)$ in terms of this action.

\item Finally, under additional assumption that $X$ is smooth (so that we can apply $\QCoh$-version of formal deformation theory) we will show that the canonical action of the $\mathcal LX$ on $E$ is closely related to the canonical action of the Lie algebra $\mathbb T_X[-1]$ given by the Atiyah class of $E$ deducing that
$$\ch(E,t) = \Tr(\exp(\At(E))\circ t).$$
In particular, if $t=\Id_E$ we will show that $\ch(E, \Id_E)$ coincides with the classical topological Chern character (defined using the splitting principle).
\end{itemize}

\subsection{Self-duality of quasi-coherent sheaves and Chern character}
We start with a short reminder of self-duality of $\QCoh$. Recall the following
\begin{Theor}[{\cite[Chapter 3, Proposition 3.4.2 and Chapter 6, Proposition 4.3.2]{GaitsRozI}}]\label{QCohSelfDual}
\
\begin{enumerate}[1.]
\item For any two $X,Y \in \Sch_{\aft}$ (for the definition see \cite[Chapter 2, 3.5]{GaitsRozI}) the morphism
$$\xymatrix{
\QCoh(X) \otimes \QCoh(Y) \ar[rr] && \QCoh(X \times Y)
}$$
in $\Cat_k$ induced by the functor
$$\xymatrix{
\QCoh(X) \times \QCoh(Y) \ar[rr]^-{\boxtimes} && \QCoh(X \times Y)
}$$
is an equivalence.

\item For any $X \in \Sch_{\aft}$ the morphisms
$$\xymatrix{
\Vect_k \ar[rr]^-{\Delta_* \mcal{O}_X} && \QCoh(X \times X) \simeq \QCoh(X) \otimes \QCoh(X)
}$$
and
$$\xymatrix{
\QCoh(X) \otimes \QCoh(X) \simeq \QCoh(X \times X) \ar[rr]^-{\Gamma \circ \Delta^*} && \Vect_k
}$$
$$\xymatrix{
}$$
exhibit $\QCoh(X)$ as a self-dual object in $\Cat_k$.
\end{enumerate}
\end{Theor}

\begin{Cor}\label{TraceOfQCoh}
Let $X$ be an almost finite type scheme with an endomorphism $\xymatrix{X \ar[r]^-{g} & X}$. Then there is an equivalence
$$
\Tr_{2\Cat_k}(g_*) \simeq \Gamma(X, \Delta^* (\Id_X, g)_* \mcal{O}_X) \simeq \Gamma(X^g, \mathcal O_{X^g}),
$$
where $X^g$ is the derived fixed locus of $g$ defined by the pullback square
$$\xymatrix{
X^g \ar[r]^-{i} \ar[d]_-{i} & X \ar[d]^-{(\Id_X, g)}
\\
X \ar[r]_-{\Delta} & X \times X.
}$$

\begin{proof}
Unwinding the definitions and using the theorem above one finds that the trace $\Tr_{2\Cat_k}(g_*)$ is given by the composite
$$\xymatrix{
\Vect_k \ar[rr]^-{\Delta_* \mcal{O}_X} && \QCoh(X \times X) \ar[rr]^-{(\Id_X \times g)_*} && \QCoh(X\times X) \ar[rr]^-{\Gamma \circ \Delta^*} && \Vect_k
}$$
which is given by tensoring with $\Gamma(X, \Delta^* (\Id_X, g)_* \mcal{O}_X)$. Now by applying the base change for quasi-coherent sheaves (see \cite[Chapter 3, Proposition 2.2.2.]{GaitsRozI}) to the diagram above we get
$$
\Gamma(X,\Delta^* (\Id_X, g)_*  \mcal{O}_X) \simeq \Gamma(X,i_* i^*  \mcal{O}_X) \simeq \Gamma(X,i_* \mathcal O_{X^g}) \simeq \Gamma(X^g, \mathcal O_{X^g})
$$
proving the claim.
\end{proof}
\end{Cor}

\begin{Ex}\label{ex:meet_derived_loops}
Note that in the case when $g=\Id_X$ we get an equivalence $X^g \simeq \mathcal LX$ where $\mathcal{L}X$ is the {\bfseries derived loops scheme of $X$} defined as
$$
\mathcal LX := \Map(S^1, X) \simeq X \underset{X\times X}{\times} X
$$
In particular, in the case when $X$ is smooth the Hochschild-Kostant-Rosenberg isomorphism (see Corollary \ref{HochKosRos}) states that
$$\pi_i\Gamma(\mathcal LX, \mathcal O_{\mathcal LX}) \simeq \bigoplus_{p-q=i}^{\dim X} H^q(X, \Omega^p_X).$$
As we will in Theorem \ref{localization_thm} under some reasonable assumptions one can use the Hochschild-Kostant-Rosenberg isomorphism to get some understanding of $\Gamma(X^g,\mathcal{O}_{X^g})$ for $g$ which is not necessary an identity.
\end{Ex}

We now turn to the categorical Chern character. Recall from \cite[Definition 1.2.9]{We} the following
\begin{construction}\label{chern}
Given a compact sheaf $E \in \QCoh(X)$ together with an endomorphism $\xymatrix{E \ar[r]^-{t} & g_*E}$ we can form a diagram
$$\xymatrix{
\Vect_k \ar@/_/[dd]_-{- \otimes E} \ar[rr]^-{\Id_{\Vect_k}} && \Vect_k \ar@/_/[dd]_-{- \otimes E}  \ar@2[ddll]^-{T}
\\
\\
\QCoh(X) \ar@/_/[uu]_-{\Hom(E, -)} \ar[rr]_-{g_*} && \QCoh(X) \ar@/_/[uu]_-{\Hom(E, -)}
}$$
in $2\Cat_k$ with the $2$-morphism $T$ induced by the morphism $t$. The corresponding element 
$$\ch(E,t) \in \Tr_{2\Cat_k}(g_*) \simeq \Gamma(X^g, \mathcal O_{X^g})$$
obtained via the formalism of traces (see Proposition \ref{2trace}) is called a \emph{categorical Chern character of $E$}. 
\end{construction}

Recall also the following result
\begin{Prop}[{\cite[Proposition 2.2.3.]{We}}]\label{prop23}
Assume that $X$ is a quasi-compact scheme (in this case $\QCoh(X)$ is compactly generated by dualizable objects). Given an endomorphism $\xymatrix{X \ar[r]^-{g} & X}$ and a dualizable sheaf $E \in \QCoh(X)$ together with a lax $g$-equivariant structure $\xymatrix{E \ar[r]^-{t} & g_* E}$ on $E$ there is an equality
$$
\ch(E,t) \simeq \Tr_{\QCoh(X^g)}\left(\xymatrix{i^* E \ar[r]^-{\beta}_-{\sim} & i^* g^* E \ar[r]^-{i^*(b)} & i^* E}\right)
$$
in $\Gamma(X^g, \mathcal O_{X^g})$, where $b \in \Hom_{\QCoh(X)}(g^*E,E)$ is the morphism which corresponds to $t$ using the adjunction between $g^*$ and $g_*$, the morphism $\xymatrix{X^g \ar@{^(->}[r]^i & X}$ is the inclusion of the derived fixed points and the equivalence $\beta$ is induced by the equivalence $i \simeq g \circ i$.
\end{Prop}

It turns out that the equivalence $\beta$ above is non-trivial even in the case $g = \Id_X$. In fact, it is closely related to the Atiyah class of $E$, as we will see further in this section.

\subsection{Canonical $\mathcal LX$-equivariant structure}
In this subsection we will construct a canonical action of $\mathcal LX$ on any sheaf in $\QCoh(X)$ and use this structure to give an alternative description to the morphism $\beta$ from Proposition \ref{prop23}. We start with the following

\begin{Def} Let $X \in \PreStack$ be a prestack and $Y \in \PreStack_{/X}$ be a prestack over $X$. Then for a group $G \in \Grp(\PreStack_{/X})$ which acts on $Y$ define a prestack $Y/G \in \PreStack$ as the colimit
$$
Y/G:=\underset{\Delta^{\sf op}}{\colim}\left(\xymatrix{\ldots \ar@<+1.35ex>[r] \ar@<+.45ex>[r] \ar@<-.45ex>[r] \ar@<-1.35ex>[r] & G \times_X G \times_X Y \ar@<+.9ex>[r] \ar[r] \ar@<-.9ex>[r]  & G \times_X Y \ar@<-.45ex>[r]_-{q_2} \ar@<.45ex>[r]^-{a} & Y}\right)$$
in the category of prestacks.
\end{Def}

Recall now the following
  
\begin{Def} Let $X \in \PreStack$ be a prestack and $Y \in \PreStack_{/X}$ be a prestack over $X$. Then for a group $G \in \Grp(\PreStack_{/X})$ which acts on $Y$ we define a {\bfseries category of $G$-equivariant sheaves on $Y$} denoted by $\Rep_{G}(\QCoh(Y))$ simply as
$$
\Rep_{G}(\QCoh(Y)):=\QCoh(Y/G).
$$
\end{Def}
We will further frequently abuse the notation by considering $G$-equivariant sheaves on $X$ as objects of the category $\QCoh(Y)$ via the pullback functor
$$\xymatrix{
\QCoh(Y/G) \ar[rr]^-{h^*}_-{\sim} && \QCoh(Y)
}$$
where $\xymatrix{Y \ar[r]^-{h} & Y/G}$ is the natural projection map.
\begin{construction}\label{act_morphism}
Notice that by the definition of $Y/G$ we have an equivalence of functors $\xymatrix{a^* h^* \ar[r]_-{\sim} & q_2^* h^*}$. Consequently, for every $G$-equivariant quasi-coherent sheaf $\mcal{F} \in \Rep_{G}(\QCoh(Y))$ we get an equivalence $\xymatrix{a^* \mcal{F} \ar[r]^-{\alpha_{\mcal{F}}}_-{\sim} & q_2^* \mcal{F}}$ which we will further call an \emph{action morphism of $G$ on $\mathcal F$}.
\end{construction}
\begin{Rem}\label{our_character}
Suppose that the action of $G$ on $Y$ is trivial so that $a=q_2$ and therefore for any $\mathcal F\in \Rep_G(Y)$ the morphism $\alpha_{\mathcal F}$ is an automorphism. For dualizable $\mathcal F \in \QCoh(X)$ we then define a \emph{character $\chi_{\mathcal F}$ of the representation $\mathcal F$} as the trace
$$
\chi_{\mathcal F} := \Tr_{\QCoh(G\times_X Y)}(\alpha_{\mathcal F}) \in \pi_0\Gamma(G\times_X Y, \mathcal O_{G \times_X Y})
$$
of the action morphism. To justify the notation above, suppose that $X=Y=\Spec k$ and $G$ is an ordinary group scheme over a field $k$. In this case $E$ corresponds to some representation $\rho\in \Rep_G(\Vect_k)$ and for every $g\in G(k)$ the pullback $g^*\alpha_{\mathcal F}$ is equal to $\rho(g)$. In particular, since the pullback functor is symmetric monoidal we get
$$
\chi_{\mathcal F}(g): = g^* \Tr_{\QCoh(G\times_X Y)}(\alpha_{\mathcal F}) = \Tr(g^*\alpha_{\mathcal F}) = \Tr(\rho(g))
$$
i.e. $\chi_{\mathcal F}$ defined as above coincides with the usual character of representation $\rho$.
\end{Rem}

Notice now that for a morphism $\xymatrix{Z \ar[r]^-{p} & X}$ in $\PreStack$ together with a section $\xymatrix{X \ar[r]^-{s} & Z}$ of $p$ the fiber product $G:=X \times_Z X$ admits the structure of a group object over $X$: the morphisms $p$ and $s$ above realize $Z$ as a pointed object in $\PreStack_{/X}$ and the fiber product $X \times_Z X$ is simply the loop object of $Z \in \PreStack_{/X}$ at the point $s$. Moreover, we can regard $X$ as acted on by $G$ by the (necessarily) trivial action and there is an induced map $\xymatrix{X/G \ar[r] & Z}$ (which is equivalence if and only if the section $s$ is an effective epimorphism).

In this paper we are interested in the special case when $Z:=X \times X$, the morphism $p:=q_1$ is given by the projection to the first component and the section $s:=\Delta$ is given by the diagonal morphism. In this case the resulting group
$$X \underset{X \times X}{\times} X \in \Grp(\PreStack_{/X})$$
is called the \emph{inertia group} of $X$ and by definition coincides with the derived loop stack $\mathcal LX$ of $X$ from Example \ref{ex:meet_derived_loops}. In particular, by the discussion above we see that $\mcal{L}X$ naturally acts on $X$ and that we have a map $\xymatrix{B_{/X} \mcal{L}X \ar[r]^c & X \times X}$ (where $B_{/X}\mathcal{L}X \in \PreStack_{/X}$ is the delooping of $\mcal{L}X$ over $X$). We are now ready to state the following

\begin{Prop}\label{CanonicalStructure} Every $\mcal{F} \in \QCoh(X)$ admits a natural $\mcal{L}X$-equivariant structure.
\begin{proof}
The desired structure is given by the composite
$$\xymatrix{
\QCoh(X) \ar[r]^-{q_2^*}& \QCoh(X \times X) \ar[r]^-{c^*} & \QCoh(B_{/X} \mcal{L}X)=\Rep_{\mcal{L}X}(\QCoh(X))
}$$
where $\xymatrix{X \times X \ar[r]^-{q_2} & X}$ is the projection to the second component.
\end{proof}
\end{Prop}
Following \cite[Chapter 8, 6.2]{GaitsRozII} we will further call this structure the \emph{canonical $\mcal{L}X$-equivariant structure} on $\mcal{F}$.

\begin{Rem}
Notice that if one would take the first projection $q_1$ in the proposition above instead the resulting $\mcal{L}X$-equivariant structure on every quasi-coherent sheaf on $X$ will be trivial (since we consider $X\times X$ as an object over $X$ precisely via $q_1$). The canonical $\mathcal LX$-equivariant structure is, however, highly nontrivial as we will see just below.
\end{Rem}

\medskip

We are now ready to describe how the categorical Chern character can be understood using the natural action of $\mcal{L}X$ on $X$:
\begin{Prop}\label{23action}
Let $\xymatrix{E \ar[r]^-{t} & E}$ be a dualizable quasi-coherent sheaf on $X$ together with an endomorphism. Then
$$
\ch(E,t) = \Tr_{\QCoh(\mathcal LX)}\left(i^*(t) \circ \alpha_E \right) \in \pi_0 \Gamma(\mcal{L}X,\mcal{O}_{\mcal{L}X})
$$
where $\ch(E,t) \in \Tr_{2\Cat_k}(\Id_{\QCoh(X)}) \simeq \Gamma(\mcal{L}X,\mcal{O}_{\mcal{L}X})$ is the categorical Chern character and $\alpha_E$ is the action morphism of $\mcal{L}X$ on $E$ from Construction \ref{act_morphism} (here $E$ is considered as an $\mcal{L}X$-equivariant quasi-coherent sheaf on $X$ using the canonical $\mcal{L}X$-equivariant structure from Proposition \ref{CanonicalStructure}).
\begin{Rem}
Note that in the case $t=\Id_E$ the right hand side coincides by definition with the character of the canonical representation of $\mathcal LX$ on $E$ from Remark \ref{our_character}.
\end{Rem}

\begin{proof}
Let more generally $g$ be arbitrary endomorphism of $X$ and consider the pullback diagram defining $X^g$
$$\xymatrix{
X^g \ar[r]^-{i} \ar[d]_-{i} & X \ar[d]^-{(\Id_X,g)}
\\
X \ar[r]_-{\Delta} & X \times X \ar[r]_-{q_2} & X.
}$$
Unwinding the definitions, one finds that the morphism $\beta$ from Proposition \ref{prop23} can be rewritten as the composite
$$\xymatrix{
i^* E \simeq i^* \Delta^* q_2^* E \ar[r]_-{\sim} & i^* (\Id_X,g)^* q_2^* E \simeq i^* g^* E,
}$$
where the middle morphism is induced by the pullback diagram above. In particular, in the special case when $g:=\Id_X$ and $b:=t$ we get an equivalence
$$\xymatrix{
\ch(E,t) = \Tr_{\QCoh(\mcal{L}X)}(i^* E \ar[r]^-{\beta}_-{\sim} & i^* E \ar[r]^{i^*(t)} & i^*E)
}$$
and $\beta$ is equivalent to the composite
$$\xymatrix{
i^* E \simeq i^* \Delta^* q_2^* E \ar[r]_-{\sim} & i^* \Delta^* q_2^* E \simeq i^*E
}$$
with the middle morphism above induced by the pullback diagram
$$\xymatrix{
\mcal{L}X \ar[r]^-{i} \ar[d]_-{i}&  X \ar[d]^-{\Delta}
\\
X \ar[r]_-{\Delta} & X \times X \ar[r]^-{q_2} & X.
}$$
Rewriting this diagram as 
$$\xymatrix{
\mcal{L}X \ar@<-.5ex>[r]_-{i} \ar@<.5ex>[r]^-{i} & X \ar[r]^-{\Delta} & X \times X \ar[r]^-{q_2} & X
}$$
we see that $\beta$ is precisely the action morphism of $\mcal{L}X$ on $\Delta^* q_2^* E$, that is, an action morphism of  $\mcal{L}X$ on $E$ with the canonical $\mcal{L}X$-equivariant structure.
\end{proof}
\end{Prop}

\subsection{Categorical Chern character as exponential}
In this section we give a different description of $\mathcal LX$-equivariant structure on a sheaf $E \in \QCoh(X)$ on a smooth proper scheme $X$ using the formal deformation theory developed in \cite{GaitsRozII} (we invite the reader to look at Appendix \ref{app:formal_moduli_reminder} for a quick reminder). Recall that in order to conveniently work with deformation theory one needs to replace the category of quasi-coherent sheaves $\QCoh(X)$ with a closely related but different category of ind-coherent sheaves $\ICoh(X)$ (see section \ref{sect:ICoh_reminder}). However, for a smooth scheme $X$ there is a natural symmetric monoidal equivalence $\xymatrix{\QCoh(X) \ar[r]^-{\Upsilon_X}_-{\sim} & \ICoh(X)}$ (see Example \ref{icoh_of_smooth}), so we will state all needed results of \cite{GaitsRozII} using $\QCoh(X)$ instead of $\ICoh(X)$.

\smallskip  As we have seen in Proposition \ref{23action} the Chern character can be described using the group structure (over $X$) on the inertia group $\mathcal LX$ of $X$. Note that the underlying quasi-coherent sheaf of the Lie algebra $\Lie_X(\mathcal LX)$ that corresponds to $\mathcal{L}X$ via formal groups-Lie algebras correspondence (see Theorem \ref{LieAlgAndFormGroups}) is easy to understand: because of the pullback square
$$\xymatrix{
\mathcal LX \ar[d]^i\ar[r]^i & X \ar[d]^-\Delta
\\
X \ar@/^/[u]^e \ar[r]_-\Delta & X\times X
}$$\smallskip
 where $e$ is obtained by pulling back $q_1$ we get an equivalence
$$
e^*\mathbb T_{\mathcal LX/X} \simeq e^*i^*\mathbb T_{X/X\times X}\simeq \mathbb T_{X/X\times X} \simeq \mathbb T_X[-1]
$$
in $\QCoh(X)$.

\begin{Cor}[Hochschild-Kostant-Rosenberg]\label{HochKosRos}
For a smooth scheme $X$ we have
$$i_* \mathcal O_{\mathcal LX} \simeq \Sym_{\QCoh(X)}(\Omega_X^1[1]) \simeq \bigoplus_{p=0}^{\dim X} \Omega_{X}^p[p] \qquad\qquad HH_i(X)\simeq \bigoplus_{p-q=i}^{\dim X} H^{q}(X, \Omega_X^p).$$
\begin{proof}
Since by definition $HH_i(X) := \pi_i\Gamma(\mathcal LX, \mathcal O_{\mathcal LX}) \simeq \pi_i\Gamma(X, i_*\mathcal O_{\mathcal LX})$ the second equivalence immediately follows from the first one. To get the first equivalence, note that applying the exponent map (see Theorem \ref{formal_exponent}) we get an equivalence
$$\xymatrix{
\mathbb{V}(\mathbb T_X[-1]) \ar[rr]^-{\exp_{\mcal{L}X}}_-{\sim} && \mcal{L}X
}$$
of formal schemes over $X$. But by smoothness of $X$ we have $(\mathbb T_X[-1])^\vee \simeq \Omega_X^1[1] \in \Coh^{<0}(X)$ and hence by Example \ref{ex:vector_stacks_vs_rel_spec} we get
$$\mathbb V(\mathbb T_X[-1]) \simeq \Spec_{/X}\left(\Sym_{\QCoh(X)}(\Omega_X[1])\right)$$
thus obtaining an equivalence $i_* \mcal{O}_{\mcal{L}X} \simeq i_*\mathcal O_{\mathbb V(\mathbb T_X[-1])} \simeq \Sym_{\QCoh(X)}(\Omega_X[1])$ in $\QCoh(X)$.
\end{proof}
\end{Cor}

\medskip

Now the Hochschild-Kostant-Rosenberg theorem above shows us that for a perfect sheaf with an endomorphism $(E,t)$ the categorical Chern character $\ch(E, t)$ is in fact an element of
$$
\pi_0 \Gamma(\mathcal LX, \mathcal O_{\mathcal LX}) \simeq \bigoplus_{p=0}^{\dim X} H^p(X, \Omega^p_X).
$$
In particular, it makes sense to ask for a more concrete description of $\ch(E,t)$. In order to do this, let us pick a sheaf $E\in \QCoh(X)$ and analyze the canonical $\mcal{L}X$-equivariant structure on $E$. Since action of a formal group correspond to the action of the corresponding Lie algebra  (see Theorem \ref{LieAlgAndFormGroups}) in the case when our formal group is the inertia group $\mcal{L}X$ we get an equivalence
$$\xymatrix{
\Rep_{\mcal{L}X}\big{(} \QCoh(X) \big{)} \simeq \Mod_{\mathbb T_X [-1]}\left(\QCoh(X)\right) .
}$$
Consequently, we see that the canonical $\mcal{L}X$-equivariant structure on $E\in \QCoh(X)$ can be equivalently understood in terms of the corresponding action of $\mathbb T_X [-1] \in \LAlg(\QCoh(X))$ on $E$ (which we will further also call \emph{canonical}). 

In order to understand the canonical action of $\bb{T}_X[-1]$, let us for a second discuss a more general situation. Suppose that we have a Lie algebra $\mathfrak g \in \QCoh(X)$ which acts on some $\mcal{F} \in \QCoh(X)$. The associative algebra structure on $\EEnd_{\QCoh(X)}(\mcal{F}) \in \Alg(\QCoh(X))$ endows it with the structure of a Lie algebra (we will further denote by $\gl_{\mcal{F}} \in \LAlg(\QCoh(X))$ the corresponding Lie algebra) and the action of $\mathfrak g$ on $\mcal{F}$ induces a map $\xymatrix{\mathfrak g \ar[r] & \mathfrak{gl}_{\mcal{F}}}$ in $\LAlg(\QCoh(X))$. Now by Lie algebras-formal groups correspondence from \ref{LieAlgAndFormGroups} we can ``integrate'' this map to obtain a map $\xymatrix{\widehat G  \ar[r]^-{\rho} & \GL_{\mcal{F}}}$ of formal groups, where  $\Lie_X(\widehat G) \simeq \mathfrak g$ and $\Lie_X(\widehat{\GL}_{\mathcal{F}}) \simeq \gl_{\mathcal{F}}$. It particular, by functoriality this gives a commutative diagram
$$\xymatrix{
\mathbb V(\mathfrak{g}) \ar[rr]^-{\mathbb V(\Lie_X(\rho))}\ar[d]_{\exp_{\widehat G}}^{\sim} & & \mathbb V(\gl_{\mcal{F}}) \ar[d]_{\sim}^{\exp_{\widehat{\GL}_{\mcal{F}}}}
 \\
\widehat G \ar[rr]_-\rho & & \widehat \GL_{\mcal{F}}
}$$
of formal moduli problems over $X$.

Now applying the above procedure to the canonical action of $\mathbb T_X[-1]$ on $E$ we arrive at the following

\begin{Def}\label{AtiyahClass}
Let $X$ be a smooth scheme. Define \emph{Atiyah class $\At(E)$} of $E \in \QCoh(X)$ as the top horizontal map in the commutative diagram
$$\xymatrix{
\mathbb V\big{(}\mathbb T_X[-1])\big{)} \ar[rr]^-{\At(E)} \ar[d]_-{\exp_{\mathcal LX}}^-\sim && \mathbb V(\gl_{E})  \ar[d]_{\sim}^{\exp_{\widehat{\GL}_E}}
\\
\mathcal LX \ar[rr]_-{\alpha_E} && \widehat \GL_{E}.
}$$

By the Lie algebras--formal groups correspondence Atiyah class of $E$ corresponds by definition to the canonical $\mathbb T_X[-1]$-module structure on $E$. By taking dual of the Lie-module structure map $\xymatrix{\mathbb T_X[-1] \ar[r] & \EEnd_{\QCoh(X)}(E)}$ we see that $\At(E)$ corresponds to some class in $H^1(X, \EEnd_{\QCoh(X)}(E)\otimes\mathbb L_X)$. 
\end{Def}

Combining Proposition \ref{23action} with the notation above we obtain:
\begin{Cor}\label{prop:semi_eq_Chern}
Let $X$ be a smooth proper scheme and $E \in \QCoh(X)$ is a perfect sheaf with an endomorphism $\xymatrix{E \ar[r]^-{t} & E}$. Then under the Hochschild-Kostant-Rosenberg identification we have an equality
$$\ch(E,t) = \Tr_{\QCoh(\mathcal LX)}\left(\xymatrix{i^*E \ar[rr]_\sim^{\exp(\At(E))} && i^*E \ar[rr]^-{i^*(t)} && i^*E}\right)$$
of elements of $\displaystyle \bigoplus_p H^p(X, \Omega^p_X)$.
\end{Cor}

\begin{Ex} \label{ex:ch0}
Recall that for a perfect sheaf $E \in \QCoh(X)$ we have an equality
$$\ch_0(E) = \mathrm{rk}(E) = \Tr_{\QCoh(X)}(\Id_E)$$
where $\ch(E)$ here is the classical Chern character. Now as we will see below (Proposition \ref{ChernIsChern}), the categorical Chern character $\ch(E,\Id_E)$ in fact coincides with the classical one $\ch(E)$. In particular, the corollary above allows as to prove a generalization of this statement: for a perfect sheaf $(E,t)$ with an endomorphism we get an equality
$$\ch(E,t) =\Tr\left(t + t\circ \At(E) + t \circ \frac{\At(E)^{\wedge 2}}{2} + \ldots\right)$$
and therefore since $\Tr(\At(E)^{\wedge n}) \in H^n(X, \Omega^n_X)$ we obtain $\ch_0(E, t) = \Tr(t)$.
\end{Ex}

\subsection{Comparison with the classical Chern character}
We will now compare our definition of Atiyah class with a more classical one. These results are well-known to experts and we include them mostly for reader's convenience.
\begin{Def}
The prestack $\Perf$ of \emph{perfect sheaves} is defined as
$$
\Perf(R) := (\Mod_R^{\perf})^{\simeq}.
$$
\end{Def}

Note that since the operation of taking full subcategory of dualizable objects commutes with limits of symmetric monoidal categories and the operation of taking maximal subgroupoid commutes with limits of categories for any prestack $Y \in \PreStack$ we get a natural equivalence
$$
\Hom_{\PreStack}(Y,\Perf) \simeq \lim \limits_{\Spec R \in \Aff_{/Y}}\Hom_{\PreStack}(\Spec R, \Perf) =\lim \limits_{\Spec R \in \Aff_{/Y}}(\Mod_R^{\perf})^{\simeq} \simeq \big{(} \QCoh(Y)^{\perf} \big{)}^{\simeq}
$$
Let us now denote by $\mcal{E} \in \QCoh(\Perf)$ the universal perfect sheaf classified by the identity morphism $\Id_{\Perf}$. We then have
\begin{Prop}[{\cite[Chapter 8, Proposition 3.3.4.]{GaitsRozII}}]
There is a canonical equivalence
$$\xymatrix{
\mathbb T_{\Perf}[-1] \ar[r]_-{\sim} & \EEnd_{\QCoh(\Perf)}(\mathcal E)
}$$
of Lie algebras.
\end{Prop}

\begin{Rem} It is not hard to show that the underlying sheaves of $\bb{L}_{\Perf}$ and $\EEnd_{\QCoh(\Perf)}(\mathcal E)[-1]$ are equivalent. Indeed, the pullback diagram
$$\xymatrix{
\mathcal{L}\Perf \ar[r]^-{i} \ar[d]^-{i} & \Perf \ar[d]^-{\Delta}
\\
\Perf \ar[r]_-{\Delta} \ar@/^1.0pc/@{-->}[u]^-{e} & \Perf \times \Perf
}$$
induces an equivalence
$$
\bb{L}_{\Perf} \simeq \bb{L}_{\Perf/\Perf \times \Perf}[-1] \simeq e^* i^* \bb{L}_{\Perf/\Perf \times \Perf}[-1]  \simeq e^* \bb{L}_{\mathcal{L}\Perf/\Perf}[-1].
$$
Now to calculate $e^* \bb{L}_{\mathcal{L}\Perf/\Perf}[-1]$ we note that by definition for any morphism
$$\xymatrix{
\Spec R \ar[r]^-{\eta} & \Perf
}$$
classifying $E \in \Mod_R^{\perf}$ and an $R$-module $M \in \Mod_{R}^{\leq 0}$ the space $\Hom_{\Mod_{R}}(\eta^* e^*  \bb{L}_{\mathcal{L}\Perf/\Perf},M)$ is equivalent to the space of lifts
$$\xymatrix{
\Spec R \ar[d]_-{t} \ar[r]^-{\eta} & \Perf \ar[r]^-{e} & \mathcal{L}\Perf \ar[d]^-{i}
\\
\Spec R[M] \ar@{-->}[urr] \ar[rr]_-{\zeta} && \Perf.
}$$
Now the morphism $\zeta$ classifies the module $u^* E$ where $\xymatrix{\Spec R[M] \ar[r]^-{u} & \Spec R}$ is the projection map (so that $t^* u^* \simeq \Id_{\Mod_{R}}$) and $\mathcal{L}X$ classifies a perfect sheaf together with an automorphism we get
\begin{align*}
& \Hom_{\ssl\Mod_{R}}(\eta^* e^*  \bb{L}_{\mathcal{L}\Perf/\Perf},M) \simeq \Aut_{\ssl\Mod_{R[M]}}(u^*E) \underset{\Aut_{\Mod_R}(E)}{\times} \{\Id_{E}\} \simeq\\
& \simeq \Hom_{\ssl\Mod_{R[M]}}(u^* E,u^*E)\underset{\ssl\Aut_{\Mod_R}(E)}{\times} \{\Id_{E}\} \simeq \Hom_{\ssl\Mod_{R}}(E,u_* u^* E) \underset{\ssl\Aut_{\Mod_R}(E)}{\times} \{\Id_{E}\} \simeq\\
& \simeq \Hom_{\ssl\Mod_{R}}(E,E \oplus E \otimes_R M) \underset{\ssl\Aut_{\Mod_R}(E)}{\times} \{\Id_{E}\} \simeq \Hom_{\ssl\Mod_{R}}(E,E \otimes_R M) \simeq \\
& \simeq \Hom_{\ssl\Mod_{R}}\big{(} \EEnd_{\ssl\Mod_{R}}(E),M \big{)} \simeq \Hom_{\ssl\Mod_{R}}\big{(} \eta^*\EEnd_{\ssl\QCoh(\Perf)}(\mcal{E}),M \big{)} 
\end{align*}
proving the claim, where above we use that any endomorphism of $u^*E \in \Mod_{R[M]}$ lying over $\Id_E$ is automatically invertible.
\end{Rem}

The following proposition provides a convenient description of the Atiyah class as a map of underlying sheaves:
\begin{Prop}\label{prop:more_classical_Atiyah}
Let $E \in \Perf(X)$ be a perfect sheaf on a smooth scheme $X$ classified by the map $\xymatrix{X \ar[r]^-{e} & \Perf}$. Then the induced map of tangent spaces
$$\xymatrix{
\mathbb T_X \ar[r] & e^*\mathbb T_{\Perf} \simeq \EEnd_{\QCoh(X)}(E)[1]
}$$
is equal to $\At(E)[1]$.


\begin{proof}
Note that given two stacks $\mathcal X, \mathcal Y$ admitting deformation theory together with formal groups $H \in \Grp(\FormModuli_{\mathcal{X}})$ and $G\in \Grp(\FormModuli_{\mathcal{Y}})$ over $\mathcal X$ and $\mathcal Y$ respectively the datum of commutative square
$$\xymatrix{
\mathcal X \ar[r]^-{f} \ar[d] & \mathcal Y\ar[d] 
\\
\widehat B_{/\mathcal X} H \ar[r] & \widehat B_{/\mathcal Y} G
}$$
determines a morphism
$$\xymatrix{
H \ar[r] & G_{\mathcal X}:= \mathcal X\times_{\mathcal Y}G
}$$
of formal groups over $\mcal{X}$. Unwinding the definitions, one finds that under the identifications
$$
\Lie_{\mathcal X}(H) \simeq \mathbb T_{\mathcal X/(\widehat B_{/\mathcal X} H)} \qquad \Lie_{\mathcal X}(G_{\mathcal X}) \simeq f^*\Lie_{\mathcal Y}(G) \simeq f^*\mathbb T_{\mathcal Y/(\widehat B_{/\mathcal Y} G)}
$$
the induced map of Lie algebras
$$\xymatrix{
\mathbb T_{\mathcal X/(\widehat B_{/\mathcal X} H)} \simeq \Lie_X(H) \ar[r] & f^*\Lie_Y(G) \simeq f^*\mathbb T_{\mathcal Y/(\widehat B_{/\mathcal Y} G)}
}$$
coincides with the natural map of relative tangent sheaves. Applying this observation to the diagram
\begin{align}\label{BLX_BLPerf_relation}
\xymatrix{
X\ar[r]^-{e} \ar[d]_-{\Delta} & \Perf\ar[d]^-{\Delta}
\\
(X\times X)_{\widehat\Delta} \ar[r]_-{e\times e} & (\Perf\times\Perf)_{\widehat\Delta}
}
\end{align}
and using the equivalence $\mathbb{T}_X \simeq \mathbb{T}_{X/(X\times X)_{\widehat\Delta}}[1]$ (and similarly for $\Perf$) we deduce that the map of tangent sheaves $\xymatrix{\mathbb T_X \ar[r] & e^*\mathbb T_{\Perf} \simeq \EEnd_{\QCoh(X)}(E)[1]}$ we are interested in is the shift of the morphism that underlies the map of Lie algebras induced by the morphism $\xymatrix{\mathcal{L}X \ar[r] & X \times_{\Perf}\widehat{\mathcal L}\Perf}$ of formal groups over $X$, where $\widehat{\mathcal L}\Perf \in \Grp(\FormModuli_{\Perf})$ is the completion of $\mathcal L\Perf$ along the constant loops. Since $\mathcal L\Perf$ classifies a perfect sheaf with an automorphism we get an equivalence $X \times_{\Perf}\widehat{\mathcal L}\Perf \simeq \widehat{\GL}_{/X}(E)$ of formal groups over $X$. Moreover, by commutativity of the diagram (\ref{BLX_BLPerf_relation}), the induced group morphism $\xymatrix{\mathcal LX \ar[r]^-{\alpha} & \widehat{\GL}_{/X}(E)}$ is given precisely by the canonical action. The result now follows from the definition \ref{AtiyahClass} of the Atiyah class.
\end{proof}
\end{Prop}

\smallskip
Recall now the following 

\begin{construction}[Classical algebraic Chern character]\label{classical_Chern__precise}
Let $X$ be a smooth proper scheme and $\mathcal L$ be a line bundle on $X$. We define
$$\ch(\mathcal L) := \exp(c_1(\mathcal L))\quad \in \quad \bigoplus_p H^p(X, \Omega_X^p),$$
where $c_1(\mathcal L) \in H^1(X, \Omega_X^1)$ is the algebraic (Hodge version of the) first Chern class defined via the map $\xymatrix{\mathcal O_X^* \ar[r]^-{d\log} & \Omega_X^1}$. For a general perfect sheaf $E$ on $X$ we define $\ch(E) \in \bigoplus_p H^p(X, \Omega_X^p)$ by additivity and the splitting principle.
\end{construction}

\begin{warning} \label{topological_vs_algebraic_ch}
Let $k=\mathbb C$, $X$ a smooth proper scheme over $k$ and $E$ a vector bundle on $X$. Then one has a priori two different notions of the Chern character: the one constructed above and the topological Chern character $\xymatrix{\ch^{\mathrm{top}}(E) \in H^*(X(\mathbb C), \mathbb Z) \ar@{^{(}->}[r] & H^*(X(\mathbb C), \mathbb C) \simeq \bigoplus_{p,q} H^{q,p}(X)}$. It turns out these two notions \textit{do not} coincide, but are very closely related to each other: a simple computation on $\mathbb P^1$ shows that $c_1^{\mathrm{top}}(\mathcal O_{\mathbb P^1}(1)) = -2\pi i \cdot c_1(\mathcal O_{\mathbb P^1}(1))$ and hence in general $\ch_k^{\mathrm{top}}(E) = (-2\pi i)^k \ch_k(E)$.
\end{warning}

In particular, as a well-known corollary we recover the fact that the classical Chern character of a sheaf can be expressed in terms of the Atiyah class:
\begin{Cor}\label{AtiyaClassOfLineBundleIsChern}
Let $X$ be a smooth proper scheme.
\begin{enumerate}
\item Let $\mathcal M$ be a line bundle on $X$. Then the Atiyah class $\At(\mathcal M) \in H^1(X, \EEnd_{\QCoh(X)}(\mathcal M)\otimes\Omega_X^1) \simeq H^1(X, \Omega_X^1)$ of $\mathcal M$ coincides with the first Chern class $c_1(\mathcal M)$ of $\mathcal M$.

\item Let $E$ be a perfect sheaf on $X$. Then the classical Chern character of $E$ is equal to $\Tr \exp(\At(E))$.
\end{enumerate}

\begin{proof}
For the first claim, note that by Proposition \ref{prop:more_classical_Atiyah} the Atiyah class of any perfect sheaf $E \in \QCoh(X)^{\perf}$ can be equivalently described as the shift of the induced map on tangent spaces
$$\xymatrix{
\mathbb T_X \ar[r] & e^* \mathbb T_{\Perf} \simeq \EEnd_{\QCoh(X)}(E)[1]},
$$
where the map $\xymatrix{X \ar[r]^-{e} & \Perf}$ classifies $E \in \QCoh(X)^{\perf}$. Now if $E=\mathcal M$ is a line bundle, the classifying map $\xymatrix{X \ar[r] & \Perf}$ factors through $B\mathbb G_m$ and so it is enough to prove the statement for the universal line bundle $\mathfrak L$ on $B\mathbb{G}_m$. But the canonical map $\xymatrix{ B\mathbb G_m \ar[r]^-{u} & \Perf}$ induces an equivalence $\mathbb T_{B\mathbb G_m} \simeq u^*\mathbb T_{\Perf}$ and hence the Atiyah class obtained as the shift of
$$\xymatrix{
\At(\mathfrak L)\colon \mathcal O_{B\mathbb G_m}[1] \simeq \mathbb T_{B\mathbb G_m} \ar[r] & \EEnd_{\QCoh(B\mathbb G_m)}(\mathfrak L)[1] \simeq \mathcal O_{B\mathbb G_m}[1]
}$$
is just the identity map, which corresponds to $c_1(\mathfrak L)$.

For the second claim, as both $\ch(E)$ and $\Tr(\exp(\At(E)))$ are additive in triangles and commute with pullbacks by the splitting principle it is sufficient to prove the equality for $E = \mathcal M$ being a line bundle. But since by the previous part $\At(\mathcal M) = c_1(\mathcal M)$, we get
$$\Tr(\exp(\At(\mathcal M)))=\exp(\At(\mathcal M))=\exp(c_1(\mathcal M))=\ch(\mathcal M)$$
as claimed.
\end{proof}
\end{Cor}
\begin{Rem}[Chern-Weil theory and the Atiyah class]\label{ChernWeild_and_Atiyah}
Since $\ch(E,t)$ is $k$-linear in the second argument, one expects it to be an element of some $k$-linear cohomology theory of $X$ like de Rham or Hodge cohomology. Our description of $\ch(E,t)$ is closer to the differential-geometric approach to characteristic classes. Namely, recall that if $X$ is a smooth proper scheme over the field $\mathbb C$ of complex numbers and $E$ is a vector bundle over $X$ one can give the following description of the algebraic Chern character of $E$ (see Warning \ref{topological_vs_algebraic_ch} for the difference between algebraic and topological Chern characters): choose a smooth connection $\nabla$ on $E$ and let $F_\nabla \in \EEnd(E)\otimes \Omega_X^2$ be the corresponding curvature form. Then
$$\ch(E) = \Tr e^{F_\nabla} \in H^*_{\mathrm{dR}}(X, \mathbb C).$$
Now assume additionally that $\nabla$ is of type $(1,0)$. Then the curvature $F_{\nabla}$ splits into a sum $F^{2,0}_{\nabla} + F^{1,1}_{\nabla}$. It follows that $\Tr e^{F^{1,1}_{\nabla}}$ is a representative of $\ch(E)$ in Hodge cohomology $\bigoplus_p H^{p,p}(X) \simeq \bigoplus_p F^p H_{\mathrm{dR}}^p(X, \mathbb C)/F^{p+1} H_{\mathrm{dR}}^p(X, \mathbb C)$. Finally, it turns out the class $F^{1,1}_\nabla \in H^1(X, \Omega_X^1 \otimes \EEnd(E))$ is independent of a choice of $\nabla$ and in fact coincides with the Atiyah class $\At(E)$ (see \cite[Proposition 4]{Atiyah_AtiyahClass}), linking the previous corollary and the Chern-Weil theory.
\end{Rem}

Finally we obtain a concrete description of the categorical Chern character of a sheaf with the trivial equivariant structure
\begin{Prop}\label{ChernIsChern}
Let $E$ be a dualizable object of $\QCoh(X)$. Then under the Hochschild-Kostant-Rosenberg isomorphism \ref{HochKosRos}
$$\pi_0 \Tr_{2\Cat_k}(\Id_{\QCoh(X)}) \simeq \pi_0\Gamma(\mathcal LX, \mathcal O_{\mathcal LX}) \xymatrix{\ar[r]^-{\mathrm{HKR}}_-\sim &} \bigoplus_{p=0}^{\dim X} H^p(X, \Omega_X^p)$$ 
the categorical Chern character $\ch(E, \Id_E)$ coincides with the classical one $\ch(E)$.

\begin{proof}
By Proposition \ref{23action} we have $\ch(E, \Id_E) = \Tr_{\QCoh(\mcal{L}X)}(\alpha_E)$ and by definition of the Atiyah class we have $\exp_{\mathcal LX}^*(\alpha_E) = \exp(\At(E))$. Since by construction $\mathrm{HKR}=\exp_{\mathcal LX}^*$, we conclude by the second part of Corollary \ref{AtiyaClassOfLineBundleIsChern}.
\end{proof}
\end{Prop}

\section{Trace of pushforward functor via ind-coherent sheaves}\label{sect:icoh}
Recall from \cite[Chapter 5, 5.3]{GaitsRozI} the quasi-coherent sheaves functor in fact has an appropriate $2$-categorical functoriality, in a sense that it can be lifted to a symmetric monoidal functor from a symmetric monoidal $(\infty,2)$-category of correspondences (see Appendix \ref{app:corr} for a discussion of traces and correspondences). This allows us to reformulate many $2$-categorical questions about quasi-coherent sheaves to questions about the category of correspondences, where they can be in most cases answered by direct diagram chasing. This observation, for example, gives a direct proof that the morphism of traces $\Tr(f^*,\Id_{f^*})$ induced by the diagram
$$\xymatrix{
\QCoh(X) \ar@/_/[dd]_-{f^*} \ar[rr]^-{\Id_{\QCoh(X)}} && \QCoh(X)\ar@/_/[dd]_-{f^*}
\\
\\
\QCoh(Y) \ar@/_/[uu]_-{f_*} \ar[rr]_-{\Id_{\QCoh(Y)}} && \QCoh(Y) \ar@/_/[uu]_-{f_*}
}$$
in $\Cat_k$ coincides with the classical pullback of global sections (see Remark \ref{Pullback_in_cohom}). However, the same argument does not apply to the morphism of traces
$$\xymatrix{
\displaystyle \bigoplus_{p=0}^{\dim X} H^p(X, \Omega_X^p)  \simeq \Tr(\Id_{\QCoh(Y)}) \ar[rrr]^-{\Tr(f_*,\Id_{f_*})} &&& \Tr(\Id_{\QCoh(X)}) \simeq \displaystyle  \bigoplus_{p=0}^{\dim Y} H^p(Y, \Omega_Y^p)
}$$
we are interested in. Indeed, due to the post factum knowledge that the answer should involve the Todd class, one should not expect to obtain a concrete description of it in a purely formal way. 

But as is mentioned in previously, apart from $\QCoh(X)$ there is another important $\ii$-category we can associate to $X$: the $\ii$-category $\ICoh(X) \in 2\Cat_k$ of ind-coherent sheaves on $X$. As a toy example, for a smooth classical scheme $X$ there is a functor $\xymatrix{\ICoh(X) \ar[r]^-{\Psi_X}_-{\sim} & \QCoh(X)}$ which identifies the category of ind-coherent sheaves with the category of quasi-coherent sheaves (see Example \ref{icoh_of_smooth}) as a plain category, but does not preserve monoidal structure: via this equivalence the natural $\overset{!}{\otimes}$-monoidal structure on $\ICoh(X)$ is given by
$$\mathcal F \overset{!}{\otimes} \mathcal G \simeq \mathcal F\otimes\mathcal G \otimes \omega_X^{-1}.$$
where $\omega_X \in \QCoh(X)$ is the dualizing sheaf.

In this section we will prove that the morphism of traces induced by pushforward could be understood relatively easy in the setting of ind-coherent sheaves:

\begin{itemize}
\item In subsection \ref{sect:ICoh_reminder} we will review relevant facts about the category of ind-coherent sheaves.

\item Similar to that of $\QCoh$, using (Serre) self-duality of $\ICoh$ we will prove that there is an equivalence
$$
\Tr_{2\Cat_k}(g_*) \simeq \Gamma(X^g, \omega^{\ICoh}_{X^g}).
$$

\item Then we will show that the morphism of traces 
$$\xymatrix{\Gamma(X^g, \omega_{X^g}^{\ICoh}) \simeq \Tr_{2\Cat_k}(g_{X*})\ar[rr]^-{\Tr_{2\Cat_k}(f_*)} && \Tr_{2\Cat_k}(g_{Y*}) \simeq \Gamma(Y^g, \omega_{Y^g}^{\ICoh})}$$
induced by the diagram
$$\xymatrix{
\ICoh(X) \ar[d]_-{f_*} \ar[rr]^-{(g_X)_*} && \ICoh(X) \ar[d]^-{f_*}
\\
\ICoh(Y) \ar[rr]_-{(g_Y)_*} && \ICoh(Y)
}$$
in $2\Cat_k$ coincides with the natural pushforward of distributions. In particular if both $X,Y$ are smooth and proper with trivial equivariant structure, from this we will deduce that the morphism of traces
$$\xymatrix{
\displaystyle\left(\bigoplus_{p=0}^{\dim X} H^p(X, \Omega_X^p)\right)^{\vee}  \simeq \Tr_{2\Cat_k}(\Id_{\ICoh(X)}) \ar[rr]^-{\Tr_{2\Cat_k}(f_*)} && \Tr(\Id_{\ICoh(Y)}) \simeq \displaystyle \left( \bigoplus_{p=0}^{\dim Y} H^p(Y, \Omega_Y^p) \right)^\vee
}$$
coincides with the pushforward in homology.
\end{itemize}

After dealing with this, we will use further sections to investigate how one can describe the relation between the morphism of traces induced by pushforward in the setting of $\ICoh$ and in the setting of $\QCoh$.

\subsection{Reminder on Ind-coherent sheaves}\label{sect:ICoh_reminder}
In this subsection we review some basic facts and constructions related to ind-coherent sheaves. We refer reader to \cite[Part II]{GaitsRozI} and \cite{Gaits} for further details. We start with the following
\begin{Def}\label{icoh_def}
For $X \in \Sch_{\aft}$ (see \cite[Chapter 4, 1.1.1]{GaitsRozI}) define the category of {\bfseries ind-coherent sheaves on $X$} denoted by $\ICoh(X)$ simply as
$$
\ICoh(X):=\Ind(\Coh(X)),
$$
where we denote by $\Coh(X)$ the category of coherent sheaves on $X$, i.e. the full subcategory of $\QCoh(X)$ consisting of those $\mathcal F \in \QCoh(X)$ such that $\mathcal H^i(\mathcal F)$ are non-zero only for finitely many $i$ and are coherent over the sheaf of algebras $\mathcal H^0(\mathcal O_X)$ in the usual sense. 
\end{Def}

Properties of the ind-coherent sheaves construction we need in this paper can be summarized by the following
\begin{Prop}\label{icoh}\
\\
1) (\cite[Chapter 4, Proposition 2.1.2, Proposition 2.2.3]{GaitsRozI}) The assignment of ind-coherent sheaves can be lifted to a functor
$$\xymatrix{
\Sch_{\aft} \ar[rr]^-{\ICoh_*} && \Cat_k
}$$
such that, moreover, for every morphism $\xymatrix{X \ar[r]^-{f} & Y}$ in $\Sch_{\aft}$ the diagram
$$\xymatrix{
\ICoh(X) \ar[r]^-{\Psi_X}\ar[d]_-{f_*} & \QCoh(X) \ar[d]^-{f_*}
 \\
\ICoh(Y) \ar[r]_-{\Psi_Y} & \QCoh(Y)
}$$
commutes, where $\xymatrix{\ICoh(X) \ar[r]^-{\Psi_X} & \QCoh(X)}$ is obtained by ind-extending the natural inclusion $\Coh(X) \subseteq \QCoh(X)$ (and similar for $Y$).
\\
\\
2) (\cite[Chapter 4, Corollary 5.1.12]{GaitsRozI}) The assignment of ind-coherent sheaves can be lifted to a functor
$$\xymatrix{
\Sch_{\aft,\sf proper}^{\op} \ar[rr]^-{\ICoh^!} && \Cat_k,
}$$
such that, moreover, given a proper morphism $\xymatrix{X \ar[r]^-{f} & Y}$ in $\Sch_{\aft}$ the induced pullback functor $f^!:=\ICoh^!(f)$ is right adjoint to $f_*$.
\\
\\
3) (\cite[Chapter 4, Proposition 6.3.7; Chapter 5, Theorem 4.2.5]{GaitsRozI}) For every $X \in \Sch_{\aft}$ the category $\ICoh(X)$ is symmetric monoidal and self-dual as an object of $\Cat_k$ (see Theorem \ref{IndCohSelfDual} below for a concrete description of duality maps). Moreover, for every proper $\xymatrix{X \ar[r]^-{f} & Y}$ the induced functor $f^!$ is symmetric monoidal. We will further denote the monoidal structure on $\ICoh(X)$ by $\sotimes$ and the monoidal unit, the so-called \emph{$\ICoh$-dualizing sheaf}, by $\omega_X^{\ICoh} \in \ICoh(X)$. It is straightforward to see that there is in fact an equivalence $\omega_X^{\ICoh} \simeq p^! k$, where $\xymatrix{X \ar[r]^-{p} & \ast}$ is the projection and $k \in \ICoh(\ast) \simeq \Vect_k$.
\\
\\
4) (\cite[Chapter 6, 0.3.5, 3.2.5]{GaitsRozI}) The functor $\xymatrix{\QCoh(X) \ar[r]^-{\Upsilon_X} & \ICoh(X)}$ obtained from $\Psi_X$ using self-dualities of $\QCoh(X)$ and $\ICoh(X)$ is symmetric monoidal and for every proper morphism $\xymatrix{X \ar[r]^-{f} & Y}$ in $\Sch_{\aft}$ the diagram
$$\xymatrix{
\QCoh(X) \ar[rr]^-{\Upsilon_X} && \ICoh(X) 
\\
\QCoh(Y) \ar[rr]_-{\Upsilon_Y}\ar[u]^-{f^*} && \ICoh(Y) \ar[u]_-{f^!}
}$$
commutes. 
\end{Prop}

In our main case of interest the categories of quasi-coherent sheaves and ind-coherent sheaves are not that far away to each other. Recall first the following 

\begin{construction}\label{const:qcoh_omega}
Since $X \in \Sch_{\aft}$ is quasi-compact, the global sections functor
$$\xymatrix{\QCoh(X) \ar[rr]^-{p_*} && \Vect_k}$$
is continuous, hence admits a right adjoint $p^!$. We define \emph{$\QCoh$-dualizing sheaf $\omega_X \in \QCoh(X)$} by setting $\omega_X := p^!(k)$.
\end{construction}

\begin{Ex}
If $X$ is smooth, then by classical Serre duality there is an equivalence
$$\omega_X \simeq \Omega_X^{\dim X}[\dim X]$$
in $\QCoh(X)$. In particular, $\omega_X$ is dualizable in this case.
\end{Ex}

\begin{Ex}\label{icoh_of_smooth}
Let $X$ be a smooth classical scheme. Since $X$ is quasi-compact and separated there is an equivalence $\QCoh(X) \simeq \Ind(\QCoh(X)^{\perf})$ of categories (see e.g. \cite[Proposition 3.19]{BFN}), and by smoothness we also get $\Coh(X) \simeq \QCoh(X)^{\perf}$. Consequently, it follows that the canonical functor $\xymatrix{\ICoh(X) \ar[r]^-{\Psi_X} & \QCoh(X)}$ is an equivalence of categories and thus the functor $\xymatrix{\QCoh(X) \ar[r]^-{\Upsilon_X} & \ICoh(X)}$ is a symmetric monoidal equivalence of categories. In particular we can identify $\ICoh(X)$ with $\QCoh(X)$ with the twisted monoidal structure
$$\mathcal F \sotimes G \simeq \mathcal F \otimes \mathcal G \otimes \omega_X^{-1}.$$
\end{Ex}

\medskip

We will further need comparison between $\omega_X$ and $\omega_X^{\ICoh}$:
\begin{Prop}\label{compare_dualizing}
Let $X$ be a proper derived scheme. Then there is an equivalence $\omega_X \simeq \Psi_X(\omega_X^{\ICoh})$ in $\QCoh(X)$. In particular $\Gamma(X, \omega^{\ICoh}_X) \simeq \Gamma(X, \omega_X)$ where here $\Gamma=p_*$ is the pushforward along projection morphism $\xymatrix{X \ar[r]^-{p} & \ast}$ in the setting of ind-coherent and quasi-coherent sheaves respectively.

\begin{proof}
The first statement follows from the fact that due to \cite[Proposition 7.2.2., Proposition 7.2.9(a)]{Gaits} the diagram
$$\xymatrix{
\ICoh(X) \ar[r]^-{\Psi_X} & \QCoh(X)
\\
\Vect_k \ar[u]^-{p^{!,\ICoh}} \ar[r]^-{\Psi_{\ast}}_-{\sim} & \Vect_k \ar[u]_-{p^{!,\QCoh}}
}$$
commutes. The assertion about global sections follows from commutativity of the square
$$\xymatrix{
\ICoh(X) \ar[r]^-{\Psi_X}\ar[d]_-{p_*} & \QCoh(X) \ar[d]^-{p_*}
 \\
\Vect_k \ar[r]^{\Psi_*}_\sim & \Vect_k.
}$$
\end{proof}
\end{Prop}
\begin{Rem}
Note that unlike the $\QCoh$-dualizing sheaf, the $\ICoh$-dualizing sheaf is defined for much bigger class of prestacks (the quasi-coherent dualizing sheaf $\omega_X \in \QCoh(X)$ exists only if $\mathcal O_X \in \QCoh(X)$ is compact). However for the comparison of the morphism of traces induced by pushforward in $\QCoh$-setting and $\ICoh$-setting it is more convenient to work with $\QCoh$-version of dualizing sheaf.
\end{Rem}

\subsection{Computing the trace of pushforward}
In this subsection we discuss morphism of traces in the setting of ind-coherent sheaves. We first note that similar to quasi-coherent sheaves, ind-coherent sheaves are self-dual as an object of $\Cat_k$:

\begin{Theor}[{\cite[Chapter 4, Proposition 6.3.4; Chapter 5, Theorem 4.2.5]{GaitsRozI}}]\label{IndCohSelfDual}\
\begin{enumerate}[1.]
\item For any two $X,Y \in \Sch_{\aft}$ (for the definition see \cite[Chapter 2, 3.5]{GaitsRozI}) the morphism
$$\xymatrix{
\ICoh(X) \otimes \ICoh(Y) \ar[rr] && \ICoh(X \times Y)
}$$
in $\Cat_k$ induced by the functor
$$\xymatrix{
\ICoh(X) \times \ICoh(Y) \ar[rr]^-{\sboxtimes} && \ICoh(X \times Y)
}$$
is an equivalence.

\item For any $X \in \Sch_{\aft}$ the morphisms
$$\xymatrix{
\Vect_k \ar[rr]^-{\Delta_* \omega^{\ICoh}_X} && \ICoh(X \times X) \simeq \ICoh(X) \otimes \ICoh(X)
}$$
and
$$\xymatrix{
\ICoh(X) \otimes \ICoh(X) \simeq \ICoh(X \times X) \ar[rr]^-{\Gamma \circ \Delta^!} && \Vect_k
}$$
$$\xymatrix{
}$$
exhibit $\ICoh(X)$ as a self-dual object in $\Cat_k$.
\end{enumerate}
\end{Theor}

The proof of the following corollary is similar to that of Corollary \ref{TraceOfQCoh} (where we use \cite[Chapter 4, Proposition 5.2.2]{GaitsRozI} to perform base change for ind-coherent sheaves)

\begin{Cor}\label{TraceOfICoh}
Let $X$ be an almost finite type scheme with an endomorphism $g$. Then
$$
\Tr_{2\Cat_k}(g_*) \simeq \Gamma(X,\Delta^! (\Id_X, g)_* \omega^{\ICoh}_X) \simeq \Gamma(X^g, \omega^{\ICoh}_{X^g}).
$$
In particular, in the case when $X$ is proper by Proposition \ref{compare_dualizing} we get an equivalence 
$$
\Tr_{2\Cat_k}(g_*) \simeq \Gamma(X^g, \omega^{\ICoh}_{X^g}) \simeq \Gamma(X^g, \omega_{X^g}).
$$ 
\end{Cor}
\begin{Rem}
Let $\xymatrix{Z \ar[r]^-{p} & \ast}$ be an almost finite type scheme. Note that then
$$\Gamma(Z, \omega_Z) \simeq \Hom_{\QCoh(Z)}(\mathcal O_Z, p^! k) \simeq \Hom_{\Vect_k}(p_*\mathcal O_Z, k) \simeq \Gamma(Z, \mathcal O_Z)^\vee.$$
In particular, using the previous proposition we obtain an equivalence
$$\Tr_{2\Cat_k}(g_*) \simeq \Gamma(X^g, \omega_{X^g}) \simeq \Gamma(X^g, \mathcal O_{X^g})^\vee.$$
\end{Rem}

As a corollary, we get
\begin{Cor}\label{cor:tr_id_icoh_smooth}
Let $X$ be smooth and proper scheme. Then
$$\Tr_{2\Cat_k}(\Id_{\ICoh(X)}) \simeq \Gamma(\mathcal LX, \omega_{\mathcal LX}) \simeq \Gamma(\mathcal LX, \mathcal O_{\mathcal LX})^\vee \overset{\mathrm{HKR}^\vee}{\simeq} \left(\bigoplus_{p=0}^{\dim X} H^p(X, \Omega_X^p) \right)^\vee$$
where $\mathrm{HKR}$ is the Hochschild-Kostant-Rosenberg equivalence (see Corollary \ref{HochKosRos}).
\end{Cor}

\medskip

We now turn to the computation of morphism of traces. Main result of this section is the following
\begin{Prop}\label{TraceInICoh}
Let $\xymatrix{(X,g_X) \ar[r]^-{f} & (Y, g_Y)}$ be an equivariant proper morphism in $\Sch_{\aft}$. Then the induced morphism of traces
$$\xymatrix{
\Gamma(X^{g_X}, \omega^{\ICoh}_{X^{g_X}}) \simeq \Tr_{2\Cat_k}(g_{X*})\ar[rr]^-{\Tr_{2\Cat_k}(f_*)} && \Tr_{2\Cat_k}(g_{Y*}) \simeq \Gamma(Y^{g_Y}, \omega^{\ICoh}_{Y^{g_Y}}) 
}$$
can be obtained by applying the global sections functor $\Gamma(Y^{g_Y},-)$ to the morphism
$$\xymatrix{
(f^g)_* \omega^{\ICoh}_{X^{g_X}} \simeq (f^g)_*(f^g)^! \omega^{\ICoh}_{Y^{g_Y}} \ar[rr] && \omega^{\ICoh}_{Y^{g_Y}}
}$$
in $\ICoh(X)$ induced by the counit of the adjunction $(f^g)_* \dashv (f^g)^!,$ where $\xymatrix{X^{g_X} \ar[r]^-{f^g} & Y^{g_Y}}$ is the induced by $f$ morphism on derived fixed points.

\begin{proof}
The proof is a direct consequence of the fact that the self-duality of ind-coherent sheaves arise from the category of correspondences. Namely, by \cite[Chapter 5, Theorem 4.1.2]{GaitsRozI} the ind-coherent sheaves functor can be lifted to a symmetric monoidal functor
$$\xymatrix{
{\sf Corr}(\Sch_{\aft})^{\sf proper} \ar[rr] && 2\Cat_k,
}$$
where ${\sf Corr}(\Sch_{\aft})$ is a symmetric monoidal $(\infty,2)$-category whose objects are $X \in \Sch_{\sf aft}$, morphisms from $X$ to $Y$  are spans
$$\xymatrix{
X & W \ar[r] \ar[l] & Y
}$$
in $\Sch_{\sf aft}$ (with the composition given by pullbacks), $2$-morphisms are commutative diagrams 
$$\xymatrix{
& W_1 \ar[dr] \ar[dl] \ar[dd]^-{h}
\\
X && Y
\\
& W_2 \ar[ur] \ar[ul]
}$$
where $h$ is proper and the monoidal structure is given by the cartesian product. Informally speaking, the extension of the ind-coherent sheaves to the category of correspondences is given by mapping the span $\xymatrix{X & W \ar[r]^-{t} \ar[l]_-{s} & Y}$ to the morphism $\xymatrix{\ICoh(X) \ar[r]^-{s^!} & \ICoh(W) \ar[r]^-{t_*} & \ICoh(Y)}$  in $2\Cat_k$. We refer to \cite[Chapter 7, Chapter 5]{GaitsRozI} for a throughout discussion of the category of correspondences and to Corollary \ref{morph_of_tr_in_icoh} for a complete proof of the proposition.
\end{proof}
\end{Prop}

\begin{Rem}\label{Pullback_in_cohom}
One can similarly show that the morphism of traces
$$\xymatrix{
\Gamma(Y^{g_Y}, \mathcal O_{Y^{g_Y}}) \simeq \Tr_{2\Cat_k}(g_{Y*})\ar[rr]^-{\Tr_{2\Cat_k}(f^*)} && \Tr_{2\Cat_k}(g_{X*}) \simeq \Gamma(X^{g_X}, \mathcal O_{X^{g_X}}) 
}$$
is induced by the natural map $\xymatrix{\mathcal O_{Y^{g_Y}} \ar[r] & (f^g)_*\mathcal O_{X^{g_X}}}$. It follows that for $X$ and $Y$ smooth and proper with trivial equivariant structure, under the Hochschild-Kostant-Rosenberg equivalence
$$\pi_i\Gamma(\mathcal LX,\mathcal O_{\mathcal LX}) \simeq \bigoplus_{p-q=i} H^q(X,\Omega^p_X) \qquad \pi_i\Gamma(\mathcal LY,\mathcal O_{\mathcal LY}) \simeq \bigoplus_{p-q=i} H^q(Y,\Omega^p_Y)$$
the morphism of traces is exactly the pullback in cohomology.

Note, however, that the strategy above does not give directly the description of the morphism of traces
$$\xymatrix{
\Gamma(\mathcal LX, \mathcal O_{\mcal{L}X})  \simeq \Tr_{2\Cat_k}(\Id_{\QCoh(X)}) \ar[rr]^-{\Tr_{2\Cat_k}(f_*)} && \Tr_{2\Cat_k}(\Id_{\QCoh(Y)}) \simeq \Gamma(\mathcal LY, \mathcal O_{\mcal{L}Y})
}$$
on quasi-coherent sheaves. The reason is that in this case the formalism of traces uses the right adjoint to $f_*$ which is $f^!$. However, the functor $f^!$ does not come from the $\QCoh$ functor out of the category of correspondences (which uses the adjoint pair $f^* \dashv f_*$ instead).
\end{Rem}
Finally, using the identification of Corollary \ref{cor:tr_id_icoh_smooth} we obtain
\begin{Cor}\label{tr_of_pushforward_via_homology}
Let $X,Y$ be smooth and proper with the trivial equivariant structure. Then under the $\mathrm{HKR}$-identification the morphism of traces
$$\xymatrix{\displaystyle\left(\bigoplus_{p=0}^{\dim X} H^p(X, \Omega_X^p)\right)^{\vee} \simeq \Tr_{2\Cat_k}(\Id_{\ICoh(X)}) \ar[rr]^-{\Tr_{2\Cat_k}(f_*)} && \Tr_{2\Cat_k}(\Id_{\ICoh(Y)}) \simeq \displaystyle \left( \bigoplus_{p=0}^{\dim Y} H^p(Y, \Omega_Y^p) \right)^\vee}$$
coincides with the pushforward in homology (which we define as the map dual to the pullback of cohomology).

\begin{proof}
Immediately follows from the fact that for a map of schemes $\xymatrix{Z \ar[r]^-{f} & W}$ the induced map
$$\xymatrix{
\Gamma(Z, \mathcal O_Z)^\vee \simeq \Gamma(Z, \omega_Z) \ar[r] & \Gamma(W, \omega_W) \simeq \Gamma(W, \mathcal O_W)^\vee
}$$
coincides by construction of dualizing sheaves with the map dual to the pullback of functions.
\end{proof}
\end{Cor}

\section{Orientations and traces}\label{sect:CY}
Let $\xymatrix{X \ar[r]^-{f} & Y}$ be a morphism between smooth proper schemes with trivial equivariant structure. Our goal in this section is to give some description of the morphism of traces induced by the pushforward functor $\xymatrix{\QCoh(X) \ar[r]^-{f_*} & \QCoh(Y)}$ (which is automatically equivariant). Since the diagram
$$\xymatrix{
\QCoh(X) \ar[r]^-{f_*} \ar[d]_-{-\otimes \mathcal O_X} & \QCoh(Y) \ar[d]^-{-\otimes \mathcal O_Y}
 \\
\ICoh(X) \ar[r]_-{f_*} & \ICoh(Y)
}$$
commutes and we already have a description of the morphism of traces induced by the functor $\xymatrix{\ICoh(X) \ar[r]^-{f_*} & \ICoh(Y)}$ (where we take identity endomorphisms on both $\ICoh(X)$ and $\ICoh(Y)$), it is enough to understand the morphism of traces induced by $\xymatrix{\QCoh(X) \ar[r]^{-\otimes \mathcal O_X} & \ICoh(X)}$ and analogously for $Y$. We start by introducing the following 

\smallskip

\begin{Def}\label{CalabiYau}
For an almost finite type scheme $Z$ an \emph{orientation} on $Z$ is a choice of an equivalence $\mathcal O_Z \simeq \omega_Z$ in $\QCoh(Z)$.
\end{Def}
\begin{Rem}
Let $\orient\colon \mathcal O_Z \simeq \omega_Z$ be an orientation on $Z$. Then in particular we obtain a self-duality equivalence
$$
\Gamma(Z, \mathcal O_Z) \xymatrix{\ar[r]_\sim^-{\orient} & } \Gamma(Z, \omega_Z) \simeq \Gamma(Z, \mathcal O_Z)^\vee
$$
which is moreover a morphism of $\Gamma(Z, \mathcal O_Z)$-modules. Note that the space of orientations on $Z$ is a torsor over $\Gamma(Z, \mathcal O_Z)^{\times}$, i.e. after a choice of particular orientation, all other orientations are in bijection with invertible functions on $Z$.
\end{Rem}
The relevance of the constructions above to the comparison of traces is explained by 

\begin{Rem}\label{CY_gives_global_sections}
If $X$ is smooth proper scheme, then any orientation $t$ on $\mathcal LX$ induces an equivalence
$$
\pi_0 \Tr_{2\Cat_k}(\Id_{\ICoh(X)}) \simeq \pi_0\Gamma(\mathcal{L}X, \omega_{\mathcal{L}X}) \xymatrix{\ar[r]_\sim^-{\orient^{-1}} & } \pi_0\Gamma( \mathcal LX, \mcal{O}_{\mathcal LX}) \stackrel{\mathrm{HKR}}{\simeq} \bigoplus_{p=0}^{\dim X} H^p(X, \Omega_X^p).
$$
\end{Rem}

\subsection{Serre orientation}
In this and the next subsections we will introduce several orientations on the derived fixed schemes and discuss some of their properties. We start by recalling that there is a well-known isomorphism
$$\bigoplus_{p,q} H^q(X, \Omega_X^p) \simeq \left(\bigoplus_{p,q} H^q(X, \Omega_X^p)\right)^\vee$$
given by the Poincar\'e duality. The next construction shows that this equivalence is in fact induced by an orientation:
\begin{construction}[Serre orientation] \label{SerrePoincareDuality}
Let $X$ be a smooth proper scheme and denote by $\xymatrix{\mathbb V(\mathbb T_X[-1]) \ar[r]^-{j} & X}$ the evident projection map so that in particular we have an equivalence of sheaves
$$j_* \mathcal O_{\mathbb V(\mathbb T_X[-1])} \simeq \Sym_{\QCoh(X)}(\Omega_X[1])$$
(see Corollary \ref{HochKosRos}). By projecting this equivalence to the top exterior summand we obtain a map 
$$\xymatrix{
j_*\mathcal O_{\mathbb V(\mathbb T_X[-1])} \ar[r] & \omega_X
}$$
in $\QCoh(X)$ which using the adjunction $j_* \dashv j^!$ (as the morphism $j$ is proper) induces an equivalence 
$$\xymatrix{
\mathcal O_{\mathbb V(\mathbb T_X[-1])} \ar[r]_-{\sim} & j^! \omega_X  \simeq \omega_{\mathbb V(\mathbb T_X[-1])}
}$$ 
and hence endows $\mathbb V(\mathbb T_X[-1]) \simeq \mcal{L}X$ with an orientation. We will further call this equivalence the \emph{Serre orientation}.
\end{construction}

The following proposition is a formal consequence of the construction
\begin{Prop}\label{Serre_CY_is_local_Poincare}
Let $X$ be smooth proper scheme. Then the equivalence
$$
\bigoplus_{p-q=i} H^q(X, \Omega_X^p) \xymatrix{\ar[r]_-{\sim}^-{\mathrm{HKR}} & \pi_i\Gamma(\mathcal LX, \mathcal O_{\mathcal LX}) \ar[r]_-{\sim}^-{\orient_S} & \pi_i\Gamma(\mathcal LX, \omega_{\mathcal LX}) \ar[r]_-{\sim} & \pi_i\Gamma(\mathcal LX, \mathcal O_{\mathcal LX})^\vee \ar[r]_-{\sim}^-{\mathrm{HKR}^\vee} & }\bigoplus_{q-p=i} H^q(X, \Omega_X^p)^\vee,
$$
where $\orient_S$ is induced by the Serre orientation coincides with the one induced by the classical Serre duality. In particular, it sends a form $\displaystyle\eta \in \bigoplus_{p,q}^{\dim X} H^q(X, \Omega_X^p) $ to the functional $\int_X - \wedge \eta$, i.e. it is given by the usual Poincar\'e pairing.
\end{Prop}

Using the proposition above we obtain 

\begin{Prop}\label{TraceInICohForSmooth}
Let $\xymatrix{X \ar[r]^-{f} & Y}$ be a morphism of smooth proper schemes. Then the induced morphism of traces
\begin{align}\label{form:classical_pushforwad_is_tr_icoh}
\xymatrix{
\displaystyle \bigoplus_{p=0}^{\dim X} H^p(X, \Omega_X^p) \simeq \Tr_{2\Cat_k}(\Id_{\ICoh(X)})\ar[rr]^-{\Tr_{2\Cat_k}(f_*)} && \Tr_{2\Cat_k}(\Id_{\ICoh(Y)}) \simeq \displaystyle \bigoplus_{p=0}^{\dim Y} H^p(Y, \Omega_Y^p)
}
\end{align}
coincides with the pushforward in cohomology (the Poincar\'e dual of the pullback map), where the first and the last equivalences are obtained from \ref{CY_gives_global_sections} using the Serre orientation \ref{SerrePoincareDuality}.

\begin{proof}
By Corollary \ref{tr_of_pushforward_via_homology} the morphism of traces
$$\xymatrix{\displaystyle\left(\bigoplus_{p=0}^{\dim X} H^p(X, \Omega_X^p)\right)^{\vee} \simeq \pi_0\Tr_{2\Cat_k}(\Id_{\ICoh(X)}) \ar[rr]^-{\Tr_{2\Cat_k}(f_*)} && \pi_0\Tr_{2\Cat_k}(\Id_{\ICoh(Y)}) \simeq \displaystyle \left( \bigoplus_{p=0}^{\dim Y} H^p(Y, \Omega_Y^p) \right)^\vee}$$
is dual to the pullback in cohomology. It follows from the Proposition \ref{Serre_CY_is_local_Poincare} that the composition (\ref{form:classical_pushforwad_is_tr_icoh}) first takes Poincar\'e dual of the form, then pushes it forward in homology and then again takes the Poincar\'e dual.
\end{proof}
\end{Prop}

\subsection{Canonical orientation}
Another important CY-structure is given by the
\begin{construction}[Canonical orientation]\label{FixedPointsAreCalabi}
Let $X$ be a smooth proper scheme. Given an endomorphism $\xymatrix{X \ar[r]^-g & X}$ the \emph{derived fixed-points of $g$} defined as the pullback
$$\xymatrix{
X^g \ar[r]^-{i} \ar[d]_-{i} & X \ar[d]^-{(\Id_X,g)}
\\
X \ar[r]_-{\Delta} & X \times X
}$$
admits an orientation which is given by the series of equivalences
$$\mathcal O_{X^g} \simeq i^* \omega_X \otimes i^* \omega_X^{-1} \simeq i^* \omega_X \otimes i^* \omega_{X/X \times X}  \simeq i^* \omega_X \otimes \omega_{X^g/X} \simeq i^! \omega_X \simeq \omega_{X^g}.$$
We will further call this orientation on $X^g$ \emph{canonical}.
\end{construction}

To see why the canonical orientation is relevant, we prove the following
\begin{Prop}\label{TraceIsFromCalabiYau}
For a classical smooth scheme $X$ the morphism of traces
$$\xymatrix{
\Gamma(\mcal{L}X,\mcal{O}_{\mcal{L}X}) \ar[rrr]^-{\Tr_{2\Cat_k}(- \otimes \mathcal O_X)}_-{\sim} &&& \Gamma(\mcal{L}X, \omega_{\mcal{L}X})
}$$
induced by the diagram
$$\xymatrix{
\QCoh(X) \ar[d]_-{- \otimes \mathcal O_X} \ar[rr]^-{\Id_{\QCoh(X)}} && \QCoh(X) \ar[d]^-{- \otimes \mathcal O_X}
\\
\ICoh(X) \ar[rr]_-{\Id_{\ICoh(X)}} && \ICoh(X)
}$$
is obtained by applying the global sections functor $\Gamma(\mcal{L}X,-)$ to the canonical orientation $\xymatrix{\mathcal O_{\mathcal LX} \ar[r]^-{\orient_C}_-{\sim} & \omega_{\mathcal LX}}$ on $\mathcal LX$.
\begin{proof} The proof essentially boils down to coherence of various operations, see Theorem \ref{TraceIsFromCalabiYau_proof} in appendix for a full proof.
\end{proof}
\end{Prop}

\subsection{Group orientations}
The previous subsection raises the question of how one can understand the canonical orientation more explicitly. Our goal now is to show how one can obtain Serre and canonical orientations on $\mcal{L}X$ using various formal group structures on $\mathcal{L}X$. To begin, we first need to fix some concrete way how we trivialize the tangent sheaf to a group:
\begin{construction}\label{trivialization_of_tangent}
Let $\widehat G \in \Grp(\FormModuli_{/X})$ be a formal group over smooth proper scheme $X$ such that the corresponding Lie algebra $\mathfrak g := \Lie_X (\widehat G)$ lies in $\Coh(X)^{<0}$. Consider the pullback diagram
$$\xymatrix{
\widehat G \ar[r]^-{i} \ar[d]^-{i} & X \ar[d]
\\
X \ar@/^/[u]^-{e}  \ar[r] & \widehat{B}_{/X} \widehat G 
}$$
where $\widehat{B}_{/X} \widehat G$ is the completion along $X$ of $B_{/X}\widehat G$. Since the relative tangent sheaf by its universal property is stable under pullbacks we get a trivialization
$$\mathbb T_{\widehat G/X} \simeq i^*\mathbb T_{X/\widehat{B}_{/X} \widehat G} \simeq i^* e^* i^* \mathbb{T}_{X/\widehat{B}_{/X}\widehat G} \simeq i^* e^* \mathbb{T}_{\widehat G/X} \simeq i^*\mathfrak g.$$
\end{construction}
\begin{Rem}
Intuitively, the trivialization $\mathbb T_{\widehat G/X} \simeq i^* \mathfrak g$ above is given by the left-invariant vector fields.
\end{Rem}
Now trivialization of the tangent sheaf to a group sometimes allows to construct orientations. To explain how, recall that for an eventually coconnective almost of finite type morphism of derived schemes $\xymatrix{Z \ar[r]^-{f} & W}$, there is a natural equivalence $f^!- \simeq \omega_f\otimes f^*-$ of functors, and, moreover, if the morphism $f$ is a regular embedding one can explicitly identify $\omega_f$ with the shifted determinant of the normal bundle of $f$ (see e.g. \cite[Corollary 7.2.5.]{Gaits} and \cite[Chapter 9, Section 7]{GaitsRozII} respectively). We will need not only the existence of the equivalences from \textit{loc.\ cit.} but also their constructions, so we review here relevant parts of the theory in $\QCoh$-language. We start with the following general
\begin{Prop}
Let $\mathscr C = \QCoh(X)$ and let $V \in \mathscr C$ be a dualizable object such that $\Sym^{d+1}_{\mathscr C}(V) \simeq 0$ and $\Sym^d_{\mathscr C}(V) \not\simeq 0$ for some $d\ge 0$. Then:

\begin{enumerate}
\item We have $\Sym_{\mathscr C}^{>d}(V) \simeq 0$ and $\Sym_{\mathscr C}(V) \in \mscr{C}$ is dualizable.

\item The top symmetric power $\Sym^d(V)$ is an invertible object of $\mathscr C$ with the inverse equivalent to $\Sym^d_{\mathscr C}(V^\vee)$.

\item The multiplication map followed by the projection on the top summand
$$\xymatrix{
\Sym_{\mathscr C}(V) \otimes \Sym_{\mathscr C}(V) \ar[r] & \Sym_{\mathscr C}(V) \ar[r] & \Sym^d_{\mathscr C}(V)
}$$
is a perfect pairing, i.e. the induced map $\Sym_{\mathscr C}(V)\otimes \Sym_{\mathscr C}^d(V)^{-1} \to \Sym_{\mathscr C}(V)^\vee$ is an equivalence.
\end{enumerate}

\begin{proof}
\begin{enumerate}[1.]
\item Assume that $\Sym_{\mathscr C}^n(V) \simeq 0$ for some $n\in \mathbb Z_{>0}$. Since we are in characteristic zero we have an equivalence $\Sym^n_{\mathscr C}(V) \simeq (V^{\otimes n})^{\Sigma_n}$ and the natural map $\xymatrix{(V^{\otimes n+1})^{\Sigma_{n+1}} \ar[r] & (V^{\otimes n})^{\Sigma_n}}$ admits a section $\mathrm{Nm}^{\Sigma_{n+1}}_{\Sigma_n}/[\Sigma_{n+1} : \Sigma_n]$. Hence $\Sym^{n+1}_{\mathscr C}(V)$ is a direct summand in $(V^{\otimes n+1})^{\Sigma_n} \simeq \Sym_{\mathscr C}^n(V)\otimes V \simeq 0$, so $\Sym^{n+1}_{\mathscr C}(V) \simeq 0$. Finally, using characteristic zero assumption again we deduce that $\Sym^k_{\mathscr C}$ is a colimit over a finite diagram and hence $\Sym_{\mathscr C}(V) \simeq \bigoplus_{k=0}^d \Sym^k_{\mathscr C}(V)$ is dualizable as a finite colimit of dualizable objects.

\item Let us show that the evaluation map $\xymatrix{e:\Sym^d_{\mathscr C}(V)\otimes \Sym^d_{\mathscr C}(V^\vee) \ar[r] & I_{\mathscr C}}$ is an equivalence. Since $\mathscr C=\QCoh(X)$ is a limit of module categories, it is enough to assume $\mathscr C \simeq \Mod_R$ for some connective $k$-algebra $R$. Moreover, by derived Nakayama's lemma we can assume $R$ is discrete. Further, by the usual Nakayama's lemma it is enough to prove the statement for all residue fields of $R$, in which case the statement is clear.

\item Similar to the previous point.

\end{enumerate}
\end{proof}
Let now $A \in \CAlg(\mathscr{C})$ be a commutative algebra object in a presentably symmetric monoidal $k$-linear category $\mathscr C$. Note that the forgetful functor $\xymatrix{\Mod_A(\mathscr C) \ar[r]^-{i} & \mathscr C}$ admits both left $i^*$ and right $i^!$ adjoints explicitly given by
$$i^*(\mathcal F) \simeq A\otimes \mathcal F \qquad\qquad i^!(\mathcal F) \simeq \HHom_{\mathscr{C}}(A ,\mathcal F),$$
where $ \HHom_{\mathscr{C}}(-,-)$ is the inner hom in $\mathscr C$. We then have
\end{Prop}
\begin{Cor}[Grothendieck's formula]\label{simple_Grothendieck_regular_emb}
In the notations above let $A:=\Sym_{\mathscr{C}}(V)$ where $V \in \mathscr C = \QCoh(X)$ is a dualizable object such that moreover $\Sym^{d+1}_{\mathscr C}(V) \simeq 0$ and $\Sym_{\mathscr C}^d(V) \not\simeq 0$ for some $d\ge 0$. Then for any $\mathcal F \in \mathscr{C}$ there is a natural equivalence
$$i^!(\mathcal F) \simeq \omega \otimes_{A} i^*(\mathcal F)$$
in $\Mod_A(\mathscr C)$, where $\omega := A \otimes  \Sym^d(V[1])^{-1}$.

\begin{proof}
By the previous proposition $A$ is dualizable and $A^\vee \simeq \omega$. Hence
$$
i^!(\mathcal F) \simeq  \HHom_{\mathscr{C}}(A, \mathcal F) \simeq A^\vee \otimes \mathcal F \simeq A^\vee\otimes_A (A\otimes \mathcal F) \simeq \omega \otimes_{A} i^*(\mcal{F})
$$
as claimed.
\end{proof}
\end{Cor}
By applying The Grothendieck's formula to the case $\mcal{F}=\mcal{O}_X$, $V = \Omega_X[1]$ we thus obtain an equivalence
$$
T^{ab}\colon \omega_{\mathbb V(\mathbb T_X[-1])/X} \simeq i^*\left(\Sym^{\mathrm{top}}(\Omega_X[1])^{-1}\right) \simeq \Sym^{\mathrm{top}}\left(i^*(\mathbb T_X[-1])\right) \simeq \Sym^{\mathrm{top}}(\mathbb T_{\mathbb V(\mathbb T_X[-1])/X}),
$$
where in the last equivalence we use the abelian group structure on $\mathbb V(\mathbb T_X[-1])$ to identify $i^*(\mathbb T_X[-1])\simeq \mathbb T_{\mathbb V(\mathbb T_X[-1])/X}$. Consequently, given any other trivialization of the relative tangent bundle $\alpha\colon \mathbb T_{\mathbb V(\mathbb T_X[-1])/X} \simeq i^*(\mathbb T_X[-1])$ we can precompose $\Sym^{\mathrm{top}}(\alpha)$ with $T^{ab}$ to obtain an equivalence $\omega_{\mathbb V(\mathbb T_X[-1])/X} \simeq i^*\omega_X^{-1}$. This suggests the following

\smallskip 

\begin{construction}[Loop group and abelian orientations]\label{GrpCY_loop}
Let $X$ be a smooth proper scheme. For any group structure on $\mathbb V(\mathbb T_X[-1])$ by Construction \ref{trivialization_of_tangent} we obtain a trivialization of the relative tangent sheaf $\mathbb T_{\mathbb V(\mathbb T_X[-1])/X}$ and hence by discussion above an equivalence $\omega_{\mathbb V(\mathbb T_X[-1])/X} \simeq i^*\omega_X^{-1}$. By multiplying both sides with $i^*\omega_X$ we obtain an orientation
$$\omega_{\mathbb V(\mathbb T_X[-1])} \simeq \omega_{\mathbb V(\mathbb T_X[-1])/X} \otimes i^*\omega_X \simeq i^*\omega_X^{-1} \otimes i^*\omega_X \simeq \mathcal O_{\mathbb V(\mathbb T_X[-1])}.$$

If group structure on $\mathbb V(\mathbb T_X[-1])$ is pulled back from $\mathcal LX$ via the exponent map $\exp_{\mathcal LX}$, we will call the orientation above \emph{loop group orientation}. In the case when group structure on $\mathbb V(\mathbb T_X[-1])$ is abelian, we will call the corresponding orientation \emph{abelian orientation}.
\end{construction}

\smallskip

The orientations on $\mathbb V(\mathbb T_X[-1])$ from above are in fact not new ones:
\begin{Prop}\label{Calabi_yau_comparison} Let $X$ be smooth proper scheme. Then:
\begin{enumerate}[1.]
\item The abelian group orientation on $\mathcal{L}X$ coincides with the Serre orientation from construction \ref{SerrePoincareDuality}.

\item The loop group orientation on $\mathcal{L}X$ coincides with the canonical structure from construction \ref{FixedPointsAreCalabi}.
\end{enumerate}

\begin{proof}
\begin{enumerate}[1.]
\item Since  $\mathbb V(\mathbb T_X[-1])$ is affine over $X$, the space of orientations $\mathcal O_{\mathbb V(\mathbb T_X[-1])} \simeq \omega_{\mathbb V(\mathbb T_X[-1])}$ on $\mathbb V(\mathbb T_X[-1])$ is equivalent to the space of $i_*\mathcal O_{\mathbb V(\mathbb T_X[-1])}$-linear equivalences $i_*\mathcal O_{\mathbb V(\mathbb T_X[-1])} \simeq i_*\omega_{\mathbb V(\mathbb T_X[-1])}$ in $\QCoh(X)$. Also note that
$$
i_*\mathcal O_{\mathbb V(\mathbb T_X[-1])} \simeq \Sym_{\QCoh(X)}(\Omega_X[1]) \quad\text{and}\quad i_*\omega_{\mathbb V(\mathbb T_X[-1])} \simeq \HHom_{\QCoh(X)}\big{(}\Sym(\Omega_X[1]), \omega_X\big{)}.
$$
In particular, unwinding definitions one finds that the Serre orientation is induced by the non-degenerate pairing
$$\xymatrix{\Sym_{\QCoh(X)}(\Omega_X[1]) \otimes \Sym_{\QCoh(X)}(\Omega_X[1]) \ar[r] & \Sym_{\QCoh(X)}(\Omega_X[1]) \ar[r] & \Sym^{\mathrm{top}}_{\QCoh(X)}(\Omega_X[1])} \simeq \omega_X,$$
where the first map above is the multiplication and the second is the projection on the top summand. But this is precisely the same pairing which we used to construct an equivalence $\omega_{\mathbb V(\mathbb T_X[-1])/X} \simeq \Sym^{\mathrm{top}}\left(i^*(\mathbb T_X[-1])\right)$
(see Proposition \ref{simple_Grothendieck_regular_emb}).

\item Consider more generally arbitrary formal group structure on $\mathbb V(\mathbb T_X[-1])$ induced by a pullback diagram
$$\xymatrix{
\mathbb V(\mathbb T_X[-1]) \ar[r]^-i \ar[d]_-i & X \ar[d]
 \\
X \ar[r] & \widehat B
}$$
for some $\widehat B \in \FormModuli_{X//X}$. This diagram induces an equivalence $C\colon \omega_{\mathbb V(\mathbb T_X[-1])/X} \simeq i^*(\omega_{X/\widehat B})$ of relative dualizing sheaves. Since $\omega_{X/\widehat B} \simeq \omega_X^{-1}$, by multiplying $C$ with $i^*\omega_X$ we obtain an orientation $\orient_{\widehat B}$ on $\mathbb V(\mathbb T_X[-1])$. As an example, if we consider $B = (X\times X)_{\widehat \Delta}$ the orientation obtained this way by definition coincides with the canonical one. Let now $\alpha\colon \mathbb T_{\mathbb V(\mathbb T_X[-1])/X} \simeq i^*(\mathbb T_X[-1])$ be the trivialization of the tangent sheaf obtained from the pullback diagram above. Unwinding the definitions, one finds that the composite equivalence
\begin{align}\label{formula:almost_grp_CY}
\xymatrix{\omega_{\mathbb V(\mathbb T_X[-1])/X} \ar[r]^-{T^{ab}}_-\sim & \Sym^{\mathrm{top}}(\mathbb T_{\mathbb V(\mathbb T_X[-1])/X}) \ar[rr]^-{\Sym^{\mathrm{top}}(\alpha)}_-\sim && \Sym^{\mathrm{top}}(i^* \mathbb T_X[-1]) \ar[r]_-\sim & i^*(\omega_X^{-1})}
\end{align}
coincides with $C$. But by definition (\ref{formula:almost_grp_CY})$\otimes i^*(\omega_X)$ is the group orientation and $C\otimes i^*(\omega_X) = \orient_{\widehat B}$, hence the group orientation coincides with $\orient_{\widehat B}$ as claimed.
\end{enumerate}
\end{proof}
\end{Prop}

\section{The Todd class}\label{sect:todd}
From \Cref{TraceIsFromCalabiYau} and \Cref{Calabi_yau_comparison} we know that the morphism of traces
$$\Gamma(\mathcal LX, \mathcal O_{\mathcal LX}) \simeq \xymatrix{\Tr_{2\Cat_k}(\Id_{\QCoh(X)}) \ar[rr]^-{\Tr_{2\Cat_k}(-\otimes \mathcal O_X)} && \Tr_{2\Cat_k}(\Id_{\ICoh(X)})} \simeq \Gamma(\mathcal LX, \omega_{\mathcal LX})$$
is given by $\Gamma(\mathcal LX, \orient_C)$, where $\orient_C$ is the canonical orientation from \Cref{FixedPointsAreCalabi}. In \Cref{prop:tr_mor_is_gtodd} we will prove that the composite equivalence
$$\xymatrix{
\mathcal O_{\mathcal{L}X} \ar[r]^-{\orient_C}_-{\sim} & \omega_{\mathcal LX} \ar[r]_-{\sim}^-{\orient_S^{-1}} & \mathcal O_{\mathcal LX}
}$$
is given by the determinant (we refer the reader to \cite[Section 3.1]{STV_determinant} for the construction of the determinant map $\xymatrix{\det: \Perf(X) \ar[r] & \Pic(X)}$) of the derivative of the exponential map 
$$
d\exp_{\mathcal LX}\colon i^*\mathbb T_X[-1] \simeq \mathbb T_{\mathbb V(\mathbb T_X[-1])/X} \simeq \exp_{\mathcal LX}^* \mathbb T_{\mathcal LX/X} \simeq i^* \mathbb T_X[-1],
$$
where the first equivalence above is via abelian group structure on $\mathbb V(\mathbb T_X[-1])$, and the second one uses exponent and loop group structure on $\mathcal LX$. In this and next sections we will prove that the determinant of the morphism above coincides with the classical Todd class $\td_X$ which is defined as a multiplicative characteristic class. Our proof is motivated by the following 

\begin{Ex}
Let $G$ be a real Lie group with the corresponding Lie algebra $\mathfrak g$. In a small enough neighborhood of $0$ we then have two trivializations of $\mathbb T_G$ induced by the group structure on $G$ and abelian group structure on $\mathfrak g$ (via the exponential map $\xymatrix{\mathfrak g \ar[r]^-{\exp_G} & G}$). One can then compute (see Lemma \ref{derivative_of_exp_fseries} for the proof in the formal power series setting) that for $x\in\mathfrak g$ close enough to $0$ and $e(x):= \exp_G(x)$  the change of trivialization isomorphism
$$\xymatrix{
\mathfrak g \simeq \mathbb T_{\mathfrak g, 0} \ar[rr]_-\sim^-{+x} && \mathbb T_{\mathfrak g, x} \ar[rr]_-\sim^-{(d\exp_G)_x} && \mathbb T_{G, e(x)} \ar[rr]_-\sim^-{(dL_{e(x)^{-1}})_{e(x)}} && \mathbb T_{\mathfrak g, 0} \simeq \mathfrak g
},$$
(where for $g\in G$ we denote by $\xymatrix{ G \ar[r]^-{L_g} & G}$ the left translations maps by $g$) is given by the linear operator $(1 - e^{-\ad_{\mathfrak g}(x)})/\ad_{\mathfrak g}(x)$. Note that in this way we obtain an $\EEnd(\mathfrak g)$-valued function on $\mathfrak g$
$$x \xymatrix{\ar@{|->}[r] &} \frac{1 - e^{-\ad_{\mathfrak g}(x)}}{\ad_{\mathfrak g}(x)}.$$
\end{Ex}
In this section we will imitate the example above:
\begin{itemize}
\item Given a map $\xymatrix{\mathfrak g \ar[r]^-{\rho} & \EEnd_{\mathscr C}(E)}$ in a $k$-linear presentably symmetric monoidal category $\mathscr C$ and a power series $f\in k[[t]]$ we will construct an $\EEnd_{\mathscr C}(E)$-valued "formal function on $\mathfrak g$" (which is by definition simply a map $\xymatrix{\Sym_{\mscr{C}}(\mathfrak g) \ar[r]^-{f(\rho)} & \EEnd_{\mathscr C}(E)}$ in $\mathscr C$) which informally sends an element $x\in \mathfrak g$ to $f(\rho(x))$. In the special case $\mathscr C=\QCoh(X)$ and $\mathfrak g = \mathbb T_X[-1]$ we will give an interpretation of multiplicative characteristic classes in these terms.

\item Using the interpretation from the previous step the problem of comparing Todd classes reduces to proving that
$$\det\left(\frac{1 - e^{-\ad_{\mathbb T_X[-1]}}}{\ad_{\mathbb T_X[-1]}}\right) = \det(d\exp_{\mathcal LX}),$$
where the left hand side is obtained by applying determinant to the formal function constructed from $f(t) = (1-e^{-t})/t$ and $\rho = \ad_{\mathbb T_X[-1]}$. In order to prove this, we show that both sides make sense for any Lie algebra $\mathfrak g \in \LAlg(\mscr{C})$. Moreover, since both sides are functorial with respect to continuous symmetric monoidal functors the equivalence above can be checked in the classifying category for Lie algebras $\mathcal U_{\Lie}$ (see Construction \ref{univ_cat}). We will show that $\mathcal U_{\Lie}$ admits a set of functors to $\Vect_k$ which detects non-zero morphisms, hence reducing the problem to ordinary $\mathfrak{gl}_n$ in $\Vect_k$ for which the statement is well-known.
\end{itemize}

\subsection{Group-theoretic Todd class}
We return our discussion to the trivialization of tangent bundle to a group: note that if a formal moduli problem over $X$ has two different structures of a formal group, Construction \ref{trivialization_of_tangent} gives two a priori different trivializations of the tangent sheaf. In order to conveniently measure the difference between these trivializations we introduce the following

\begin{construction}
Let $\mathcal Y \in (\FormModuli_{/X})_*$ be a pointed formal moduli problem over smooth proper scheme $X$ such that the pullback of the tangent sheaf $(\mathbb{T}_{\mathcal Y/X})_{|X} \in \QCoh(X)$ is perfect (so that $\mathbb{T}_{\mathcal Y/X} \in \QCoh(\mathcal Y)$ is itself perfect) and let $s_1,s_2$ be two formal group structures on $\mathcal Y$. Then by Construction \ref{trivialization_of_tangent} we get two trivialization of the tangent sheaf $\mathbb T_{\mathcal Y/X}$ and thus by applying the first trivialization and then inverse of the second trivialization we get an automorphism $\gamma_{s_1,s_2}: \mathbb{T}_{\mathcal Y/X} \simeq \mathbb{T}_{\mathcal Y/X}$ in $\QCoh(\mathcal Y)$. We define the \emph{group-theoretic Todd class of $\mathcal Y$} denoted by $\td_{\mathcal Y,s_1,s_2} \in \Gamma(\mathcal Y, \mathcal O_{\mathcal Y})^\times$ as the determinant $\td_{\mathcal Y,s_1,s_2}:=\det(\gamma_{s_1,s_2})$.
\end{construction}
\begin{construction} \label{GTodd}
In the case when $\widehat G$ is a formal group over $X$ it has two distinguished group structures: the abelian one $s_1$ coming from the exponent map and the given one $s_2$. In this case we will further use the notation
$$d\exp_{\widehat G} := \gamma_{s_1,s_2} \qquad\text{and}\qquad \td_{\widehat G} := \det(d\exp_{\widehat G}) \in \Gamma(\widehat G, \mathcal O_{\widehat G})^\times$$
for the corresponding automorphism and its determinant.
\end{construction}

Note that from Proposition \ref{Calabi_yau_comparison} we get
\begin{Prop} \label{prop:tr_mor_is_gtodd}
Let $\xymatrix{\mathcal O_{\mathcal LX} \ar[r]^-{\orient_S}_{\sim} & \omega_{\mathcal LX}}$ be the Serre orientation (construction \ref{SerrePoincareDuality}) and $\xymatrix{\mathcal O_{\mathcal LX} \ar[r]^-{\orient_C}_{\sim} & \omega_{\mathcal LX}}$ be the canonical orientation (construction \ref{FixedPointsAreCalabi}). Then the composite automorphism
$$\xymatrix{
\mathcal O_{\mathcal{L}X} \ar[r]^-{\orient_C}_-{\sim} & \omega_{\mathcal LX} \ar[r]_-{\sim}^-{\orient_S^{-1}} & \mathcal O_{\mathcal LX}
}$$
is given by $\td_{\mathcal LX}$, where we consider $\mathcal LX$ as a formal group over $X$ using the derived loops group structure.

\begin{proof}
By construction and \Cref{Calabi_yau_comparison} under the equivalence
$$\End_{\QCoh(\mathcal LX)}(\mathcal O_{\mathcal LX}) \simeq \End_{\QCoh(\mathcal LX)}(i^* \omega_X^{-1}) \simeq \End_{\QCoh(\mathcal LX)}\left(\Sym^{\mathrm{top}}(i^*\mathbb T_X[-1])\right)$$
the composite $\orient_C^{-1} \circ \orient_S$ is given by the induced map on $\Sym^{\mathrm{top}}$ applied to the composite autoequivalence
$$d\exp_{\mathcal LX} \colon \xymatrix{i^*\mathbb T_X[-1] \ar[r]_-\sim & \mathbb T_{\mathcal LX/X} \ar[r]_-\sim & i^* \mathbb T_X[-1]},$$
which under the equivalence $\Sym^{d}(V[-1]) \simeq \Lambda^d(V)[-d]$ is given by $\Lambda^{\mathrm{top}}(d\exp_{\mathcal LX}[1])$. On the other hand
$$\td_{\mathcal LX} = \det(d\exp_{\mathcal LX}) = \det(d\exp_{\mathcal LX}[1])^{-1} = \Lambda^{\mathrm{top}}(d\exp_{\mathcal LX}[1])^{-1}.$$
It follows $\td_{\mathcal LX} = (\orient_C^{-1} \circ \orient_S)^{-1} = \orient_S^{-1} \circ \orient_C$.
\end{proof}
\end{Prop}

Thus using Propositions \ref{TraceIsFromCalabiYau} and \ref{TraceInICohForSmooth} we conclude

\begin{Cor}\label{Tr_Q_to_I_via_group_Todd}
Let $X$ be a smooth proper scheme. Then under Serre orientation equivalence $\orient_S^{-1}$ the morphism of traces
$$\Gamma(\mathcal LX, \mathcal O_{\mathcal LX}) \simeq \xymatrix{\Tr_{2\Cat_k}(\Id_{\QCoh(X)}) \ar[rr]^-{\Tr_{2\Cat_k}(-\otimes \mathcal O_X)} && \Tr_{2\Cat_k}(\Id_{\ICoh(X)})} \simeq \Gamma(\mathcal LX, \omega_{\mathcal LX}) \overset{\orient_S^{-1}}{\simeq} \Gamma(\mathcal LX, \mathcal O_{\mathcal LX})$$
is given by multiplication with $\td_{\mathcal LX}$.
\end{Cor}

\subsection{Formal functions}
In this subsection we will introduce a construction of formal functions mentioned in the introduction to this section. We first recall the following
\begin{construction}
Let $\mathscr C$ be a presentably symmetric monoidal category and let $V \in \mathscr{C}$ be an object of $\mathscr C$. We define \emph{symmetric algebra with divided powers $\Sym_{\mathrm{dp}, \mathscr C}(V)$ on $V$} as
$$\Sym_{\mathrm{dp}, \mathscr C}(V) := \bigoplus_{n=0}^\infty (V^{\otimes n})^{\Sigma_n}.$$
Note that there exists a canonical map $\xymatrix{\Sym_{\mathrm{dp}, \mathscr C}(V) \ar[r] & \Free_{\mathscr C}^{\mathbb E_1}(V)}$ where $\Free^{\mathbb{E}_1}_{\mathscr{C}}$ denotes the free associative algebra functor in $\mathscr{C}$.
\end{construction}
\begin{Rem}
There is also a canonical norm map $\xymatrix{\Sym_{\mathscr C}(V) \ar[r] & \Sym_{\mathrm{dp}, \mathscr C}(V)}$ which is an equivalence in characteristic zero.
\end{Rem}

\begin{construction}\label{power_series}
Let $\mathscr C$ be a $k$-linear presentably symmetric monoidal category and $\mathfrak g \in \mathscr{C}$ be an object together with a map $\xymatrix{\mathfrak g \ar[r]^-{\rho} & \EEnd_{\mathscr C}(E)}$. Also let

$$f(t)=\sum_{n=0}^\infty f_n t^n \in k[[t]]$$
be a power series with $f_0 = 1$. Define then an "$\EEnd_{\mathscr C}(E)$-valued formal function $f(\rho)$ on $\mathfrak g$" as the composite
$$\xymatrix{
f(\rho): \Sym_{\mathrm{dp}, \mathscr C}(\mathfrak g) \ar[rr]^-{\Sym_{\mathrm{dp}}(\rho)} && \Sym_{\mathrm{dp}, \mathscr C}(\EEnd_{\mathscr C}(E)) \ar@{->}[r] & \Free^{\mathbb{E}_1}_{\mathscr{C}} (\EEnd_{\mathscr C}(E)) \ar[r]^-f & \EEnd_{\mathscr C}(E)
}$$
where the map $\xymatrix{\Free^{\mathbb{E}_1}_{\mathscr{C}}(\EEnd_{\mathscr C}(E)) \ar[r]^-{f} & \EEnd_{\mathscr C}(E)}$ is defined on the $n$-th component as the composition map $\xymatrix{\circ:\EEnd_{\mathscr C}(E)^{\otimes n} \ar[r] & \EEnd_{\mathscr C}(E)}$ followed by the multiplication by the coefficient $f_n$.
\end{construction}

\begin{variant}\label{power_series_variant}
By adjunction
$$
\Hom_{\mathscr C}\big{(}\Sym_{\mathrm{dp}, \mathscr C}(\mathfrak g), \EEnd_{\mathscr C}(E)\big{)} \simeq \Hom_{\mathscr C}\big{(}\Sym_{\mathrm{dp}, \mathscr C}(\mathfrak g)\otimes E, E\big{)}
$$
the morphism $f(\rho)$ corresponds to some map $\xymatrix{\Sym_{\mathrm{dp},\mathscr C}(\mathfrak g)\otimes E \ar[r] & E}$ which we by abuse of notations will also denote by $f(\rho)$. 
\end{variant}

\begin{Ex}
Let $\mathfrak g$ be a Lie algebra in $\Vect_k$ and let $\xymatrix{\mathfrak g \ar[r]^-{\rho} & \EEnd_{\Vect_k}(E)}$ be a representation of $\mathfrak g$. Then by definition for any power series $f(t) \in k[[t]]$ as above and $x\in \mathfrak g$ we have $f(\rho)(x^{\otimes n})= f_n \rho(x)^{\circ n}$.
\end{Ex}

We now give a more geometric description of $f(\rho)$ in the case $\mathscr C=\QCoh(X)$. 

\begin{construction}\label{geometric_formal_function_of_g_vee}
Let $\mathfrak g$ be a Lie algebra such that $\mathfrak g \in \Coh(X)^{<0}$ (where $X$ is smooth proper scheme) and let $\xymatrix{\mathfrak g \ar[r]^-{\rho} & \EEnd_{\QCoh(X)}(E)}$ be some representation where $E \in \QCoh(X)^{\perf}$. Then the composite
$$\xymatrix{
\Sym_{\QCoh(X)}(\mathfrak g) \ar[r]_-{\sim} & \Sym_{\mathrm{dp},\QCoh(X)}(\mathfrak g) \ar[rr]^-{\rho(g)} && \EEnd_{\QCoh(X)}(E)
}$$
induces a map
$$\xymatrix{
\mathcal{O}_X \ar[r] & \Sym_{\QCoh(X)}(\mathfrak g^\vee) \otimes \EEnd_{\QCoh(X)}(E)
}$$
and hence an element in
\begin{gather*}
\Hom_{\QCoh(X)}\big{(}\mathcal O_X, \Sym_{\QCoh(X)}(\mathfrak g^\vee)\otimes \EEnd_{\QCoh(X)}(E)\big{)} \simeq \\
\simeq \Gamma\big{(}X, j_*j^*\EEnd_{\QCoh(X)}(E)\big{)} \simeq \Gamma\big{(}\mathbb V(\mathfrak g), \EEnd_{\QCoh(\mathbb V(\mathfrak g))}(j^*E)\big{)} \simeq \End_{\QCoh(\mathbb V(\mathfrak g))}(j^*E),
\end{gather*}
which we will by abuse of notation denote by the same symbol $f(\rho) \in \Aut_{\QCoh(\mathbb V(\mathfrak g))}(j^* E)$ (where $\xymatrix{\mathbb{V}(\mathfrak g) \ar[r]^-{j} & X}$ is the projection map and  $f(\rho)$ is invertible since $f_0\ne 0$ and $\mathbb V(\mathfrak g)$ is a nil-thickening of $X$). Moreover, since pullbacks preserve perfect objects, the sheaf $j^*E \in \QCoh(\mathbb V(\mathfrak g))$ is also perfect. Consequently, we can take the determinant of the automorphism $f(\rho)$ above to obtain an element
$$
c^f_{\mathfrak g}(E) := \det\big{(}f(\rho)\big{)} \in \Aut_{\QCoh(\mathbb V(\mathfrak g))}(\mathcal O_{\mathbb V(\mathfrak g)}) \simeq \Gamma(\mathbb V(\mathfrak g), \mathcal O_{\mathbb V(\mathfrak g)})^\times .
$$
\end{construction}

\subsection{$\td_{\mathcal LX}$ and multiplicative characteristic classes}
We now show that the construction above is closely related to the theory of multiplicative characteristic classes. Namely, recall the following
\begin{Def}
For a power series $f\in 1 + tk[[t]]$ define a multiplicative characteristic class
$$\xymatrix{
K_0(X) \ar[rr]^-{c^f} && \displaystyle \bigoplus_{p=0}^{\dim X} H^p(X, \Omega_X^p)
}$$
by setting it to be $f(c_1(\mathcal M))$ on line bundles and extending to all vector bundles by multiplicativity and the splitting principle.
\end{Def}

We now prove the following
\begin{Prop}\label{invariant_polynomial_det}
Let $X$ be a smooth algebraic variety, $E \in \QCoh(X)$ be a perfect sheaf considered as a $\mathbb T_X[-1]\simeq \Lie_X(\mathcal LX)$-module via the canonical $\mathcal LX$-equivariant structure from \ref{CanonicalStructure}. Then for any power series $f$ as above the determinant
$$c^f_{\mathbb T_X[-1]}(E) \in \Gamma\left(\mathbb V(\mathbb T_X[-1]), \mathcal O_{\mathbb V(\mathbb T_X[-1])} \right) \simeq \bigoplus_{p=0}^{\dim X} H^p(X, \Omega_X^p)$$
is equal to $c^f(E)$.

\begin{proof}
Let us denote the canonical action of $\mathbb T_X[-1]$ on $E$ by $a$. Since both $c^f_{\mathbb T_X[-1]}(E)=\det\left(f(a)\right)$ and $c^f(E)$ commute with pullbacks and map direct sums to products, by the splitting principle it is enough to prove the statement in the case when $E:=\mathcal M$ is a line bundle. In this case via the equivalence $\widehat{\GL}(\mathcal M) \simeq \widehat{\mathbb G_m}$ the morphism $\det\left(f(a)\right)$ corresponds to $f(a)$ and so it is left to show that the map
$$\xymatrix{
\mathcal O_X \ar[rr]^-{f(a)^\vee} && \Sym_{\QCoh(X)}(\mathbb T_X[-1])^\vee \simeq \Sym_{\QCoh(X)}(\Omega_X[1]).
}$$
dual to
$$\xymatrix{
f(a): \Sym_{\QCoh(X)}(\mathbb T_X[-1]) \ar[rrr]^-{\Sym_{\QCoh(X)}(a)} &&& \Sym_{\QCoh(X)}(\mathcal O_X) \ar[r]_-\sim & \Free^{\bb{E}_1}(\mathcal O_X) \ar[r]^-f & \mathcal O_X
}$$
coincides with $c^f(\mcal{M})$, where we use above that $\EEnd_{\QCoh(X)}(\mathcal M) \simeq \mathcal O_X$ as $\mcal{M}$ is a line bundle. Since by Proposition \ref{AtiyaClassOfLineBundleIsChern} the representation $\xymatrix{\mathbb T_X[-1] \ar[r]^-{a} & \mathcal M}$ classifies $\At(\mathcal M) \simeq c_1(\mathcal M)$ (see Corollary \ref{AtiyaClassOfLineBundleIsChern}) we have $a^\vee = c_1(\mathcal M)$, and so unwinding the construction we find $f(a)^\vee = f(c_1(\mathcal M))$ which is by definition $c^f(\mathcal M)$.
\end{proof}
\end{Prop}

\subsection{Comparison with the classical Todd class}
Recall by Corollary \ref{Tr_Q_to_I_via_group_Todd} that the morphism of traces
\begin{align}\label{formula:morphism_of_tr_QICoh}
\Gamma(\mathcal LX, \mathcal O_{\mathcal LX}) \simeq \xymatrix{\Tr_{2\Cat_k}(\Id_{\QCoh(X)}) \ar[rr]^-{\Tr_{2\Cat_k}(-\otimes \mathcal O_X)} && \Tr_{2\Cat_k}(\Id_{\ICoh(X)})} \simeq \Gamma(\mathcal LX, \omega_{\mathcal LX}) \overset{\orient_S^{-1}}{\simeq} \Gamma(\mathcal LX, \mathcal O_{\mathcal LX})
\end{align}
is closely related to $d\exp_{\mathcal LX}$. The following theorem provides a description of $d\exp_{\widehat G}$ for an arbitrary formal group $\widehat G$ over $X$
\begin{Theor}\label{lite_BCH}
Let $\widehat G$ be a formal group over $X$ such that $\mathfrak g := \Lie_X(\widehat G) \in \Coh^{<0}$. Then
$$
d\exp_{\widehat G} = \frac{1-e^{-\ad_{\mathfrak g}}}{\ad_{\mathfrak g}}.
$$
\end{Theor}
\begin{Rem}
In the theorem above one can drop any assumptions on $\mathfrak g$ and smoothness assumption on $X$ if instead of $\QCoh$-version one considers $\mathfrak g$ as a Lie algebra in $\ICoh(X)$.
\end{Rem}

The proof of (a generalization of) this theorem is the content of the next section. Here we will only use this theorem to deduce that the group theoretic Todd class $\td_{\mathcal LX}$ (see Construction \ref{GTodd}) coincides with the classical one and thus will give a concrete description of the morphism of traces (\ref{formula:morphism_of_tr_QICoh}). Namely, recall that the classical Todd class $\td_X$ of $X$ is defined as
$$\td_X := c^{1/f}(\mathbb T_X) = c^f(\mathbb T_X[-1]) \qquad\text{where}\qquad f(t) = \frac{1 - e^{-t}}{t}.$$

\begin{Cor}\label{Loop_Todd_vs_classical}
Let $X$ be a smooth proper scheme. Then $\td_{\mathcal LX} = \td_X$.

\begin{proof}
Note that by Proposition \ref{invariant_polynomial_det} above we have
$$\td_X = \det\left(\frac{1 - e^{-\ad_{\mathbb T_X[-1]}}}{\ad_{\mathbb T_X[-1]}}\right),$$
where $\ad_{\mathbb T_X[-1]}$ is the adjoint representation of $\mathbb T_X[-1]$. Consequently, since the group-theoretic Todd class $\td_{\mathcal LX}$ was defined as the determinant of $d\exp_{\mathcal LX}$ it is enough to prove that
$$
d\exp_{\mathcal LX} = \frac{1 - e^{-\ad_{\mathbb T_X[-1]}}}{\ad_{\mathbb T_X[-1]}},
$$
which is a special case of the Theorem \ref{lite_BCH} above.
\end{proof}
\end{Cor}

Finally, we can describe the morphism of traces (\ref{formula:morphism_of_tr_QICoh}) in the classical terms
\begin{Cor}\label{morphism_of_Tr_is_classical_Todd}
Let $X$ be a smooth proper scheme. Then under the Serre orientation and $\mathrm{HKR}$ identifications the morphism of traces
$$\bigoplus_{p,q} H^{q,p}(X) \overset{\mathrm{HKR}}{\simeq} \xymatrix{\pi_* \Tr_{2\Cat_k}(\Id_{\QCoh(X)}) \ar[rrr]^-{\Tr_{2\Cat_k}(-\otimes \mathcal O_X)} &&& \pi_* \Tr_{2\Cat_k}(\Id_{\ICoh(X)})} \overset{\mathrm{HKR}\circ \orient_S^{-1}}{\simeq} \bigoplus_{p,q} H^{q,p}(X)$$
is given by multiplication with $\td_X$.

\begin{proof}
By Corollary \ref{Tr_Q_to_I_via_group_Todd} we know that the morphism of traces is given by multiplication with the group theoretic Todd class $\td_{\mathcal LX}$, and by the previous Corollary $\td_{\mathcal LX} = \td_X$. 
\end{proof}
\end{Cor}

\section{Abstract exponential}
The goal of this section is to prove Theorem \ref{lite_BCH} by first extending it to arbitrary Lie algebras in any stable presentably symmetric monoidal $k$-linear category and then by reducing this more general statement to the case of $\mathfrak{gl}_V$ for $V\in \Vect_k^\heartsuit$. More concretely, let $\widehat G$ be a formal group over $X$ with the corresponding Lie algebra $\mathfrak g$. Note that the morphism $\xymatrix{i^* \mathfrak g \ar[r]^-{d\exp_{\widehat G}} & i^* \mathfrak g}$ via a series of adjunctions similar to that in Construction \ref{geometric_formal_function_of_g_vee} corresponds to some morphism
$$\xymatrix{
\widetilde{d\exp_{\widehat G}} \colon \Sym_{\QCoh(X)}(\mathfrak g) \otimes \mathfrak g \ar[r] & \mathfrak g
}$$
and the same is true for $(1-e^{-\ad_{\mathfrak g}})/\ad_{\mathfrak g}$ by Variant \ref{power_series_variant}. We will show in this section that these two morphisms are actually equal. Before proceeding to the proof, we will first discuss theory of tangent comodules of arbitrary cocommutative coalgebra (which plays the role of tangent space to a formal moduli problem) and review some generalities about operads.

\subsection{Tangent comodule}
In this subsection we describe the construction of relative tangent sheaf in an abstract setting. Informally, the construction will be dual to the cotangent complex formalism developed in \cite[Section 7.3]{HA}. However, since in \cite[Section 7.3]{HA} for most of the statements the category $\mscr{C}$ is assumed to be presentable and we are rather interested in the case of $\mscr{C}^{\op}$ (so that $\CAlg(\mscr{C}^{\op})=\coCAlg(\mscr{C})^{\op}$ is the opposite to the category of cocommutative coalgebras in $\mscr{C}$) we explain the construction here in full details. We start by introducing the following formal
\begin{Def}\label{costabilization}
Let $\mscr{C}$ be a finitely cocomplete category (i.e. $\mscr C$ admits all finite colimits). Define then a {\bfseries costabilization of $\mscr{C}$} denoted by $\coStab(\mscr{C})$ as
$$\xymatrix{
\coStab(\mscr{C}):=\Stab(\mscr{C}^{\op})^{\op}.
}$$
\end{Def}

\begin{Rem}\label{adjoint_to_tang} Notice that the category $\coStab(\mscr{C})$ is always stable. Moreover, since by definition the costabilization $\coStab(\mscr{C})$ of $\mscr{C}$ can be concretely described as the limit of the diagram
$$\xymatrix{
... \ar[r]^-{\Sigma} & \mscr{C}_{/\emptyset} \ar[r]^-{\Sigma} & \mscr{C}_{/\emptyset}\ar[r]^-{\Sigma} & \mscr{C}_{/\emptyset}
}$$
in $\Cat_{\infty}$, where $\emptyset \in \mscr{C}$ is the initial object, we see that if the category $\mscr{C}$ is presentable then so is $\coStab(\mscr{C})$. In particular, in this case the evident projection functor
$$\xymatrix{
\Sigma^{\infty}_{\sf co}: \coStab(\mscr{C})=\Stab(\mscr{C}^{\op})^{\op} \ar[rr]^-{(\Omega^{\infty}_{\mscr{C}^{\op}})^{\op}} && (\mscr{C}^{\op})^{\op} \simeq \mscr{C}
}$$
by the adjoint functor theorem admits a right adjoint which we will further denote by $\Omega^{\infty}_{\sf co}$.
\end{Rem}

\begin{Rem}\label{pushforward}
Let $\mscr{C}$ be a finitely cocomplete category and $\xymatrix{A \ar[r]^-{f} & B}$ be a morphism in $\mscr{C}$. Then the induced functor $\xymatrix{\mscr{C}_{A/} \ar[rr]^-{- \sqcup_A B} && \mscr{C}_{B/}}$ preserves colimits and consequently induces a colimit-preserving functor $\xymatrix{\coStab(\mscr{C}_{A/}) \ar[rr]^-{f_*} && \coStab(\mscr{C}_{B/})}$. Moreover, in the case $\mscr{C}$ is presentable by the adjoint functor theorem we see that the functor $f_*$ admits a right adjoint $f_* \dashv f^*$.
\end{Rem}

\smallskip

\begin{Ex}\label{costab_is_comodules} Let $\mscr{C} \in \CAlg(\Pr^{\Ll, \st}_{\infty})$ be a stable presentably symmetric monoidal category and let $C \in \coCAlg(\mscr{C})$ be a cocommutative coalgebra object in $\mscr{C}$. We then get an equivalence
$$
\coStab(\coCAlg(\mscr{C})_{C/})=\Stab\big{(} (\coCAlg(\mscr{C})_{C/})^{\op} \big{)}^{\op}= \Stab\big{(} \CAlg(\mscr{C}^{\op})_{/C} \big{)}^{\op} \simeq \Mod_C(\mscr{C}^{\op})^{\op}= \coMod_C(\mscr{C})
$$
where the middle equivalence follows from \cite[Theorem 7.3.4.13]{HA}.
\end{Ex}

\smallskip

Example \ref{costab_is_comodules} and Remark \ref{adjoint_to_tang} above motivate the following

\smallskip

\begin{Def}\label{tangent_complex}
Let $\mscr{C}$ be a presentable category and $A \in \mscr{C}$ be an object. Define then a {\bfseries tangent complex to $A$} denoted by $\bb{T}_A \in  \coStab(\mscr{C}_{A/})$ as the image of $\xymatrix{(A \ar[r]^-{\Id_A} & A)} \in \mscr{C}_{A/}$ under the functor
$$\xymatrix{
\mscr{C}_{A/} \ar[rr]^-{\Omega^{\infty}_{\sf co}} && \coStab(\mscr{C}_{A/}).
}$$
\end{Def}

\begin{Rem}\label{all_is_tangent}
Let $\mscr{C}$ be a presentable category and $\xymatrix{A \ar[r]^-{f} & B}$ be a morphism in $\mscr{C}$. Since the diagram
$$\xymatrix{
 \coStab(\mscr{C}_{A/}) \ar[d]_-{\Sigma^{\infty}_{\sf co}} \ar[rr]^-{f_*} &&  \coStab(\mscr{C}_{B/}) \ar[d]^-{\Sigma^{\infty}_{\sf co}}
 \\
 \mscr{C}_{A/} \ar[rr]_-{- \sqcup_A B} && \mscr{C}_{B/}
}$$
by construction commutes, we see that the diagram of right adjoints 
$$\xymatrix{
 \coStab(\mscr{C}_{A/}) && \ar[ll]_-{f^*}\coStab(\mscr{C}_{B/}) 
 \\
 \mscr{C}_{A/} \ar[u]^-{\Omega^{\infty}_{\sf co}} && \mscr{C}_{B/} \ar[ll] \ar[u]_-{\Omega^{\infty}_{\sf co}}
}$$
 also commutes. In particular, we get an equivalence $\xymatrix{\Omega^{\infty}_{\sf co}(A \ar[r]^-{f} & B) \simeq f^* \bb{T}_B}$ in $\coStab(\mscr{C}_{A/})$. 
\end{Rem}

Using the remark \ref{all_is_tangent} we can introduce the following
\begin{Def}\label{relative_cotang}
Let $\mscr{C}$ be a presentable category and $\xymatrix{A \ar[r]^-{f} & B}$ be a morphism in $\mscr{C}$. Define then the {\bfseries relative tangent complex} denoted by $\bb{T}_{A/B} \in \coStab(\mscr{C}_{A/})$ as the fiber 
$$\xymatrix{
\bb{T}_{A/B}:={\sf fib}(\bb{T}_A \ar[r] & f^*\bb{T}_B),
}$$
where the morphism $\xymatrix{\bb{T}_A \ar[r] & f^*\bb{T}_B}$ is obtained by applying the functor $\xymatrix{\mscr{C}_{A/} \ar[r]^-{\Omega^{\infty}_{\sf co}} & \coStab(\mscr{C}_{A/})}$ to the morphism $\xymatrix{A \ar[r]^-{f} & B}$ in $\mscr{C}_{A/}$.
\end{Def}

We end this subsection with the following

\begin{Prop}\label{tangent_to_pullback}
Let $\mscr{C}$ be a presentable category and
$$\xymatrix{
A \ar[r]^-{f} \ar[d]_-{g} & B \ar[d]^-{t}
\\
C \ar[r]_-{h} & D 
}$$
be a pullback square in $\mscr{C}$. Then there is a canonical equivalence
$$\mathbb T_{A/C} \simeq f^* \mathbb T_{B/D}$$
in $\coStab(\mscr{C}_{A/})$.

\begin{proof}
Since the diagram above can be also considered as a pullback square in $\mscr{C}_{A/}$ and the functor $\xymatrix{\mscr{C}_{A/} \ar[r]^-{\Omega^{\infty}_{\sf co}} & \coStab(\mscr{C}_{A/})}$ preserves limits (being right adjoint) using the remark \ref{all_is_tangent} we get a pullback
$$\xymatrix{
\bb{T}_A \ar[r] \ar[d] &  f^* \bb{T}_B \ar[d]
\\
g^* \bb{T}_C \ar[r] & f^*t^* \bb{T}_D
}$$
in $ \coStab(\mscr{C}_{A/})$. In particular, passing to the fibers of the vertical morphisms we get an equivalence
$$\xymatrix{
\bb{T}_{A/C}={\sf fib}(\bb{T}_A \ar[r] & g^* \bb{T}_C) \simeq {\sf fib}(f^* \bb{T}_B \ar[r] & f^* t^* \bb{T}_D) \simeq f^* {\sf fib}(\bb{T}_B \ar[r] & t^* \bb{T}_D)=f^* \bb{T}_{B/D}
}$$
as claimed.
\end{proof}
\end{Prop}

\subsection{Monoidal category built from an operad}
In this subsection we describe the construction of the universal category one can assign to an operad we need in the proof of Proposition \ref{proof_in_vect}. 

\begin{Prop}\label{univ_cat}
Let $\mcal{O}$ be an $\infty$-operad in $\Vect_k$ (we refer the reader to \cite[Chapter 6, Section 1]{GaitsRozII} and \cite[Chapter 2]{HA} for a discussion of $\infty$-operads). Then there exists a symmetric monoidal $k$-linear category $\mcal{U}_{\mcal{O}} \in \CAlg(\Cat_k)$ such that for any symmetric monoidal $k$-linear category $\mscr{C} \in \CAlg(\Cat_k)$ there is a natural equivalence
$$
\Funct_{\CAlg(\Cat_k)}(\mcal{U}_{\mcal{O}},\mscr{C}) \simeq \Alg_{\mcal{O}}(\mscr{C}),
$$
where $\Funct_{\CAlg(\Cat_k)}(\mcal{U}_{\mcal{O}},\mscr{C})$ is the full subcategory of $\Funct^{\otimes}(\mcal{U}_{\mcal{O}},\mscr{C})$ spanned by $k$-linear symmetric monoidal functors.
\begin{proof}
The category $\mcal{U}_{\mcal{O}}$ is simply the category of $\Vect_k$-presheaves $\Psh_k({\sf Env}_k(\mcal{O}))$ on the $\Vect_k$-enriched analogue ${\sf Env}_k(\mcal{O})$ of the monoidal envelope of $\mcal{O}$ from \cite[Proposition 2.2.4.1]{HA} with the Day convolution monoidal structure (see \cite[Example 2.2.6.17]{HA} and \cite{Glas}). The desired universal property follows from $\Vect_k$-enriched version of \cite[Proposition 2.2.4.9]{HA} and \cite[Lemma 2.13]{Glas}.
\end{proof}
\end{Prop}

\begin{Rem}
Suppose that the operad $\mcal{O}$ has only one color. Then the category $\mcal{U}_{\mcal{O}}$ can be informally constructed as follows: first, one considers the $\ii$-category ${\sf Env}_k(\mcal{O})$ whose objects are natural numbers and morphisms are described as
$$\Hom_{{\sf Env}_k(\mcal{O})}(n,m) \simeq \bigoplus_{f: n \rightarrow m} \bigotimes_{i=1}^m \mcal{O}(f^{-1}(i)).$$
Moreover, the category ${\sf Env}_k(\mcal{O})$ has a natural symmetric monoidal structure given by addition. The category $\mcal{U}_{\mcal{O}}$ is then the category of $\Vect_k$-presheaves on ${\sf Env}_k(\mcal{O})$ with the symmetric monoidal structure determined by the fact that it preserves colimits and that the Yoneda's embedding $\xymatrix{{\sf Env}_k(\mcal{O}) \ar[r] & \mcal{U}_{\mcal{O}}}$ is symmetric monoidal. The equivalence
$$\xymatrix{
\Funct_{\CAlg(\Cat_k)}(\mcal{U}_{\mcal{O}},\mscr{C}) \ar[r]_-{\sim} & \Alg_{\mcal{O}}(\mscr{C})
}$$
for $\mscr{C} \in \CAlg(\Cat_k)$ is then simply given by the evaluation at $\Hom_{{\sf Env}_k(\mcal{O})}(-,1) \in \mcal{U}_{\mcal{O}}$.
\end{Rem}

\subsection{Proof of Theorem \ref{lite_BCH}}\label{sect:end_of_our_to_classical_comparison}
In this section we prove a generalization of the Theorem \ref{lite_BCH} to an arbitrary $k$-linear presentably symmetric monoidal category $\mathscr C$ with the monoidal unit $I \in \mathscr C$. Our first step is to define $d\exp_{\widehat G}$ in this context. Note that given a formal group $\widehat G \in \Grp(\FormModuli_{/X})$ one can trivialize its tangent sheaf by applying the base change to the pullback diagram

\begin{align}\label{skdhbchjsdbchjsbhc}
\xymatrix{
\widehat G \ar[rrr]^i \ar[d]_\Delta &&& X \ar[d]^e
\\
\widehat G \times \widehat G \ar[rr]_-{\Id_{\widehat G}\times \inv_G} && \widehat G\times \widehat G \ar[r]_-m & \widehat G
}
\end{align}
where the lower horizontal map informally sends a pair $(g,h)$ to $g\cdot h^{-1}$. The morphism $d\exp_{\widehat G}$ is then defined by comparing two trivializations of the tangent sheaf that come from two formal group structures on $\widehat G$ (the initial one and the abelian one).

 Using the formalism of tangent comodules one can emulate the same construction in algebraic setting. Note first that by Yoneda's lemma the diagram analogous to (\ref{skdhbchjsdbchjsbhc}) is fibered for any group object in any category admitting finite limits. Consequently, we can introduce the following

\begin{construction}
Let $\mathscr C$ be a $k$-linear symmetric presentably monoidal category and $\mathfrak g \in \LAlg(\mathscr C)$ be a Lie algebra in $\mathscr C$ (so that we get a group object $U(\mathfrak g)$ in the category $\coCAlg(\ICoh(X))$ of cocommutative coalgebras in $\mathscr{C}$). Consider the pullback square of cocommutative coalgebras
$$\xymatrix{
U(\mathfrak g) \ar[rr]\ar[d] & & I \ar[d]
\\
U(\mathfrak g)\otimes U(\mathfrak g) \ar[rr] & & U(\mathfrak g).
}$$
Using Proposition \ref{tangent_to_pullback} this diagram induces an equivalence of $U(\mathfrak g)$-comodules $\mathbb T_{U(\mathfrak g)}\simeq U(\mathfrak g)\otimes \mathfrak g$. As in Construction \ref{GTodd} by comparing the trivial Lie algebra structure on $\mathfrak g$ with the given one we obtain an autoequivalence $\xymatrix{d\exp_{\mathfrak g}:U(\mathfrak g)\otimes \mathfrak g \ar[r]_{\sim} & U(\mathfrak g)\otimes \mathfrak g}$ of $U(\mathfrak g)$-comodules and hence by adjunction a morphism $\xymatrix{\widetilde{d\exp}_{\mathfrak g}: U(\mathfrak g)\otimes \mathfrak g \ar[r] & \mathfrak g}$ in $\mathscr C$.
\end{construction}

\smallskip 

If we take $\mathscr C = \QCoh(X)$ and $\mathfrak g = \Lie_{X} \widehat G$ in the construction above the morphism $\widetilde{d\exp}_{\mathfrak g}$ is precisely $\widetilde{d\exp}_{\widehat G}$ from the introduction to this section. Consequently, to conclude it is left to prove 

\smallskip

\begin{Prop}\label{proof_in_vect}
Let $\mathscr C$ be a stable symmetric presentably monoidal $k$-linear category, $\mathfrak g \in \LAlg(\mscr{C})$ be a Lie algebra in $\mathscr C$. Then
$$
\widetilde{d\exp}_{\mathfrak g} = \frac{1-e^{-\ad_{\mathfrak g}}}{\ad_{\mathfrak g}},
$$
where the morphism on the right is obtained from Variant \ref{power_series_variant} by applying $f=(1-e^{-t})/t$ to the adjoint representation of $\mathfrak g$.

\begin{proof}
We first argue that it is sufficient to prove the equality holds in $\Vect_k$ for discrete free Lie algebras. Indeed, since both of the morphisms $\widetilde{d\exp_{\mathfrak g}}$ and $(1-e^{-\ad_{\mathfrak g}})/{\ad_{\mathfrak g}}$ are functorial with respect to continuous monoidal functors, to show that $\widetilde{d\exp}_{\mathfrak g} - (1-e^{-\ad_{\mathfrak g}})/{\ad_{\mathfrak g}}=0$ it is sufficient to prove the statement for the Lie algebra $\mathfrak{g}:=\Hom_{{\sf Env}_k(\Lie)}(-,I) \in \Lie(\mcal{U}_{\Lie})$, where ${\sf Env}_k(\Lie)$ and $\mcal{U}_{\Lie}$ are the universal categories from Proposition \ref{univ_cat}. Now since
$$
\Hom_{\mcal{U}_{\Lie}}(U(\mathfrak g) \otimes \mathfrak g, \mathfrak g) \simeq \Hom_{\mcal{U}_{\Lie}}( \bigoplus_{n \geq 0} \mathfrak g^{\otimes n+1}, \mathfrak g) \simeq \prod_{n \geq 0} \Hom_{\mcal{U}_{\Lie}}(\mathfrak g^{\otimes n+1}, \mathfrak g) \simeq  \prod_{n \geq 0} \Lie(n+1)
$$
we see that if $\widetilde{d\exp}_{\mathfrak g} - (1-e^{-\ad_{\mathfrak g}})/{\ad_{\mathfrak g}} \neq 0$, then there exists $n \in \bb{N}$ such that ${\sf pr}_{\Lie(n+1)}(\widetilde{d\exp}_{\mathfrak g} - (1-e^{-\ad_{\mathfrak g}})/{\ad_{\mathfrak g}}) \neq 0$. In particular, since the evident map 
$$\xymatrix{
\Lie(n+1) \ar[r] & \Hom_{\Vect_k}\big{(} {\sf Free}(n+1)^{\otimes n+1},{\sf Free}(n+1) \big{)}
}$$
is injective, where ${\sf Free}(n+1) \in \LAlg(\Vect_k^\heartsuit) \subset \LAlg(\Vect_k)$ is the discrete free Lie algebra in $\Vect_k$ on $(n+1)$-dimensional vector space, we see that $\widetilde{d\exp}_{\mathfrak g} - (1-e^{-\ad_{\mathfrak g}})/{\ad_{\mathfrak g}}$ should be nonzero for ${\sf Free}(n+1)$. Consequently, to conclude the statement of the proposition it is sufficient to prove that $\widetilde{d\exp}_{\mathfrak g} = (1-e^{-\ad_{\mathfrak g}})/{\ad_{\mathfrak g}}$ for discrete free Lie algebras in $\Vect_k$.

Now since for a free, discrete Lie algebra the adjoint representation is faithful, we can further assume that $\mathfrak g = \gl_V$ for $V\in \Vect_k^{\heartsuit}$. In geometric terms, we want to compute the derivative $\xymatrix{d\exp_{\widehat{\GL}_V}\colon \exp_{\widehat{\GL}_V*} \mathbb T_{\mathbb V(\gl_V)} \ar[r] & \mathbb T_{\widehat{\GL}_V}}$ of the exponential map $\xymatrix{\exp_{\widehat{\GL}_V}\colon \mathbb V(\gl_V) \ar[r] & \widehat{\GL}_V}$. Since by construction tangent complexes are determined by their restriction to $A$-points, it is enough to compute the induced map for any $\xymatrix{\Spec A \ar[r]^-{x} & \mathbb V(\gl_V)}$. Moreover, since formal moduli are determined by their restriction to Artinian local $k$-algebras, we can assume that this is the case. Since by construction the $\ICoh$-tangent sheaf is defined as the Serre dual of $\Pro(\QCoh)$-cotangent sheaf $\widehat{\mathbb L}_{-}$ (see \cite[Chapter 1, Section 4.4]{GaitsRozII}), it is enough to compute the induced map on pro-cotangent sheaves (we need to work with the pro-categories since $\gl_V$ is an ind-scheme for infinite dimensional $V$). Let $M$ be a connective $A$-module and let $A\oplus M$ be the corresponding trivial square zero extension. By definition the space $\Hom_{\Pro(\Mod_A)}(x^* \widehat{\mathbb L}_{\mathbb V(\mathfrak{gl}_V)}, M)$ classifies lifts in the diagram
$$\xymatrix{
\Spec A \ar[r]^-x\ar@{_(->}[d] & \mathbb V(\mathfrak{gl}_V) 
\\
\Spec A\oplus M \ar@{-->}[ru]
}$$
and the differential map $\xymatrix{\Hom_{\Pro(\Mod_A)}(x^*\widehat{\mathbb L}_{\mathbb V(\gl_V)}, M) \ar[r] & \Hom_{\Pro(\Mod_A)}(x^*\exp_{\widehat{\GL}_V}^* \widehat{\mathbb L}_{\widehat{\GL}_V}, M)}$ correspond to the postcomposition with $\exp_{\widehat{\GL}_V}$. Unwinding the definitions, one finds that $x^*\widehat{\mathbb L}_{\mathbb V(\gl_V)} \simeq "\prolim" A\otimes_k V_i^\vee \in \Pro(\Mmod_A)$ where $\{V_i\}$ is the diagram of finite dimensional $k$-vector subspaces of $V$ and analogously for $x^*\widehat{\mathbb L}_{\widehat{\GL}_V}$. Hence both $x^*\widehat{\mathbb L}_{\mathbb V(\gl_V)}$ and $x^*\exp_{\widehat{\GL}_V}^* \widehat{\mathbb L}_{\widehat{\GL}_V}$ are pro-free of finite rank $A$-modules (as $V\in \Vect_k^{\heartsuit}$), and so by Yoneda's lemma it is enough to understand the morphism $\Hom_{\Pro(\Mod_A)}(\exp_{\widehat{\GL}_V,x}^*, M)$ for all free of finite rank modules $M$. Further, since each such module is a direct sum of $A$, we can assume $M=A$ (in this case $A\oplus M \simeq A[\epsilon] := A\otimes_k k[\epsilon]$, where $k[\epsilon] := k[\epsilon]/(\epsilon^2)$ with $\deg(\epsilon) = 0$).

Now unwinding the definitions one finds that for a local Artinian augmented $k$-algebra $A$ we have
$$\mathbb V(\gl_V)(A) \simeq \gl_V(\mathfrak m_A) := \End_{\Mod_A}(V\otimes_k A) \underset{\End_k(V)}{\times}\{0\} \qquad \widehat{\GL}_V(A) \simeq \widehat\GL_V(\mathfrak m_A) :=\Aut_{\Mod_A}(V\otimes_k A) \underset{\Aut_k(V)}{\times} \{\Id_V\}$$
and the exponential map $\exp_{\widehat{\GL}_V}(A)$ sends a matrix $X$ to $\sum_{n=0}^\infty \frac{X^n}{n!}$. Hence the general case follows from the following well-known lemma:

\begin{Lemma}\label{derivative_of_exp_fseries}
Let $V$ be a discrete $k$-vector space and let $A$ be a local Artinian augmented $k$-algebra. Then for each $X,Y \in \gl_V(\mathfrak m_A)$ we have an equality
$$e^{-X}e^{X+\epsilon Y} = 1 + \epsilon \cdot \frac{1-e^{-\ad_X}}{\ad_X}(Y)$$
in $\widehat{\GL}_V(A[\epsilon])$.

\begin{proof}
It is sufficient to prove the statement in the universal case when $A=\mathbb Q[\epsilon]\langle\langle X,Y \rangle\rangle$ the free ring of non-commutative power series on two variables over $\mathbb Q[\epsilon]$. We have
\begin{align}\label{formula:der_of_exp}
e^{-X} e^{X+\epsilon Y} = e^{-X} \sum_{n=0}^\infty \frac{(X+\epsilon Y)^n}{n!} = e^{-X}\cdot\sum_{n=0}^\infty\frac{1}{n!}\left(X^n + \epsilon\sum_{k=0}^{n-1} X^k Y X^{n-1-k}\right) = 1 + e^{-X} \cdot \sum_{n=0}^\infty \frac{\epsilon}{n!}\sum_{k=0}^{n-1} X^k Y X^{n-1-k}.
\end{align}
Note that
$$\ad_X\left(\sum_{k=0}^{n-1} X^k Y X^{n-1-k}\right) = \sum_{k=0}^{n-1}\left( X^{k+1} Y X^{n-1-k} - X^k Y X^{n-k} \right)= X^n Y - YX^n.$$
Hence by applying $\ad_X$ to both sides of (\ref{formula:der_of_exp}) we obtain
$$
\ad_X\left(e^{-X} e^{X+\epsilon Y}\right) = \ad_X \left( e^{-X} \cdot \sum_{n=0}^\infty \frac{\epsilon}{n!}\sum_{k=0}^{n-1} X^k Y X^{n-1-k}\right)=e^{-X} \sum_{n=0}^\infty \frac{\epsilon}{n!}(X^n Y-Y X^n)=$$
$$
= \epsilon \cdot e^{-X} \cdot \left(e^X Y - Y e^X\right) = \epsilon \cdot \left(Y - \Ad_{e^{-X}}(Y)\right) = \epsilon \cdot\left(1-e^{-\ad_X}\right)(Y)
$$
It follows that $e^{-X}e^{X+\epsilon Y}$ and $(1 + \epsilon \cdot (1-e^{-\ad_X})/\ad_X)(Y)$ differ by something commuting with $X$, i.e. there exists a formal power series $f(X)$ such that
$$e^{-X}e^{X+\epsilon Y} - \left(1 + \epsilon \cdot \frac{1-e^{-\ad_X}}{\ad_X}(Y)\right) = f(X).$$
Setting $Y=0$ in the last equality we find that $f(X) = 0$.
\end{proof}
\end{Lemma}
\end{proof}
\end{Prop}

\section{Equivariant Grothendieck-Riemann-Roch}\label{sect:eq_RR}
In this section we finally bring the results of previous sections together to give a proof of the classical Grothendieck-Riemann-Roch theorem as well as its equivariant analogue.
\subsection{Abstract GRR theorem}
We start by introducing the context we are interested in:
\begin{Def}\label{equiv_pushforward}
Let $(X, g_X)$ and $(Y, g_Y)$ be a pair of derived schemes with endomorphisms. An \emph{equivariant morphism $\xymatrix{(X, g_X) \ar[r]^-{f} & (Y, g_Y)}$} is a commutative diagram
$$\xymatrix{
X \ar[d]_-{f} \ar[r]^-{g_X} & X \ar[d]^-{f}
\\
Y \ar[r]_-{g_Y} & Y
}$$
where $\xymatrix{X \ar[r]^-{f} & Y}$ is some morphism of schemes. In this setting we will further denote by $\xymatrix{X^{g_X} \ar[r]^-{f^g} & Y^{g_Y}}$ the induced map on fixed points. 
\end{Def}

\begin{Rem}\label{push_of_equiv_str} Note that for a lax $g_X$-equivariant sheaf $\mathcal F \xymatrix{\ar[r]^t &} (g_{X})_* \mathcal F$ its pushforward $f_* \mcal{F} \in \QCoh(Y)$ to $Y$ automatically admits a $g_Y$-lax equivariant structure given by the composite
$$\xymatrix{
f_* \mathcal F \ar[r]^-{f_*(t)} & f_*(g_{X})_* \mathcal F \ar[r]_-{\sim} & (g_{Y})_* f_*\mathcal F.
}$$
We will further use the notation $f_*(\mathcal F, t)$ for $f_* \mathcal F \in \QCoh(Y)$ together with the above lax equivariant structure.
\end{Rem}

The definition above motivates the following
\begin{Def}
Let $(X,g)$ be a smooth scheme with an endomorphism. We will denote by $K_0^g(X)$ the usual $K_0$-group of the category of lax $g$-equivariant perfect sheaves on $X$.
\end{Def}

\begin{Rem} Note that for a proper equivariant morphism $\xymatrix{ (X,g_X) \ar[r]^-{f} & (Y,g_Y)}$ between smooth schemes\footnote{In fact it is enough to assume that $Y$ is smooth.} the induced pushforward functor $f_*$ is exact and preserves perfect sheaves (as $\Coh(Y) \simeq \QCoh(Y)^{\perf}$ and proper morphism preserves coherent sheaves) and therefore there is an induced morphism 
$$\xymatrix{
K_0^{g_X}(X) \ar[r] & K_0^{g_Y}(Y)
}$$
which we will also denote by $f_*$.
\end{Rem}

Motivated by Corollary \ref{morphism_of_Tr_is_classical_Todd} we also introduce the following

\begin{Not}\label{equiv_todd}
Given a smooth proper scheme together with an endomorphism $\xymatrix{X \ar[r]^-{g} & X}$ we will further denote the canonical orientation $\mathcal{O}_{X^g} \simeq \omega_{X^g}$ (see \Cref{FixedPointsAreCalabi}) on $X^g$ by $\td_g$ and call it an \emph{equivariant Todd distribution on $(X,g)$}.
\end{Not}

Now the abstract formalism of traces readily gives us the
\begin{Prop}[Abstract Grothendieck-Riemann-Roch]\label{thm:abstract_GRR}
Let $\xymatrix{(X, g_X) \ar[r]^-{f} & (Y, g_Y)}$ be an equivariant morphism between smooth proper schemes. Then the diagram
$$\xymatrix{
K_0^{g_X}(X) \ar[rr]^-{\ch(-,-)\td_{g_X}} \ar[d]_-{f_*} && \pi_0 \Gamma(X^{g_X}, \omega_{X^{g_X}}) \ar[d]^{(f^g)_*} 
\\
K_0^{g_Y}(Y) \ar[rr]_-{\ch(-,-)\td_{g_Y}} && \pi_0 \Gamma(Y^{g_Y}, \omega_{X^{g_Y}})
}$$
is commutative, i.e. for any perfect $g_X$-lax equivariant sheaf $(E,t)$ on $X$ there is an equality
$$
(f^g)_*(\ch(E,t)\td_{g_X}) = \ch(f_*(E,t))\td_{g_Y}
$$
in $\Gamma(Y^{g_Y},\omega_{Y^{g_Y}})$, where $\ch(E,t) \in \Gamma(X^{g_X},\mathcal{O}_{X^{g_X}})$ here is the categorical Chern character \ref{chern}.

\begin{proof}
By passing to the induced morphism of traces \ref{2trace} in the commutative diagram
$$\xymatrix{
\Vect_k \ar[rr]^-E && \QCoh(X)\ar@(ul,ur)^{g_{X*}} \ar[d]^{\otimes \mathcal O_X}_\sim \ar[rr]^-{f_*} && \QCoh(Y)\ar@(ul,ur)^{g_{Y*}} \ar[d]^-{\otimes \mathcal O_Y}_\sim 
\\
 && \ICoh(X)\ar@(dl,dr)_{g_{X*}} \ar[rr]_-{f_*} && \ICoh(Y)\ar@(dl,dr)_{g_{Y*}}
}$$
by functoriality we obtain a commutative diagram
$$\xymatrix{
k \ar[rr]_-{\ch(E,t)} \ar@/^2pc/[rrrr]^{\ch(f_*(E,t))} && \Gamma(X^{g_X}, \mathcal O_{X^{g_X}}) \ar[d]^{\cdot \td_{g_X}}_-{\sim} \ar[rr] && \Gamma(Y^{g_Y}, \mathcal O_{Y^{g_Y}}) \ar[d]^{\cdot \td_{g_Y}}_-{\sim}
\\
 && \Gamma(X^{g_X}, \omega_{X^{g_X}}) \ar[rr]_-{(f^g)_*} && \Gamma(Y^{g_Y}, \omega_{Y^{g_Y}})
}$$
where we use Corollaries \ref{TraceOfQCoh}, \ref{TraceOfICoh}, \ref{morph_of_tr_in_icoh} and \Cref{TraceIsFromCalabiYau_proof} to identify morphisms of traces in the above diagram.
\end{proof}
\end{Prop}

\subsection{Equivariant GRR theorem}
Now it may be hard to apply Proposition \ref{thm:abstract_GRR} in practice as in general we don't have a good description of $\Gamma(X^{g_X}, \mathcal O_{X^{g_X}})$ and of the Todd distribution $\td_{g_X}$. Fortunately, under some reasonable assumptions  one can use ideas of localization to express it in a more computable form.
\begin{assumption}\label{reasonable_eq_assumptions}
We will further assume that $(X,g)$ is a smooth scheme with an endomorphism $\xymatrix{X \ar[r]^-{g} & X}$ such that the reduced classical scheme $\overline{X^g} := \mathcal H^0(X^g)^{\red}$ is smooth (but not necessarily connected). We will denote by $\xymatrix{\overline{X^g} \ar[r]^-{j} & X}$ the canonical embedding and by $\mathcal N_g^\vee$ its conormal bundle. Note that the action of $g$ on $\Omega_X^1$ in particular restricts to an endomorphism $\xymatrix{\mathcal N^\vee_g \ar[r]^-{g^*_{|\mathcal N_g}} & \mathcal N^\vee_g}$. We will sometimes call $\overline{X^g}$ the \emph{classical fixed locus of $g$}.
\end{assumption}

Note that the embedding $\xymatrix{\overline{X^g} \ar@{^(->}[r]^-{j} & X}$ is equivariant with respect to the trivial equivariant structure on $\overline{X^g}$ and the given one on $X$ thus induces a morphism $j^g\colon \mathcal L\overline{X^g} \xymatrix{\ar[r] &} X^g$. The theorem below gives a criterion when $j^g$ is an equivalence:

\begin{Theor}[Localization theorem]\label{localization_thm}
Let $(X,g)$ be as in the previous notation. Then the pullback map $j^g$ is an equivalence if and only if the determinant $\det(1 - g^*_{|\mathcal N_g^\vee}) \in \Gamma(\overline{X^g}, \mathcal O_{\overline{X^g}})$ is an invertible function.

\begin{proof}
Since the map $\xymatrix{j^g\colon \mathcal L \overline{X^g} \ar[r] & X^g}$ is a nil-isomorphism, by \cite[Chapter 1, Proposition 8.3.2]{GaitsRozII} the morphism $j^g$ is an equivalence if and only if the induced map on cotangent spaces $\xymatrix{\alpha\colon (j^g)^* \mathbb L_{X^g} \ar[r] & \mathbb L_{\mathcal L\overline{X^g}}}$ is an equivalence. Moreover, since the inclusion $\xymatrix{\overline{X^g} \ar[r] & \mathcal L\overline{X^g}}$ is a nil-isomorphism, it is enough to prove that $\alpha_{|\overline{X^g}}$ is an equivalence.

Consider now the following commutative diagram of derived schemes
$$\xymatrix{
\overline{X^g} \ar[d]_-j \ar[r]^-\Delta & \overline{X^g} \times \overline{X^g} \ar[d]^-{j\times j} & \ar[l]_-\Delta \overline{X^g} \ar[d]^-j 
\\
X \ar[r]_-\Delta & X\times X & \ar[l]^-{(\Id_X, g)} X.
}$$
By definition the limit of the top row is $\mathcal L\overline{X^g}$, the limit of the bottom row is $X^g$ and $j^g$ is precisely the induced map on the limits. By applying the absolute cotangent complex functor and pulling everything back to $\mathcal L\overline{X^g}$ and then further pulling back along $\xymatrix{\overline{X^g} \ar[r] & \mathcal L\overline{X^g}}$ we then obtain a commutative diagram of sheaves
$$\xymatrix{
\Omega^1_{\overline{X^g}} & \ar[l]_-\nabla \Omega^1_{\overline{X^g}} \oplus \Omega^1_{\overline{X^g}} \ar[r]^-\nabla & \Omega^1_{\overline{X^g}}
 \\
\ar[u]^-{j^*} j^* \Omega^1_X & \ar[l]^-\nabla j^* \Omega^1_X\oplus j^* \Omega^1_X \ar[u]_{j^*\oplus j^*} \ar[r]_-{(1, g^*)} & j^* \Omega^1_X \ar[u]_{j^*}
}$$
(where by $\nabla$ we denote the codiagonal map) in $\QCoh(\overline{X^g})$. By the \Cref{cotangent_of_limit} below, the pushout of the top row is $(\mathbb L_{\mathcal L\overline{X^g}})_{|\overline{X^g}}$, the pushout of the bottom row is $(\mathbb L_{X^g})_{|\overline{X^g}}$ and the induced map between pushouts is $\alpha_{|\overline{X^g}}$. It follows $\alpha_{|\overline{X^g}}$ is an equivalence if and only if the pushout of fibers of the vertical maps
$$\xymatrix{
\mathcal N_g^\vee && \ar[ll]_-\nabla \mathcal N_g^\vee \oplus \mathcal N_g^\vee \ar[rr]^-{(1, g^*_{|\mathcal N_g^\vee})} && \mathcal N_g^\vee
}$$
is nullhomotopic. The above pushout may be computed as the cofiber of the map $\xymatrix{\mathcal N_g^\vee \oplus \mathcal N_g^\vee \ar[r] & \mathcal N_g^\vee \oplus \mathcal N_g^\vee}$ given in block-matrix form by
\begin{align}\label{infty_2_cats_and_matrices}
\begin{pmatrix}
1 & 1 \\
1 & g^*_{|\mathcal N_g^\vee}
\end{pmatrix}.
\end{align}
Finally, the matrix \eqref{infty_2_cats_and_matrices} is invertible if and only if $g^*_{|\mathcal N_g^\vee} - 1$ is invertible if and only if the determinant $\det(1 - g^*_{|\mathcal N_g^\vee})$ is invertible.
\end{proof}
\end{Theor}
\begin{Lemma}\label{cotangent_of_limit}
Let
$$\xymatrix{
X \ar[r]^p\ar[d]_q\ar[rd]^r & Y \ar[d] \\
Z \ar[r] & W
}$$
be a fibered square of derived schemes. Then the induced square
$$\xymatrix{
r^*\mathbb L_W \ar[r]\ar[d] & p^*\mathbb L_Y \ar[d] \\
q^* \mathbb L_Z \ar[r] & \mathbb L_X
}$$
is a pushout in $\QCoh(X)$.

\begin{proof}
Formal from the definition of $\QCoh$-cotangent complex and the universal property of the limit.\footnote{In fact, one can analogously prove the following more general statement: let $\xymatrix{X_\bullet: I \ar[r] & \PreStack}$ be a diagram of prestacks admitting cotangent complex. Then $X:=\lim_I X_i$ admits cotangent complex and the natural map $\xymatrix{\colim p_i^* \mathbb L_{X_i} \ar[r] & \mathbb L_X}$ is an equivalence where $\xymatrix{p_i: X \ar[r] & X_i}$ are the natural projections.}
\end{proof}
\end{Lemma}
\begin{Rem}
The theorem above tells us that if the determinant $\det(1-g^*_{|\mathcal N_g})$ is invertible (a condition which is often easy to verify in practice), then the ring $\Gamma(X^g,\mathcal{O}_{X^g})$ (which naturally appears in the abstract GRR-theorem \ref{thm:abstract_GRR}) is equivalent to a ring that we understand much better:
$$\xymatrix{
-_{|\mathcal{L}\overline{X^g}}\colon \Gamma(X^g,\mathcal{O}_{X^g}) \ar[r]^-{(j^g)^*}_-\sim & \Gamma(\mathcal{L}\overline{X^g}, \mathcal O_{\mathcal {L}\overline{X^g}}) \simeq \displaystyle \bigoplus_{p,q} H^{p,q}(\overline{X^g})[p-q].
}$$
Moreover, for a perfect $g$-equivariant sheaf $(E,t)$ on $X$ it is convenient to describe the categorical Chern character $\ch(E,t) \in \Gamma(X^g,\mathcal{O}_{X^g})$ in these terms: by Corollary \ref{prop:semi_eq_Chern} we have an equality
$$
\ch(E,t)_{|\mathcal{L}(\overline{X^g})}=\Tr_{\QCoh(\overline{X^g})} \left(\exp\left(\At(E_{|\overline{X^g}})\right) \circ t_{|\overline{X^g}}\right) \quad \in \quad \pi_0\Gamma(\overline{X^g}, \mathcal O_{\overline{X^g}}) \simeq \bigoplus_p H^{p,p}(\overline{X^g}).
$$
Since the rings $\Gamma(X^g,\mathcal{O}_{X^g})$ and $\Gamma(\mathcal{L}\overline{X^g}, \mathcal O_{\mathcal {L}\overline{X^g}})$ under assumptions of the localization theorem are canonically equivalent, by abuse of notations we will sometimes identify $\ch(E,t)$ with its image in $\pi_0\Gamma(\mathcal{L}\overline{X^g}, \mathcal O_{\mathcal {L}\overline{X^g}})$.
\end{Rem}
The conditions of \Cref{localization_thm} are sometimes automatically satisfied. To see this recall the following well-known lemma:
\begin{Lemma}
Let $X$ be a smooth variety over an algebraically closed field $k$ of characteristic zero and let $G$ be a reductive group acting on $X$. Then for each $G$-fixed point $x\in X$ one can choose a set of local coordinates $\{x_1,\ldots , x_n\}$ of $\mathcal O_{X,x} \simeq k[[x_1, \ldots, x_n]]$, such that $G$ acts linearly with respect to them. In particular, the fixed locus $X^G$ is smooth.

\begin{proof}
Let $\mathfrak m_{\mathcal O_{X,x}}$ be the maximal ideal of $\mathcal O_{X,x}$. Since the category of $G$-representations is semisimple, for each $n\ge 2$ the natural surjection $\xymatrix{\mathfrak m_{\mathcal O_{X,x}}/\mathfrak m_{\mathcal O_{X,x}}^n \ar@{->>}[r] & \mathfrak m_{\mathcal O_{X,x}}/\mathfrak m_{\mathcal O_{X,x}}^2}$ admits a $G$-equivariant section. Passing to the limit $n \to \infty$ we obtain a $G$-equivariant section $\xymatrix{s: \mathfrak m_{\mathcal O_{X,x}}/\mathfrak m_{\mathcal O_{X,x}}^2 \ar[r] & \mathfrak m_{\mathcal O_{X,x}}}$. Let $\{\overline{x}_1,\ldots \overline{x}_n\}$ be a basis of $\mathbb T_{X,x} \simeq \mathfrak m_{X,x}/\mathfrak m_{\mathcal O_{X,x}}^2$ and put $x_i := s(\overline{x}_i)$. By Nakayama's lemma $\{x_1, \ldots , x_n\}$ are local coordinates at the point $x$ and $G$ by construction acts linearly withe respect to them.
\end{proof}
\end{Lemma}
\begin{Cor}
Let $X$ be a smooth scheme equipped with an automorphism $g$ of finite order. Then the fixed locus $X^g$ is smooth and the natural map $\xymatrix{\mathcal L {X^g} \ar[r] & X^g}$ is an equivalence, i.e. the derived fixed locus $X^g$ is formal.
\end{Cor}
\begin{Rem}
This result for finite order automorphism was also obtained in \cite[Corollary 1.12]{ACH_formality} by different methods.
\end{Rem}

In order to proceed further we introduce the following
\begin{Def}[Euler classes]\label{def:euler_class}
Let $(X,g)$ be as above. Define an \emph{Euler class of $g$} as
$$
e_g :=(j^g)^*\ch\big{(}j_*(\mathcal O_{\overline{X^g}}, \Id_{\mathcal{O}_{\overline{X^g}}})\big{)} \quad \in \quad \pi_0\Gamma\big{(}\mathcal L\overline{X^g}, \mathcal O_{\mathcal L\overline{X^g}}\big{)},
$$
where $\xymatrix{\mathcal{L}\overline{X^g} \ar[r]^-{j^g} & X^{g}}$ is the induced map on fixed loci (see \Cref{push_of_equiv_str} for the notation $j_*(-,-)$).
\end{Def}
To describe the Euler class more explicitly, we recall the following standard result, the proof of which we include here for reader's convenience:
\begin{Lemma}
Let $\xymatrix{Z \ar@{^{(}->}[r]^-{i} & X}$ be a closed embedding of smooth schemes. Then there exists a canonical isomorphism
$$\mathcal H^{-k}(i^*i_*\mathcal O_Z) \simeq \Lambda^k(\mathcal N_{Z/X}^\vee)$$
of quasi-coherent sheaves on $Z$.

\begin{proof}
The case $k=0$ is obvious, since $i$ is a closed embedding. For $k=1$ note that by applying the pullback functor $i^*$ to the exact sequence
$$\xymatrix{
0 \ar[r] & \mathcal I_Z \ar[r] & \mathcal O_X \ar[r] & i_*\mathcal O_Z \ar[r] & 0
}$$
we obtain an isomorphism $\mathcal H^{-1}(i^*i_*\mathcal O_Z) \simeq \mathcal H^0(i^*\mathcal I_Z) \simeq \mathcal I_Z/\mathcal I_Z^2$. But by smoothness assumption $\mathcal I_Z/\mathcal I_Z^2$ is isomorphic to the conormal bundle $\mathcal N_{Z/X}^\vee$ of $Z$ in $X$.

Finally, the isomorphism $\mathcal N_{Z/X}^\vee \simeq \mathcal H^{-1}(i^*i_*\mathcal O_Z)$ and multiplication induce a map of algebras in $\QCoh(Z)$
$$\xymatrix{
\alpha^*\colon \Lambda^*(\mathcal N_{Z/X}) \ar[r] & \mathcal H^{-*}(i^*i_*\mathcal O_Z).
}$$
By smoothness assumption, $Z$ is a locally complete intersection in $X$, hence locally both parts are exterior algebras and thus since $\mathcal H^{-1}(\alpha)$ is an isomorphism so is $\alpha^*$.
\end{proof}
\end{Lemma}

\begin{Cor}\label{euler_as_sum}
We have
$$e_g = \ch\left(\Sym(\mathcal N_g^\vee[1]), \Sym(g^*_{|\mathcal N_g^\vee[1]})\right) = \sum_k (-1)^k \ch\left(\Lambda^k(\mathcal N_g^\vee), \Lambda^k(g^*_{|\mathcal N_g^\vee})\right).$$

\begin{proof}
By definition
$$
e_g = (j^g)^* \ch\big{(}j_*(\mathcal O_{\overline{X^g}}, \Id_{\mathcal O_{\overline{X^g}}})\big{)}  = \ch\big{(}j^*j_*(\mathcal O_{\overline{X^g}}, \Id_{\mathcal O_{\overline{X^g}}})\big{)}.
$$
Next, by the lemma above there is a (Postnikov) filtration on the complex $j^*j_* \mathcal O_{\overline{X^g}}$ with associated graded $\Sym(\mathcal N_g^\vee[1])$. The statement then follows from the fact that $\ch(-,-)$ is additive in fiber sequences.
\end{proof}
\end{Cor}

\smallskip

\begin{Cor}\label{localization_euler}
Let $(X,g)$ be as in \Cref{reasonable_eq_assumptions}. Then the following conditions are equivalent:
\begin{enumerate}
\item The morphism $1 - (g^*)_{|\mathcal N_{g}^\vee}$ is invertible.

\item The induced map $\xymatrix{j^g: \mathcal L\overline{X^g} \ar[r] & X^g}$ is an equivalence.

\item The Euler class $e_g$ is invertible.
\end{enumerate}

\begin{proof}
The equivalence of the first two assertions is the content of \Cref{localization_thm}. To see that the first and the last conditions are equivalent, note that since the natural inclusion $\xymatrix{\overline{X^g} \ar[r] & \mathcal L\overline{X^g}}$ is a nil-isomorphism, the Euler class $e_g$ is invertible if and only if its zero term $e_{g,0}$ is. But by \Cref{euler_as_sum} and \Cref{ex:ch0} we get
$$e_{g,0} = \sum_{k=0} (-1)^k \ch_0(\Lambda^k(\mathcal N_g^\vee), g^*_{|\Lambda^k(\mathcal N_g^\vee)}) = \sum_k (-1)^k\Tr(g^*_{|\Lambda^k(\mathcal N_g^\vee)})=\det(1- g^*_{|\mathcal N_g^\vee}).$$
\end{proof}
\end{Cor}

\smallskip  Using the Euler class and corollary above, we can describe the Todd distribution $\td_g$ (see \Cref{equiv_todd}) in more concrete terms
\begin{Prop}\label{eq_Todd_in_classical_terms}
Let $(X,g)$ be as in \ref{reasonable_eq_assumptions} and assume that the Euler class $e_g$ is invertible (and so the natural morphism $\xymatrix{j^g: \mathcal L\overline{X^g} \ar[r] & X^g}$ is an equivalence). Then under the composite equivalence
$$
\pi_0 \Gamma(X^g, \omega_{X^g}) \xymatrix{\ar[r]^-{(j^g)^*}_-\sim & } \pi_0 \Gamma(\mathcal L\overline{X^g}, \omega_{\mathcal L\overline{X^g}}) \simeq \bigoplus_p H^{p,p}(\overline{X^g})
$$
the Todd distribution $\td_g \in \Gamma(X^g, \omega_{X^g})$ corresponds to $\frac{\td_{\overline{X^g}}}{e_g}$, where $\td_{\overline{X^g}}$ is the ordinary Todd class.

\begin{proof}
By applying \Cref{thm:abstract_GRR} to the canonical inclusion $\xymatrix{\overline{X^g} \ar[r]^-{j} & X}$ we obtain
$$
(j^g)_*(\td_{\overline{X^g}}) = (j^g)_*\big{(}\ch(\mathcal O_{\overline{X^g}}, \Id_{\mcal{O}_{\overline{X^g}}}) \td_{\overline{X^g}}\big{)} = \ch\big{(}j_*(\mathcal O_{\overline{X^g}}, \Id_{\mcal{O}_{\overline{X^g}}})\big{)}\td_g.
$$
Consequently, by pulling back along $j^g$ (and using that $(j^g)^*(j^g)_*$ is identity) we obtain
$$
\td_{\overline{X^g}} = (j^g)^*(j_g)_*(\td_{\overline{X^g}}) = (j^g)^*\big{(}\ch(j_*(\mathcal O_{\overline{X^g}}, \Id_{\mcal{O}_{\overline{X^g}}}))\td_g\big{)} = e_g \cdot (j^g)^*(\td_g).
$$
We conclude by dividing both parts by $e_g$.
\end{proof}
\end{Prop}

As a corollary we obtain:

\begin{Theor}[Equivariant Grothendieck-Riemann-Roch]\label{equiv_GRR}
Let $\xymatrix{(X,g_X) \ar[r]^-{f} & (Y,g_Y)}$ be an equivariant morphism between smooth proper schemes such that
\begin{itemize}
\item Reduced fixed loci $\overline{X^{g_X}}$ and $\overline{Y^{g_Y}}$ are smooth.

\item The induced morphisms on conormal bundles $1 - (g_X^*)_{|\mathcal N_{g_X}^\vee}$ and $1 - (g_Y^*)_{|\mathcal N_{g_Y}^\vee}$ are invertible.
\end{itemize}
Then for a perfect lax $g_X$-equivariant sheaf $(E,t)$ on $X$ we have an equality
$$
(\overline{f^g})_* \left(\ch(E,t)\frac{\td_{\overline{X^{g_X}}}}{e_{g_X}}\right) = \ch\big{(}f_*(E,t)\big{)}\frac{\td_{\overline{Y^{g_Y}}}}{e_{g_Y}}
$$
in $\bigoplus_p H^{p,p}(\overline{Y^{g_Y}})$.
\begin{proof}
This follows immediately from the abstract Grothendieck-Riemann-Roch \ref{thm:abstract_GRR} and the identification of $\td_g$ above.
\end{proof}
\end{Theor}

Specializing to the case when $Y=\ast$, we get
\begin{Cor}[Equivariant Hirzebruch-Riemann-Roch]
Let $(X,g)$ be as in the theorem above. Then for any lax $g$-equivariant perfect sheaf $(E, t)$ on $X$ we have
$$
\int\limits_{\overline{X^g}} \ch(E,t) \frac{\td_{\overline{X^g}}}{e_g}=\Tr_{\Vect_k} \Gamma(X,t).
$$
\end{Cor}

Specializing even further we recover
\begin{Cor}[Holomorphic Atiyah-Bott fixed point formula]\label{cor:atiayh_bott}
Assume that the graph of $g$ intersects the diagonal in $X\times X$ transversely. Then
$$
\sum_{x=g(x)} \frac{\Tr(g^*_{|E_x})}{\det(1-d_x g)}=\Tr_{\Vect_k} \Gamma(X,t).
$$

\begin{proof}
By assumption on $g$ the derived fixed locus $X^g$ is discrete, hence $\mathcal L\overline{X^g} \simeq X^g$ and the corollary above reads as
$$
\sum_{x=g(x)} \frac{\ch(E,t)_{|X^g}}{e_g}=\Tr_{\Vect_k} \Gamma(X,t).
$$
Let $x \in X^g$ be a fixed point of $g$. Since $X^g$ is discrete, for any perfect lax $g$-equivariant sheaf $\At(E_{|\overline{X^g}}) \simeq 0$, hence by Corollary \ref{prop:semi_eq_Chern} we have
$$\ch(E,t)_x = x^* \Tr_{\QCoh(X^g)}(g^*_{|E_{|X^g}}) = \Tr(g^*_{|E_x}).$$
Moreover, the conormal bundle $\mathcal N^\vee_g$ in this case is just the cotangent space at fixed points, hence by the vanishing of the Atiyah class on $X^g$ we deduce
$$(e_g)_x = (e_{g,0})_x = \det(1-d_x g).$$
\end{proof}
\end{Cor}

Specializing in the other direction, we recover
\begin{Cor}[Grothendieck-Riemann-Roch]\label{cor:GRR}
Let $X,Y$ be smooth proper schemes and $\xymatrix{X \ar[r]^-{f} & Y}$ be a morphism. Then for any perfect sheaf $E$ on $X$ we have
$$f_*(\ch(E)\td_X) = \ch(f_*E)\td_Y$$
where above $\ch(E)$ and $\ch(f_* E)$ are the classical Chern characters.
\begin{proof}
Consider the \Cref{equiv_GRR} in the case when $g_X, g_Y$, are morphisms. Note that we have $\overline{X^g} = X$ and the map
$$\xymatrix{
\mathcal LX  \simeq \mathcal L\overline{X^g} \ar[r]^-{j^g} & X^g \simeq \mathcal LX
}$$
(and analogously for $Y$) is tautologically equivalent to the identity. Moreover, since the conormal bundles $\mathcal N^\vee_{g_X}$ and $\mathcal N^\vee_{g_Y}$ are in this case trivial,  we have $e_{g_X} = 1$ and $e_{g_Y} = 1$. Hence by \Cref{equiv_GRR} we obtain
$$f_*(\ch(E, \Id_E)\td_X) = \ch(f_*(E, \Id_E))\td_Y.$$
It is left to note that by \Cref{ChernIsChern} the categorical Chern character $\ch(E,\Id_E) = \ch(E)$ coincides with the classical one. 
\end{proof}
\end{Cor}

\appendix
\section{Appendix: reminder on formal deformation theory}
\label{app:formal_moduli_reminder}
In this section we review main results of formal deformation theory developed in \cite[Chapters 5-9]{GaitsRozII} relevant to this work (hence we need to work over a field of characteristic zero). We start with the following
\begin{Def}(\cite[Chapter 5, Definition 1.1.1]{GaitsRozII})
Define an $\ii$-category of \emph{formal moduli problems over $X$} denoted by $\FormModuli_{/X}$ as the full subcategory of $\PreStack_{\laft/X}$ (see (\cite[Chapter 2, 1.6]{GaitsRozI}) on objects $\xymatrix{Y \ar[r]^-{p} & X}$ such that:
\begin{itemize}
\item The morphism $p$ is inf-schematic (\cite[Chapter 2, Definition 3.1.5]{GaitsRozII}).

\item The morphism $p$ is nil-isomorphism, i.e. the induced morphism $\xymatrix{^{red}Y \ar[rr]^-{^{red}p} && ^{red} X}$ is an equivalence.
\end{itemize}

\noindent Group objects in the category of formal moduli problems over $X$ are called \emph{formal groups over $X$}.
\end{Def}

Let now $\xymatrix{\mcal{Y} \ar[r]^-{p} & X}$ be a formal moduli problem over $X$. The functor $p_*$ being left adjoint to a symmetric monoidal functor $f^!$ is left lax-monoidal. Hence $p_*\omega_{\mathcal Y}$ is naturally a cocommutative coalgebra object of $\ICoh(X)$, that is, an object of $\coCAlg(\ICoh(X))$. Moreover since the functor 
$$\xymatrix{
\FormModuli_{/X} \ar[r] & \coCAlg(\ICoh(X))
}$$
$$\xymatrix{
\mcal{Y} \ar@{|->}[r] & p_* \omega_{\mcal{Y}}
}$$
is symmetric monoidal (e.g. $p_*\omega_{\mathcal Y\times_X \mathcal Z} \simeq p_*\omega_{\mathcal Y}\otimes p_*\omega_{\mathcal Z}$) for a formal group $\widehat G \in \Grp(\FormModuli_{/X})$ the sheaf $p_*\omega_{\widehat G}$ is a group object in the category of cocommutative coalgebras, i.e. a cocommutative Hopf algebra. In particular, we can define a functor
$$\xymatrix{
 \Grp(\FormModuli_{/X}) \ar[rr]^-{\Lie_X} && \LAlg(\ICoh(X))
 }$$
by setting
$$\Lie_X(\widehat G) := \Prim(p_*\omega_{\widehat G}) \in \LAlg(\ICoh(X))$$
for a formal group $\widehat G \in  \Grp(\FormModuli_{/X})$, where $\xymatrix{\HopfAlg(\ICoh(X)) \ar[r]^-{\Prim} & \LAlg(\ICoh(X))}$ is the functor of primitive elements.

\smallskip

Now the following crucial theorem relates groups in the category of formal moduli problems over $X$ and Lie algebras in the category of quasi-coherent sheaves on $X$:
\begin{Theor}[{\cite[Chapter 7, Theorem 3.6.2 and Proposition 5.1.2]{GaitsRozII}}]\label{LieAlgAndFormGroups}
We have:
\begin{enumerate}[1.]
\item There is an equivalence of $\ii$-categories
$$\xymatrix{
 \Grp(\FormModuli_{/X}) \ar[rr]^-{\Lie_X}_-{\sim} && \LAlg(\ICoh(X)),
 }$$
where $\LAlg(\ICoh(X))$ is the $(\infty,1)$-category of algebras in $\ICoh(X)$ over the Lie operad. Moreover, for a formal group $\widehat G\in \Grp(\FormModuli_{/X})$ the underlying ind-coherent sheaf of $\Lie_X(\widehat G) \in \LAlg(\ICoh(X))$ is equivalent to $\mathbb T_{\widehat G/X, e} := e^! \mathbb T_{\widehat G/X}$, where $\xymatrix{X \ar[r]^-{e} & \widehat G}$ is the identity section and $\mathbb{T}$ denotes tangent sheaf.

\item For $\widehat G \in  \Grp(\FormModuli_{/X})$ there is an equivalence of $\ii$-categories
$$\xymatrix{
\Rep_{\widehat G}(\ICoh(X)) \ar[r]_-{\sim} & \Mod_{\Lie_X( \widehat G)}(\ICoh(X)).
}$$
\end{enumerate}
\end{Theor}

\smallskip 

Now in classical theory of Lie groups for a (real) Lie group $G$ with Lie algebra $\mathfrak g$ there is an exponential map $\xymatrix{\mathfrak g \ar[r]^-{\exp_G} & G}$, which is a diffeomorphism in a small enough neighborhoods of $0\in\mathfrak g$ and $1_G\in G$. The same story works even better in the formal world since one does not need to consider neighborhoods. In order to formulate this statement explicitly in our setting, we first need the following
\begin{Def}\label{vectorprestack}
For $E \in \ICoh(X)$ define a \emph{vector prestack $\mathbb V(E) \in \FormModuli_{X//X}$ of $E$} by the property that for any $\mcal{Y} \in \FormModuli_{X//X}$ there is a natural equivalence
$$
\Hom_{\FormModuli_{X//X}}(\mathcal Y, \mathbb V(E)) := \Hom_{\QCoh(X)}(I(p_* \omega_{\mathcal Y}), E),
$$
where for a coaugmented coalgebra $C\in \coCAlg^{\sf coaug}(\QCoh(X))$ we let $\xymatrix{I(C):=\cofib(\mcal{O}_X \ar[r] & C)}$ to be the coaugmentation ideal of $C$.
\end{Def}

\begin{Rem}\label{specinf}
In fact vector prestack can be seen as a part of a more general construction. Namely, define a functor
$$\xymatrix{
\coCAlg(\ICoh(X)) \ar[rr]^-{\Spec^{\sf inf}_{X}} && \FormModuli_{/X}
}$$
as the right adjoint to the functor $\xymatrix{\mathcal Y \ar@{|->}[r] & p_*\omega_{\mathcal Y}}$ so that for any formal moduli problem $\mathcal Y \in \FormModuli_{/X}$ we have an equivalence
$$
\Hom_{\FormModuli_{/X}}(\mathcal Y, \Spec^{\sf inf}_{X} C) \simeq \Hom_{\coCAlg(\ICoh(X))}(p_*\omega_{\mathcal Y},C).
$$
Then it is straightforward to see that for a sheaf $E \in \ICoh(X)$ there is a natural equivalence $\mathbb V(E) \simeq \Spec^{\sf inf}_{X} \left( \Sym(E) \right)$, where
$$\xymatrix{
\ICoh(X) \ar[r]^-{\Sym} & \coCAlg^{\sf nil}(\ICoh(X))} \qquad E \mapsto \bigoplus_{n=0}^\infty (E^{\otimes n})_{\Sigma_n}
$$
is the symmetric algebra functor, endowed with its canonical cofree ind-nilpotent co commutative coalgebra structure. We refer interested reader to \cite[Chapter 7, 1.3]{GaitsRozII} for through discussion of the inf-spectrum functor $\Spec^{\sf inf}$.
\end{Rem}

\begin{Ex}\label{ex:vector_stacks_vs_rel_spec}
Unwinding the definitions one finds that for $E\in \QCoh(X)$ such that $E^\vee \in \Coh^{<0}(X)$ the vector prestack $\mathbb V(\Upsilon_X(E))$ is equivalent to the "vector bundle associated to $E$", i.e.
$$\mathbb V(\Upsilon_X(E)) \simeq \Spec_{/X}(\Sym_{\QCoh(X)}(E^\vee)).$$
In the case when $E^\vee \in \Coh^{\le 0}(X)$, there is a similar equivalence if we take the formal completion at the zero section of the right-hand side.
\end{Ex}

\smallskip In these notations we finally have
\begin{Theor}[{\cite[Chapter 7, Corollary 3.2.2.]{GaitsRozII}}]\label{formal_exponent}
Let $\widehat G \in \Grp(\FormModuli_{/X})$ be a formal group over $X$. Then there is a functorial equivalence
$$\xymatrix{
\mathbb V(\Lie_X(\widehat G)) \ar[rr]^-{\exp_{\widehat G}}_-{\sim} && \widehat G
}$$
of formal moduli problems over $X$.

\begin{proof}[Idea of the proof]
Given a Lie algebra $\mathfrak g\in \ICoh(X)$ its universal enveloping algebra is naturally a cocommutative Hopf algebra, i.e. a group object in the category $\coCAlg(\ICoh(X))$. Since the functor $\Spec^{\sf inf}$ is monoidal, we see that $\exp_X(\mathfrak g) := \Spec^{\sf inf}\left(U(\mathfrak g)\right)$ is a group object in $\FormModuli_{/X}$. In fact, one can show that the construction $\xymatrix{\mathfrak g \ar@{|->}[r] & \exp_X(\mathfrak g)}$ is the inverse to the $\Lie_X$ functor from Theorem \ref{LieAlgAndFormGroups}. But by \cite[Chapter 6, Corollary 1.7.3]{GaitsRozII} there is canonical equivalence of cocommutative coalgebras $U(\mathfrak g) \simeq \Sym_{\ICoh(X)}(\mathfrak g)$ (aka Milnor-Moore theorem), hence for $\mathfrak g = \Lie_X(\widehat G)$ we have
$$
\widehat G \simeq \exp_X(\mathfrak g) = \Spec^{\sf inf}\left(U(\mathfrak g)\right) \simeq \Spec^{\sf inf}\left(\Sym_{\ICoh(X)}(\mathfrak g)\right) = \mathbb V(\mathfrak g).
$$
\end{proof}
\end{Theor}
\begin{Rem}\label{exp_for_abelian}
Notice that for a formal group $\widehat A \in \Grp(\FormModuli_{/X})$ with abelian Lie algebra the map $\exp_{\widehat A}$ above is not only an equivalence of formal moduli problems, but moreover an equivalence of formal groups. For example, in the case $\widehat G:=\widehat{\mathbb G_m}$ the map 
$$\xymatrix{
\widehat{\mathbb G_a}\simeq \mathbb V(\Lie_X(\widehat{\mathbb G_m})) \ar[rr]^-{\exp_{\widehat{\mathbb G_m}}}_-{\sim} && \widehat{\mathbb G_m}
}$$
is the usual formal exponent.
\end{Rem}

\section{Appendix: correspondences and traces}
\label{app:corr}
\subsection{Ind-coherent sheaves and morphism of traces}
In this section we discuss how one can calculate the morphism of traces in the setting of ind-coherent sheaves using the category of correspondences. We start with the following

\begin{Theor}(\cite[Chapter 5, Theorem 2.1.4., Theorem 4.1.2]{GaitsRozI})\label{qcoh}
The ind-coherent sheaves functor can be lifted to a symmetric monoidal functor
$$\xymatrix{
\Corr(\Sch_{\aft})^{\proper} \ar[rr]^-{\ICoh} && 2\Cat_k.
}$$
Where $\Corr(\Sch_{\aft})^{\proper}$ is the $(\infty,2)$-category which can be informally described as follows:
\begin{enumerate}
\item Its objects are those of $\Sch_{\aft}$.
\item Given $X,Y \in \Sch_{\aft}$ a morphism from $X$ to $Y$ in $\Corr(\Sch_{\aft})^{\proper}$ is given by a span
$$\xymatrix{
X & W \ar[r]^-{f} \ar[l]_-{g} & Y
}$$
and the composition of morphisms is given by taking pullbacks.
\item Given two morphisms $W_1,W_2 \in \Hom_{\Corr(\Sch_{\aft})^{\proper}}(X,Y)$ a $2$-morphism from $W_1$ to $W_2$ is given by a commutative diagram
$$\xymatrix{
& W_1 \ar[dr] \ar[dl] \ar[dd]^-{h}
\\
X && Y
\\
& W_2 \ar[ur] \ar[ul]
}$$
in $\Sch_{\aft}$ where $h$ is proper.
\end{enumerate}
The symmetric monoidal structure on $\Corr(\Sch_{\aft})^{\proper}$ is given by the cartesian product of underlying objects of $\Sch_{\aft}$. Once again, we refer to \cite[Part III]{GaitsRozI} for a discussion of the category of correspondences. 
\end{Theor}
\begin{Rem}
In \cite[Chapter 7]{GaitsRozI} the category $\Corr(\Sch_{\aft})^{\proper}$ was denoted by $\Corr(\Sch_{\aft})^{\proper}_{\all,\all}$.
\end{Rem}
Informally speaking, the functor above maps $X \in \Sch_{\aft}$ to the category $\ICoh(X) \in \Cat_k$ and a morphism $\xymatrix{X & W \ar[r]^-{f} \ar[l]_-{g} & Y}$ in $\Corr(\Sch_{\aft})^{\proper}$ to the composite $\xymatrix{\ICoh(X) \ar[r]^-{g^!} & \ICoh(W) \ar[r]^-{f_*} & \ICoh(Y)}$ in $2\Cat_k$.

\smallskip

\begin{Rem}\label{magic_of_corr_icoh}
Since by (\cite[Chapter 9, Proposition 2.3.4.]{GaitsRozI}) every object $X \in \Corr(\Sch_{\aft})^{\proper}$ is self-dual via the morphisms
$$\xymatrix{
\ast && \ar[ll]_-{p} X \ar[rr]^-{\Delta} && X \times X
}$$
and
$$\xymatrix{
X \times X && \ar[ll]_-{\Delta} X \ar[rr]^-{p} && \ast
}$$
we see that the category $\ICoh(X) \in 2\Cat_k$ is also self-dual. Moreover, note that by (\cite[Chapter 9, Proposition 2.3.4.]{GaitsRozI}) every morphism in $\Corr(\Sch_{\aft})^{\proper}$ of the form
$$\xymatrix{
X && X \ar[ll]_-{\Id_X} \ar[rr]^-{f} && Y
}$$
where $f$ is proper admits a right adjoint given by 
$$\xymatrix{
Y && X \ar[rr]^-{\Id_X} \ar[ll]_-{f} && X
}$$
\end{Rem}

\smallskip

Our goal now is to understand morphism of traces in the setting of correspondences. We start by calculating classical traces in the category of correspondences
\begin{Prop}\label{Tr_in_corr}
The trace of the endomorphism $\xymatrix{X & Y \ar[l]_-{g} \ar[r]^-{f} & X}$ in $\Corr(\Sch_{\aft})^{\proper}$ is given by
$$\xymatrix{\ast & X^{f=g} \ar[l] \ar[r] & \ast},$$
where $X^{f=g}$ is defined as the pullback
$$\xymatrix{
X^{f=g} \ar[r]^-{i} \ar[d]_-{j} & Y \ar[d]^{(f,g)}
 \\
X \ar[r]_-{\Delta} & X\times X.
}$$

\begin{proof} By definition, the trace is given by the composite
$$\xymatrix{
& \ar[dl]_-{p} X \ar[dr]^-{\Delta} && \ar[dl]_-{g \times \Id_X} Y \times X \ar[dr]^-{f \times \Id_X} && X \ar[dl]_-{\Delta} \ar[dr]^-{p}
\\
\ast && X \times X && X \times X && \ast.
}$$
Since the composition in $\Corr(\Sch_{\aft})^{\proper}$ is given by taking pullback, the result follows.
\end{proof}
\end{Prop}

\smallskip

\begin{Cor} \label{tr_of_corr_icoh} Applying the functor $\xymatrix{\Corr(\Sch_{\aft})^{\proper} \ar[rr]^-{\ICoh} && 2\Cat_k}$ we see that the trace of the endomorphism $\xymatrix{\ICoh(X) \ar[r]^-{f_* g^!} & \ICoh(X)}$ in $2\Cat_k$ is given by $\Gamma(X^{f=g},\omega^{\ICoh}_{X^{f=g}})$.
\end{Cor}

\smallskip

We are now going to understand morphism of traces in the setting of correspondences (and therefore in the setting of ind-coherent sheaves). In order to simplify notation we will denote a morphism 
$$\xymatrix{
X && W \ar[ll]_-{g} \ar[rr]^-{f} && Y
}$$
by $\la \prescript{g}{X}{W}_Y^f\ra$. Here is the main

\smallskip

\begin{Prop}\label{morph_of_tr_in_corr}
Given a (not necessary commutative) diagram
$$\xymatrix{
X \ar[dd]_-{\la \prescript{\Id_X}{X} X_U^{s} \ra} \ar[rr]^-{\la \prescript{g}{X} Y_X^f \ra} && X \ar[dd]^-{\la \prescript{\Id_X}{X} X_U^{s} \ra} \ar@2[ddll]_-{T}
\\
\\
U  \ar[rr]_-{\la \prescript{b}{U} V_U^a \ra} && U 
}$$
in $\Corr(\Sch_{\aft})^{\proper}$ where the $2$-morphism $T$ is given by a choice of some commutative diagram
$$\xymatrix{
& Y \ar[dr]^-{s \circ f} \ar[dl]_-{g} \ar[dd]^-{t}
\\
X && U
\\
& X ~^{s}\!\times^{b}_U V \ar[ur]_-{a \circ \pr_V} \ar[ul]^-{\pr_X}
}$$
$$
\alpha_1: b \circ \pr_V \circ t \simeq s \circ \pr_X \circ t 
$$ 
$$
\alpha_2 :\pr_X \circ t \simeq g
$$
$$
\alpha_3: a \circ \pr_V \circ t \simeq s \circ f
$$
where $t \in \Hom_{\Sch_{\aft}}(Y, X ~^{s}\!\times^{b}_U V)$ is proper, the induced morphism of traces
$$\xymatrix{
X^{f=g} \simeq \Tr_{\Corr(\Sch_{\aft})^{\proper}}\big{(}\la \prescript{g}{X} Y_X^f \ra \big{)} \ar[rrr]^-{ \Tr\big{(}\la \prescript{\Id_X}{X} X_U^{s} \ra,T\big{)}} &&& \Tr_{\Corr(\Sch_{\aft})^{\proper}}\big{(}\la \prescript{b}{U} V_U^a \ra \big{)} \simeq U^{a=b}
}$$
is obtained as the map of pullbacks from the commutative diagram
$$\xymatrix{
X^{f=g}  \ar[dd] \ar[dr] \ar[rrr] &&& X^{f=g} \ar[dd]^(.70){s \circ i_X} \ar[dr]
\\
&  X^{f=g} \ar[dd]^(.70){\pr_V \circ t \circ j_X} \ar[rrr] &&& X^{f=g} \ar[dd]^(.70){(s \circ i_X, s \circ i_X)}
\\
U^{a=b} \ar[rrr]^-{i_U} \ar[dr]_-{j_U} &&& U \ar[dr]^{\Delta}
\\
& V \ar[rrr]^-{(b,a)} &&& U \times U,
}$$
where the commutativity of the front square is given by the equivalences
$$
b \circ \pr_V \circ t \circ j_X \stackrel{\alpha_1 \circ j_X}{\simeq} s \circ \pr_X \circ t \circ j_X \stackrel{s \circ \alpha_2 \circ j_X}{\simeq} s \circ g \circ j_X \simeq s \circ i_X
$$
and
$$
a \circ \pr_V \circ t \circ j_X \stackrel{\alpha_3 \circ j_X}{\simeq} s \circ f \circ j_X \simeq s \circ i_X.
$$
\begin{proof}
Unwinding the definitions, we see that the morphism of traces is induced by the diagram
$$\xymatrix{
\ast \ar[dd]_-{\Id_{\ast}} \ar[rrr]^-{\la \prescript{}{\ast}X^{\Delta}_{X \times X} \ra} &&& X \times X \ar@2[ddlll]  \ar[rrr]^-{\la \prescript{g \times \Id_X}{X \times X}{Y \times X}^{f \times \Id_X}_{X \times X} \ra} \ar[dd] &&& X \times X \ar[rrr]^-{\la \prescript{\Delta}{X \times X} X_{\ast} \ra} \ar[dd] \ar@2[ddlll] &&& \ast \ar[dd]^-{\Id_{\ast}} \ar@2[ddlll] 
\\
\\
\ast \ar[rrr]_-{\la \prescript{}{\ast}U^{\Delta}_{U \times U} \ra} &&& U \times U \ar[rrr]_-{\la \prescript{b \times \Id_U}{U \times U}{V \times U}^{a \times \Id_U}_{U \times U} \ra} &&& U \times U \ar[rrr]_-{\la \prescript{\Delta}{U \times U} U_{\ast}\ra} &&& \ast
}$$
where 
\begin{enumerate}
\item The two middle vertical morphisms are given by ${\la \prescript{\Id_X \times \Id_X}{X \times X}X \times X^{s \times s}_{U \times U} \ra}$ (see remark \ref{magic_of_corr_icoh}).
\item The left square is induced by the $2$-morphism
$$\xymatrix{
& X \ar[dr]^-{(s,s)} \ar[dl] \ar[dd]^-{s}
\\
\ast && U \times U
\\
& U \ar[ur]_-{\Delta} \ar[ul]
}$$
in $\Corr(\Sch_{\aft})^{\proper}$.
\item The middle square is induced by the $2$-morphism
$$\xymatrix{
& Y \times X \ar[dr]^-{(s \circ f) \times s} \ar[dl]_-{g \times \Id_X} \ar[dd]^-{t \times \Id_X}
\\
X \times X && U \times U
\\
& (X ~^{s}\!\times^{b}_U V) \times X \ar[ur]_-{(a \circ \pr_V) \times s} \ar[ul]^-{\pr_X \times \Id_X}
}$$
in $\Corr(\Sch_{\aft})^{\proper}$.
\item The right square is induced by the $2$-morphism
$$\xymatrix{
& X \ar[dr] \ar[dl]_-{\Delta} \ar[dd]
\\
X \times X  && \ast
\\
& X \times_U X \ar[ur] \ar[ul]
}$$
in $\Corr(\Sch_{\aft})^{\proper}$.
\end{enumerate}
Consequently, we see that the whole morphism of traces is given by the composite of the morphisms
$$\xymatrix{
Y ~^{(f,g)}\!\times^{\Delta}_{(X \times X)} X \ar[r] & Y ~^{(f,g)}\!\times_{(X \times X)} (X \times_U X) \simeq  X ~^{\Delta}\!\times_{(X \times X)}^{g \times \Id_X} (Y \times X) ~^{{(s \circ f) \times s}}\!\times_{U \times U}^{\Delta} U 
}$$
$$\xymatrix{
X ~^{\Delta}\!\times_{(X \times X)}^{g \times \Id_X} (Y \times X) ~^{{(s \circ f) \times s}}\!\times_{U \times U}^{\Delta} U  \ar[r] & X ~^{\Delta}\!\times_{(X \times X)}^{\pr_X \times \Id_X} \big{(} (X \times_V Y) \times X \big{)} ~^{(a \circ \pr_V) \times s}\!\times_{U \times U}^{\Delta} U 
}$$
and
$$\xymatrix{
X ~^{\Delta}\!\times_{(X \times X)}^{\pr_X \times \Id_X} \big{(} (X \times_V Y) \times X \big{)} ~^{(a \circ \pr_V) \times s}\!\times_{U \times U}^{\Delta} U \simeq X~^{(s,s)}\!\times_{(U \times U)}^{(b,a)} V \ar[r] & U~^{\Delta}\!\times_{(U \times U)}^{(b,a)} V.
}$$
In particular, we can rewrite the whole composition as the left vertical morphism from $X^{f=g}$ to $U^{a=b}$ in the commutative diagram
$$
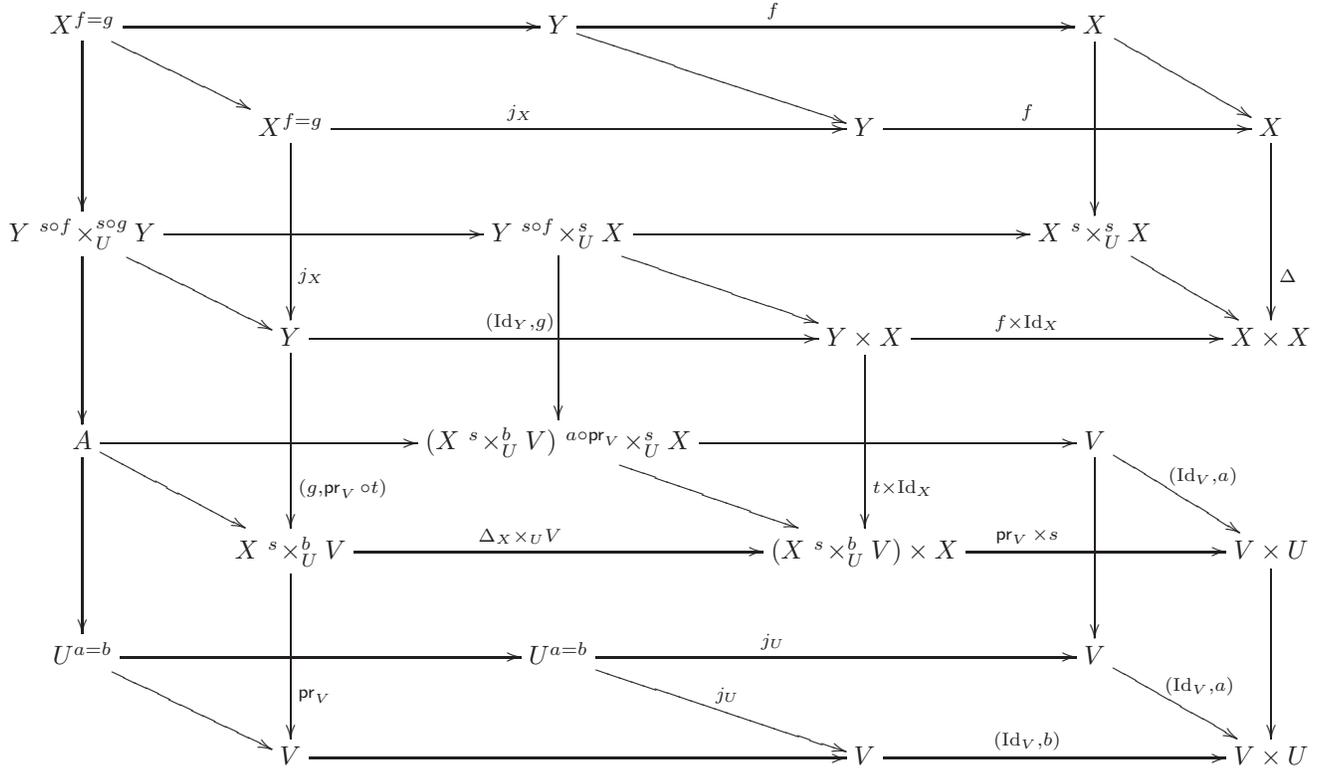
\xymatrix{
X^{f=g} \ar[dd] \ar[rr] \ar[dr] && Y \ar[dr] \ar[rr]^(.40){f} && X \ar[dr] \ar[dd]
\\
& X^{f=g} \ar[rr]^(.40){j_X} \ar[dd]^(.70){j_X}  && Y  \ar[rr]^(.40){f} && X \ar[dd]^(.70){\Delta}
\\
Y ~^{s \circ f}\!\times_{U}^{s \circ g} Y \ar[dd] \ar[rr] \ar[dr] && Y~^{s \circ f}\!\times_{U}^{s} X \ar[dd] \ar[dr] \ar[rr] && X~^{s}\!\times_{U}^{s} X \ar[dr] 
\\
& Y  \ar[rr]^(.40){(\Id_Y,g)} \ar[dd]^(.70){(g, \pr_V \circ t)}  && Y \times X \ar[dd]^(.70){t \times \Id_X} \ar[rr]^(.40){f \times \Id_X} && X \times X 
\\
A \ar[dd] \ar[rr] \ar[dr] && (X~^{s}\!\times_{U}^{b} V) ~^{a \circ \pr_V }\!\times_{U}^{s} X   \ar[dr] \ar[rr] && V \ar[dr]^-{(\Id_V,a)} \ar[dd]
\\
& X~^{s}\!\times_{U}^{b} V \ar[rr]^(.40){\Delta_X \times_U V} \ar[dd]^(.70){\pr_V}  && (X~^{s}\!\times_{U}^{b} V) \times X \ar[rr]^(.40){\pr_V \times s} && V \times U \ar[dd]
\\
U^{a=b} \ar[rr] \ar[dr] && U^{a=b} \ar[dr]^-{j_U} \ar[rr]^(.40){j_U} && V \ar[dr]^-{(\Id_V,a)}
\\
& V \ar[rr] && V \ar[rr]^(.40){(\Id_V,b)} && V \times U
}$$
where all the horizontal squares are pullback squares. Consequently, we see that we can rewrite the morphism of traces as the left vertical morphism from $X^{f=g}$ to $U^{a=b}$ in the commutative diagram
$$
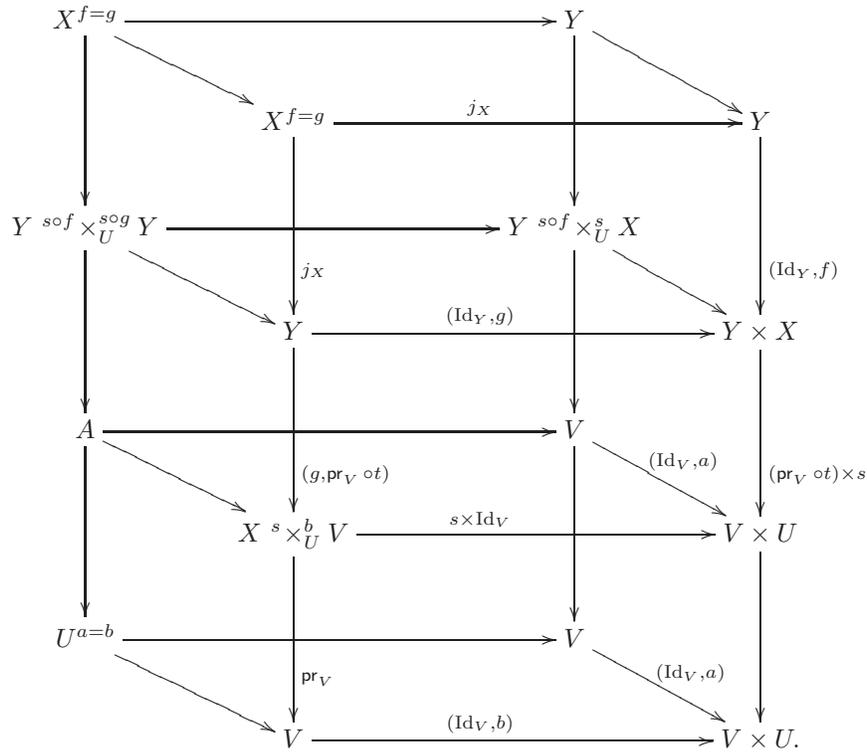
\xymatrix{
X^{f=g} \ar[dd] \ar[rrr] \ar[dr] &&& Y \ar[dd] \ar[dr]
\\
& X^{f=g} \ar[dd]^(.70){j_X} \ar[rrr]^(.40){j_X} &&& Y \ar[dd]^(.70){(\Id_Y,f)}
\\
Y ~^{s \circ f}\!\times_{U}^{s \circ g} Y  \ar[dd] \ar[dr] \ar[rrr] &&& Y ~^{s \circ f}\!\times_{U}^{s} X \ar[dd] \ar[dr]
\\
& Y \ar[dd]^(.70){(g, \pr_V \circ t)}  \ar[rrr]^(.40){(\Id_Y,g)} &&& Y \times X \ar[dd]^(.70){(\pr_V \circ t) \times s}
\\
A  \ar[dd] \ar[dr] \ar[rrr] &&& V \ar[dd] \ar[dr]^-{(\Id_V,a)}
\\
&  X ~^{s}\!\times_{U}^{b} V \ar[dd]^(.70){\pr_V} \ar[rrr]^(.40){s \times \Id_V} &&& V \times U \ar[dd]
\\
U^{a=b} \ar[rrr] \ar[dr] &&& V \ar[dr]^{(\Id_V,a)}
\\
& V \ar[rrr]^(.40){(\Id_V,b)} &&& V \times U.
}$$
It is only left to note that the above morphism can be rewritten precisely as in the statement of the proposition.
\end{proof}
\end{Prop}

\smallskip
\begin{Cor}\label{morph_of_tr_in_icoh}
Applying in the setting of the above proposition the functor 
$$\xymatrix{
\Corr(\Sch_{\aft})^{\proper} \ar[rr]^-{\ICoh} && 2\Cat_k
}$$
we see that given a commutative diagram
$$\xymatrix{
X \ar[r]^-{g_X} \ar[d]_-{f} & X \ar[d]^-{f}
\\
Y \ar[r]_-{g_Y} & Y 
}$$
the morphism of traces 
$$\xymatrix{
\Gamma(X^{g_X}, \omega^{\ICoh}_{X^{g_X}} ) \simeq \Tr_{2\Cat_k}\big{(}(g_X)_*\big{)} \ar[rr] && \Tr_{2\Cat_k}\big{(}(g_Y)_*\big{)} \simeq \Gamma(Y^{g_Y}, \omega^{\ICoh}_{Y^{g_Y}})
}$$
induced by the diagram
$$\xymatrix{
\ICoh(X) \ar@/_/[dd]_-{f_*} \ar[rr]^-{(g_X)_*} && \ICoh(X) \ar@/_/[dd]_-{f_*} \ar@2[ddll]
\\
\\
\ICoh(Y) \ar@/_/[uu]_-{f^!} \ar[rr]_-{(g_Y)_*} && \ICoh(Y) \ar@/_/[uu]_-{f^!}
}$$
is given by the counit of the adjunction $(f^g)_* \dashv (f^g)^!$ where $\xymatrix{X^{g_X} \ar[r]^-{f^g} & Y^{g_Y}}$ is the induced morphism between fixed points.
\end{Cor}

\subsection{Decorated correspondences and orientations}
Our goal in this section is to show the morphism of traces
$$\xymatrix{
\Gamma(X^g,\mathcal O_{X^g}) \simeq \Tr_{2\Cat_k}(g_*^{\QCoh}) \ar[rrr]^-{\Tr_{2\Cat_k}(- \otimes \mathcal O_X)}_-{\sim} &&& \Tr_{2\Cat_k}(g_*^{\ICoh}) \simeq \Gamma(X^g, \omega_{X^g})
}$$
induced by the diagram
$$\xymatrix{
\QCoh(X) \ar[d]_-{- \otimes \mathcal O_X} \ar[rr]^-{g_*} && \QCoh(X) \ar[d]^-{- \otimes \mathcal O_X}
\\
\ICoh(X) \ar[rr]_-{g_*} && \ICoh(X)
}$$
is induced by the canonical orientation (Construction \ref{FixedPointsAreCalabi}). Using the fact that for an eventually coconnective morphism  $\xymatrix{X \ar[r]^-f & Y}$ almost of finite type between Noetherian schemes one has a Grothendieck formula $f^!- \simeq \omega_f \otimes f^*-$ (\cite[Corollary 7.2.5.]{Gaits}) we will reduce the calculation of the morphism of traces to a simple calculation in a version of the category of correspondences.

\smallskip We start with the following
\begin{Def}\label{Dec_corr} We define an \emph{$(\infty,2)$-category of correspondences decorated by $\QCoh$} denoted by $\Corr(\Sch)^{\QCoh}$ as follows:
\begin{enumerate}
\item Its objects are those of $\Sch$.

\item Given $X,Y \in \Sch$ a morphism from $X$ to $Y$ in $\Corr(\Sch)^{\QCoh}$ is given by a span
$$\xymatrix{
X & W \ar[r]^-{f} \ar[l]_-{g} & Y
}$$
together with a quasi-coherent sheaf $\mcal{F}_W \in \QCoh(Y)$ on $W$ and the composition of morphisms is given by taking pullbacks of schemes and box products of sheaves.
\item Given two morphisms $(W_1,\mcal{F}_{W_1}), (W_2,\mcal{F}_{W_2}) \in \Hom_{\Corr(\Sch)^{\QCoh}}(X,Y)$ a $2$-morphism from $(W_1,\mcal{F}_{W_1})$ to $(W_2,\mcal{F}_{W_2})$ is given by a commutative diagram
$$\xymatrix{
& W_1 \ar[dr] \ar[dl] \ar[dd]^-{h}
\\
X && Y
\\
& W_2 \ar[ur] \ar[ul]
}$$
in $\Sch$ and a morphism $\xymatrix{h^* \mcal{F}_{W_2} \ar[r] & \mcal{F}_{W_1}}$ in $\QCoh(\mathcal{F}_{W_1})$.
\end{enumerate}
\end{Def}
\begin{Not} We will further denote by $\la \prescript{g}{X}{W}_Y^f, \mcal{F}_W \ra$ the morphism $\xymatrix{X & (W,\mcal{F}_W) \ar[l]_-{g} \ar[r]^-{f} & Y}$ in $\Corr(\Sch)^{\QCoh}$ with the correspondence given by $\xymatrix{X & W \ar[r]^-{f} \ar[l]_-{g} & Y}$ and a sheaf given by $\mathcal{F}_W \in \QCoh(W)$ and depict it as
$$\xymatrix{
& (W, \mathcal{F}_W) \ar[dl]_-{g} \ar[dr]^-{f}
\\
X && Y.
}$$
In the case we omit $\mathcal F_W$ from notation we assume that it is given by $\mathcal{O}_W$.
\end{Not}

\smallskip

Now note that the $(\infty,2)$-category $\Corr(\Sch)^{\QCoh}$ is symmetric monoidal with the monoidal structure given by the cartesian product of underlying objects of $\Sch$ (the morphisms are tensored by taking box product of the corresponding sheaves). Moreover, if $X \in \Sch$ is a scheme and $\mcal{M} \in \Pic(X)$ is a line bundle on $X$ then it is straightforward to see that the morphisms
$$\xymatrix{
& (X,\mathcal{M}) \ar[dr]^-{\Delta} \ar[dl]_-{p} &&&  (X,\mathcal{M}^{-1}) \ar[dl]_-{\Delta} \ar[dr]^-{p} 
\\
\ast && X \times X & X \times X && \ast
}$$
in $\Corr(\Sch)^{\QCoh}$ exhibit $X \in \Corr(\Sch)^{\QCoh}$ as a self-dual object. 

\medskip

\begin{Prop}\label{trace_in_dec_corr} Let $X \in \Sch$ be a scheme with an endomorphism $\xymatrix{X \ar[r]^-{g} & X}$ and $\mcal{M} \in \Pic(X)$ be a line bundle on $X$. Then there is an equivalence 
$$\xymatrix{
\Tr_{\Corr(\Sch)^{\QCoh}}^{\mathcal M}(g) \stackrel{\eta_{\mcal{M}}}{\simeq} \la {}_{\ast}(X^g)_{\ast} \ra
}$$
in $\Hom_{\Corr(\Sch)^{\QCoh}}(\ast,\ast)$, where $\Tr_{\Corr(\Sch)^{\QCoh}}^{\mcal{M}}(g)$ is the trace of  $\xymatrix{X \ar[r]^-{g} & X}$ in $\Corr(\Sch)^{\QCoh}$ with respect to the dualization data $\xymatrix{\ast & (X,\mathcal{M}) \ar[r]^-{\Delta} \ar[l]_-{p} & X \times X}$ and $\xymatrix{X \times X & (X,\mathcal{M}^{-1}) \ar[r]^-{p} \ar[l]_-{\Delta} & \ast}$ on $X \in \Corr(\Sch)^{\QCoh}$.
\begin{proof}
By definition the trace $\Tr^{\mcal{M}}_{\Corr(\Sch)^{\QCoh}}(g)$ is given by the composite
$$\xymatrix{
& (X,\mathcal{M}) \ar[dr]^-{\Delta} \ar[dl]_-{p} && X \times X \ar[dl]_-{\Id_{X} \times \Id_X} \ar[dr]^-{g \times \Id_X} &&  (X,\mathcal{M}^{-1}) \ar[dl]_-{\Delta} \ar[dr]^-{p} 
\\
\ast && X \times X && X \times X && \ast
}$$
which composing the first two morphisms can be rewritten as
$$\xymatrix{
& (X,\mcal{M}) \ar[dl]_-{p} \ar[dr]^-{(g,\Id_X)} && (X,\mcal{M}^{-1}) \ar[dr]^-{p} \ar[dl]_-{\Delta}
\\
\ast && X \times X && \ast .
}$$
The equivalence $\eta_{\mcal{M}}$ is now induced from the pullback diagram
$$\xymatrix{
X^g \ar[r]^-{i} \ar[d]_-{i} & X \ar[d]^-{\Delta}
\\
X \ar[r]_-{(\Id_X,g)} & X \times X
}$$
and the equivalence $i^* \mcal{M} \otimes i^* \mcal{M}^{-1} \simeq \mathcal{O}_{X^g}$.
\end{proof}
\end{Prop}

\medskip

We now prove the following
\begin{Prop}\label{morph_of_tr_in_dec_corr}
Let $X$ be a scheme with an endomorphism $\xymatrix{X \ar[r]^-{g} & X}$ and $\mathcal{M}_1,\mathcal{M}_2 \in \Pic(X)$ be two line bundles on $X$. Then the morphism of traces
$$\xymatrix{
\la {}_{\ast}(X^g)_{\ast} \ra \stackrel{\eta_{\mcal{M}_1}}{\simeq} \Tr^{\mcal{M}_1}_{\big{(}\Corr(\Sch)^{\QCoh}\big{)}^{2-\op}}(g) \ar[rr] && \Tr^{\mcal{M}_2}_{\big{(}\Corr(\Sch)^{\QCoh}\big{)}^{2-\op}}(g) \stackrel{\eta_{\mcal{M}_2}}{\simeq} \la {}_{\ast}(X^g)_{\ast} \ra 
}$$
in $\Hom_{\big{(}\Corr(\Sch)^{\QCoh}\big{)}^{2-\op}}(\ast,\ast)$ induced by the commutative diagram
$$\xymatrix{
\ar[d]_-{\Id_X} X \ar[r]^-{g} & X \ar[d]^-{\Id_X}
\\
X \ar[r]_-{g} & X
}$$
is given by the identity $2$-morphism
$$\xymatrix{
& X^g \ar[dd]^-{\Id_{X^g}} \ar[dl] \ar[dr]
\\
\ast && \ast
\\
& X^g \ar[ur] \ar[ul] .
}$$

\begin{proof} By \cite[Example 1.2.5]{We} the morphism of traces is induced by the diagram
$$\xymatrix{
\ast \ar[dd]_-{\Id_{\ast}} \ar[rrr]^-{\la \prescript{}{\ast}X^{\Delta}_{X \times X}, \mcal{M}_1 \ra} &&& X \times X \ar@2[ddlll]  \ar[rrrr]^-{\la \prescript{\Id_{X \times X}}{X \times X}{X \times X}^{g \times \Id_X}_{X \times X} \ra} \ar[dd] &&&& X \times X \ar[rrr]^-{\la \prescript{\Delta}{X \times X} X_{\ast}, \mcal{M}_1^{-1} \ra} \ar[dd] \ar@2[ddllll] &&& \ast \ar[dd]^-{\Id_{\ast}} \ar@2[ddlll] 
\\
\\
\ast \ar[rrr]_-{\la \prescript{}{\ast}X^{\Delta}_{X \times X}, \mcal{M}_2 \ra} &&& X \times X \ar[rrrr]_-{\la \prescript{\Id_{X \times X}}{X \times X}{X \times X}^{g \times \Id_X}_{X \times X} \ra} &&& & X \times X \ar[rrr]_-{\la \prescript{\Delta}{X \times X} X_{\ast}, \mcal{M}_2^{-1} \ra} &&& \ast
}$$
where the vertical morphisms are given by $\xymatrix{X \times X & (X \times X,\mcal{M}_1^{-1} \boxtimes \mcal{M}_2) \ar[r]^-{\Id_{X \times X}} \ar[l]_-{\Id_{X \times X}} & X \times X}$. The result now follows from equivalences
$$\xymatrix{
\la \prescript{\Delta}{X \times X} X_{\ast}, \mcal{M}_1^{-1} \ra \circ \la \prescript{\Id_{X \times X}}{X \times X}{X \times X}^{g \times \Id_X}_{X \times X} \ra \circ \la \prescript{}{\ast}X^{\Delta}_{X \times X}, \mcal{M}_1 \ra \simeq \la \prescript{}{\ast}(X^g)_{\ast}\ra
},$$
$$\xymatrix{
\la \prescript{\Delta}{X \times X} X_{\ast}, \mcal{M}_2^{-1} \ra \circ \la \prescript{\Id_{X \times X}}{X \times X}{X \times X}^{\Id_{X \times X}}_{X \times X}, \mcal{M}_1^{-1} \boxtimes \mcal{M}_2 \ra  \circ \la \prescript{\Id_{X \times X}}{X \times X}{X \times X}^{g \times \Id_X}_{X \times X} \ra \circ \la \prescript{}{\ast}X^{\Delta}_{X \times X}, \mcal{M}_1 \ra \simeq \la \prescript{}{\ast}(X^g)_{\ast}\ra
},$$
$$\xymatrix{
\la \prescript{\Delta}{X \times X} X_{\ast}, \mcal{M}_2^{-1} \ra \circ \la \prescript{\Id_{X \times X}}{X \times X}{X \times X}^{g \times \Id_X}_{X \times X} \ra \circ  \la \prescript{\Id_{X \times X}}{X \times X}{X \times X}^{\Id_{X \times X}}_{X \times X}, \mcal{M}_1^{-1} \boxtimes \mcal{M}_2 \ra \circ \la \prescript{}{\ast}X^{\Delta}_{X \times X}, \mcal{M}_1 \ra \simeq \la \prescript{}{\ast}(X^g)_{\ast}\ra
},$$
$$\xymatrix{
\la \prescript{\Delta}{X \times X} X_{\ast}, \mcal{M}_2^{-1} \ra \circ \la \prescript{\Id_{X \times X}}{X \times X}{X \times X}^{g \times \Id_X}_{X \times X} \ra \circ \la \prescript{}{\ast}X^{\Delta}_{X \times X}, \mcal{M}_2 \ra \simeq \la \prescript{}{\ast}(X^g)_{\ast}\ra
}$$
and the fact that the corresponding $2$-morphisms are given by identity maps $\Id_{X^g}$.
\end{proof}
\end{Prop}

To see why the proposition above is useful, we have the following generalization of \cite[Chapter 5, 5.3.1]{GaitsRozI}
\begin{Theor} The quasi-coherent sheaves functor can be lifted to a symmetric monoidal functor
$$\xymatrix{
\big{(}\Corr(\Sch)^{\QCoh}\big{)}^{2-\op} \ar[rr]^-{\widetilde{\QCoh}} && 2\Cat_k
}$$
where:
\begin{itemize}
\item An object $X \in \Corr(\Sch)^{\QCoh}$ is sent to the category of quasi-coherent sheaves $\QCoh(X) \in 2\Cat_k$ on it.
\item A morphism
$$\xymatrix{
& (W, \mathcal{F}_W) \ar[dl]_-{g} \ar[dr]^-{f}
\\
X && Y
}$$
in $\Corr(\Sch)^{\QCoh}$ is sent to the morphism
$$\xymatrix{
\QCoh(X) \ar[rrr]^-{f_*(\mathcal{F}_W \otimes g^*-)} &&& \QCoh(Y)
}$$
in $2\Cat_k$.
\item A $2$-morphism
$$\xymatrix{
& (W_1,\mathcal{F}_{W_1}) \ar[dr]^-{f} \ar[dl]_-{g} \ar[dd]^-{h}
\\
X && Y
\\
&  (W_2,\mathcal{F}_{W_2})  \ar[ur]_-{s} \ar[ul]^-{t}
}$$
with $\xymatrix{h^* \mcal{F}_{W_2} \ar[r]^-{\eta} & \mcal{F}_{W_1}}$ is sent to the $2$-morphism
$$\xymatrix{
s_*(\mathcal{F}_{W_2} \otimes t^*-)  \ar[r] & s_*(h_*h^*\mathcal{F}_{W_2} \otimes t^*-) \simeq  s_*h_*(h^*\mathcal{F}_{W_2} \otimes h^* t^*-) \simeq f_*(h^*\mathcal{F}_{W_2} \otimes g^*-) \ar[r]^-{\eta} & f_*(\mathcal{F}_{W_1} \otimes g^*-) 
}$$
in $2\Cat_k$.
\end{itemize}
\end{Theor}

\begin{Ex}\label{icoh_from_dec_corr} Let $X$ be a Gorenstein Noetherian scheme. Then the functor $\widetilde{\QCoh}$ sends the morphism
$$\xymatrix{
& (X,\omega_X) \ar[dr]^-{\Delta} \ar[dl]_-{p}
\\
\ast && X \times X
}$$
in $\big{(}\Corr(\Sch)^{\QCoh}\big{)}^{2-\op}$ to $\Delta_* \omega_X \in \Hom_{2\Cat_k}(\Vect_k,\QCoh(X \times X)) \simeq \QCoh(X \times X)$ and the morphism
$$\xymatrix{
& (X,\omega_X^{-1}) \ar[dl]_-{\Delta} \ar[dr]^-{p} 
\\
X \times X && \ast
}$$
to $\Gamma(X, \omega_X^{-1} \otimes \Delta^*-) \simeq \Gamma(X, \omega_{\Delta} \otimes \Delta^*-) \simeq \Gamma(X, \Delta^!-) \in \Hom_{2\Cat_k}(\QCoh(X \times X),\Vect_k)$.
\end{Ex}

\smallskip

We can now prove
\begin{Prop}\label{TraceIsFromCalabiYau_proof}
For a classical smooth scheme $X$ together with an endomorphism $\xymatrix{X \ar[r]^-{g} & X}$ the morphism of traces
$$\xymatrix{
\Gamma(X^g,\mcal{O}_{X^g}) \stackrel{\alpha_{\QCoh}}{\simeq} \Tr_{2\Cat_k}(g_*^{\QCoh}) \ar[rrr]^-{\Tr_{2\Cat_k}(- \otimes \mathcal O_X)}_-{\sim} &&& \Tr_{2\Cat_k}(g_*^{\ICoh})  \stackrel{\alpha_{\ICoh}^{-1}}{\simeq} \Gamma(X^g, \omega_{X^g})
}$$
induced by the diagram
$$\xymatrix{
\QCoh(X) \ar[d]_-{- \otimes \mathcal O_X} \ar[rr]^-{g_*} && \QCoh(X) \ar[d]^-{- \otimes \mathcal O_X}
\\
\ICoh(X) \ar[rr]_-{g_*} && \ICoh(X),
}$$
can be obtained by applying the global sections functor $\Gamma(X^g, -)$ to the canonical orientation $\xymatrix{\mathcal O_{X^g} \ar[r]^-{\orient_C}_-{\sim} & \omega_{X^g}}$ (see Construction \ref{FixedPointsAreCalabi}) on $X^g$, where the equivalences $\alpha_{\QCoh}$ and $\alpha_{\ICoh}$ above are given by Corollary \ref{TraceOfQCoh} and Corollary \ref{TraceOfICoh} respectively.

\begin{proof} Due to the equivalence $\ICoh(X) \simeq \QCoh(X)$ as $X$ is smooth and classical (see Example \ref{icoh_of_smooth}) and Example \ref{icoh_from_dec_corr} above we note that the morphism of traces we are interested in can be obtained by applying the functor $\xymatrix{\big{(}\Corr(\Sch)^{\QCoh}\big{)}^{2-\op} \ar[r]^-{\widetilde{\QCoh}} & 2\Cat_k}$ to the morphism of traces
$$\xymatrix{
\Tr^{\mcal{O}_X}_{\big{(}\Corr(\Sch)^{\QCoh}\big{)}^{2-\op}}(g) \ar[r] & \Tr^{\omega_X}_{\big{(}\Corr(\Sch)^{\QCoh}\big{)}^{2-\op}}(g) 
}$$
in $\big{(}\Corr(\Sch)^{\QCoh}\big{)}^{2-\op}$. Now using equivalences $\eta_{\mcal{O}_X}$ and $\eta_{\omega_X}$ we can form a commutative diagram
$$\resizebox{!}{0.027\textwidth}{
\xymatrix{
\Gamma(X^g,\mcal{O}_{X^g}) \ar[dr]  \ar[r]^-{\alpha_{\QCoh}}_-{\sim} &  \widetilde{\QCoh} \bigg{(} \Tr^{\mcal{O}_X}_{\big{(}\Corr(\Sch)^{\QCoh}\big{)}^{2-\op}}(g) \bigg{)} \ar[d]^-{\widetilde{\QCoh}(\eta_{\omega_X})}_-{\sim} \ar[rrr]^-{\Tr_{2\Cat_k}(- \otimes \mathcal O_X)}_-{\sim} &&& \widetilde{\QCoh} \bigg{(} \Tr^{\omega_X}_{\big{(}\Corr(\Sch)^{\QCoh}\big{)}^{2-\op}}(g) \bigg{)} \ar[r]^-{\alpha_{\ICoh}^{-1}}_-{\sim} & \Gamma(X^g, \omega_{X^g})
\\
 & \widetilde{\QCoh}(\la \prescript{}{\ast}(X^g)_{\ast}\ra) \ar[rrr] &&& \widetilde{\QCoh}(\la \prescript{}{\ast}(X^g)_{\ast}\ra) \ar[u]^-{\widetilde{\QCoh}(\eta_{\omega_X}^{-1})}_-{\sim} \ar[ur] .
}
}$$
Since the left diagonal morphism is identity by the construction and the bottom horizontal morphism is identity by Proposition \ref{morph_of_tr_in_dec_corr}, we see that the morphism of traces $\Tr_{2\Cat_k}(- \otimes \mathcal O_X)$ is given by the right diagonal morphism, that is, unwinding the definition of the morphism $\eta_{\mcal{\omega_X}}$ from Proposition \ref{trace_in_dec_corr} by the composite of
$$
 \Gamma(X^g,\mathcal{O}_{X^g}) \simeq  \Gamma(X^g,i^* \omega_X^{-1} \otimes i^* \omega_X) \simeq \Gamma(X,\omega_X^{-1} \otimes i_* i^* \omega_X) \simeq  \Gamma(X, \omega_X^{-1} \otimes \Delta^* (g \times \Id_X)_*\Delta_*\omega_X)
$$
and
$$
\Gamma(X, \omega_X^{-1} \otimes \Delta^* (g \times \Id_X)_*\Delta_*\omega_X) \simeq \Gamma(X,\Delta^! (\Id_X, g)_*  \omega_X) \simeq \Gamma(X,i_* i^!  \omega_X) \simeq \Gamma(X,i_* \omega_{X^g}) \simeq \Gamma(X^g, \omega_{X^g}).
$$
The result now follows from observation that the morphisms
$$
i_* \mcal{O}_{X^g} \simeq i_*(i^* \omega_X^{-1} \otimes i^* \omega_X) \simeq \omega_X^{-1} \otimes i_* i^* \omega_X \simeq \omega_X^{-1} \otimes \Delta^* (g,\Id_X)_* \omega_X \simeq \Delta^! (\Id_X, g)_* \omega_X \simeq  i_* i^! \omega_X \simeq i_* \omega_{X^g}
$$
and
$$
i_* i^* \mcal{O}_{X^g} \simeq i_*(i^* \omega_X^{-1} \otimes i^* \omega_X) \simeq i_*(i^* \omega_{X/X \times X} \otimes i^* \omega_X) \simeq  i_*(\omega_{X^g/X} \otimes i^* \omega_X) \simeq i_* i^! \omega_X  \simeq i_* \omega_{X^g}
$$
coincide and Construction \ref{FixedPointsAreCalabi}.
\end{proof} 
\end{Prop}

\begin{Rem} Using the fact that the functor $\xymatrix{\big{(}\Corr(\Sch)^{\QCoh}\big{)}^{2-\op} \ar[r]^-{\QCoh} & 2\Cat_k}$ sends a morphism
$$\xymatrix{
& (X,E) \ar[dr]^-{\Id_X} \ar[dl]_-{p}
\\
\ast && X 
}$$
in $\big{(}\Corr(\Sch)^{\QCoh}\big{)}^{2-\op}$ to the morphism $\xymatrix{\Vect_k \ar[r]^-{E} & \QCoh(X)}$ in $2\Cat_k$ one can also obtain a proof of \cite[Proposition 2.2.3.]{We} by calculating appropriate trace in the category of decorated correspondences and then mapping it to $2\Cat_k$.
\end{Rem}

\bigskip\bigskip

\noindent
Grigory~Kondyrev, {\sc Northwestern University; National Research University Higher School of Economics,}
\href{mailto:gkond@math.northwestern.edu}{gkond@math.northwestern.edu}

\smallskip

\noindent
Artem~Prikhodko, {\sc National Research University Higher School of Economics, Russian Federation; Center for Advanced Studies, Skoltech, Russian Federation,}
\href{mailto:artem.n.prihodko@gmail.com}{artem.n.prihodko@gmail.com}

\end{document}